\documentclass[]{article}
\usepackage[utf8]{inputenc}
\usepackage{amsmath}
\usepackage{amssymb}
\usepackage{amsthm}
\usepackage{fullpage}
\usepackage{mathrsfs}
\usepackage{tikz}
\usepackage{xcolor}
\usepackage{hyperref}
\usepackage{url}
\usepackage{braket}
\usepackage{nccmath}
\usepackage{stmaryrd}
\usepackage{xfrac}
\usepackage{authblk}
\usepackage[capitalize]{cleveref}
\usepackage{bbm}
\usepackage{makecell} 

\usepackage{algorithm}
\usepackage[noend]{algpseudocode}

\newtheorem{innercustomgeneric}{\customgenericname}
\providecommand{\customgenericname}{}
\newcommand{\newcustomtheorem}[2]{%
	\newenvironment{#1}[1]
	{%
		\renewcommand\customgenericname{#2}%
		\renewcommand\theinnercustomgeneric{##1}%
		\innercustomgeneric
	}
	{\endinnercustomgeneric}
}

\hypersetup{
	colorlinks=true,       
	linkcolor=blue,          
	citecolor=blue,        
	filecolor=blue,      
	urlcolor=blue           
}

\numberwithin{equation}{section}
\numberwithin{figure}{section}

\allowdisplaybreaks
\theoremstyle{plain}
\newtheorem{lemma}{Lemma}[section]
\newtheorem{proposition}[lemma]{Proposition}
\newtheorem{corollary}[lemma]{Corollary}

\newtheorem{theorem}[lemma]{Theorem}

\theoremstyle{definition}
\newtheorem{definition}[lemma]{Definition}
\newtheorem{remark}[lemma]{Remark}
\newtheorem{example}[lemma]{Example}
\newcustomtheorem{customthm}{Theorem}
\newcustomtheorem{customconj}{Conjecture}

\DeclareSymbolFont{bbold}{U}{bbold}{m}{n}
\DeclareSymbolFontAlphabet{\mathbbold}{bbold}

\begin{document}
	
	\title{Consistency and inconsistency in $k$-means clustering}

    \author[1]{Mo\"ise Blanchard}
    \author[2]{Adam Quinn Jaffe}
    \author[3]{Nikita Zhivotovskiy}
    
    \affil[1]{School of Industrial and Systems Engineering, Georgia Tech, Atlanta, GA}
    \affil[2]{Department of Statistics, Columbia University, New York, NY}
    \affil[3]{Department of Statistics, UC Berkeley, Berkeley, CA}

	\date{\today}
    
	\maketitle

        \begin{abstract}
            A celebrated result of Pollard proves asymptotic consistency for $k$-means clustering when the population distribution has finite variance.
            In this work, we point out that the population-level $k$-means clustering problem is, in fact, well-posed under the weaker assumption of a finite expectation, and we investigate whether some form of asymptotic consistency holds in this setting.
            As we illustrate in a variety of negative results, the complete story is quite subtle; for example, the empirical $k$-means cluster centers may fail to converge even if there exists a unique set of population $k$-means cluster centers.
            A detailed analysis of our negative results reveals that inconsistency arises because of an extreme form of cluster imbalance, whereby the presence of outlying samples leads to some empirical $k$-means clusters possessing very few points.
            We then give a collection of positive results which show that some forms of asymptotic consistency, under only the assumption of finite expectation, may be recovered by imposing some a priori degree of balance among the empirical $k$-means clusters.
        \end{abstract}

    \setcounter{tocdepth}{1}
	\tableofcontents

\section{Introduction}
\label{sec1}

A celebrated result of Pollard \cite{Pollard} states that the empirical
$k$-means cluster centers of independent, identically-distributed (i.i.d.)
samples of a random variable $X$ converge almost surely to the population
$k$-means cluster centers under the assumption of finite variance
$\mathbb{E}\|X\|^{2}<\infty $. At the same time, the strong law of large
numbers (SLLN) due to Kolmogorov states that the empirical mean of i.i.d.
samples converges almost surely to the population mean under the minimal
assumption of finite expectation $\mathbb{E}\|X\| < \infty $, and it is
easy to see that empirical- and population-level $k$-means clustering reduce
to the empirical mean and the population mean, respectively, when
$k=1$. The goal of this paper is to investigate this conspicuous gap in
the required moment conditions. In other words, is $k$-means clustering
consistent for all $k\ge 1$ under the assumption
$\mathbb{E}\|X\|<\infty $?

\begin{figure}
\includegraphics[width=1.0\linewidth]{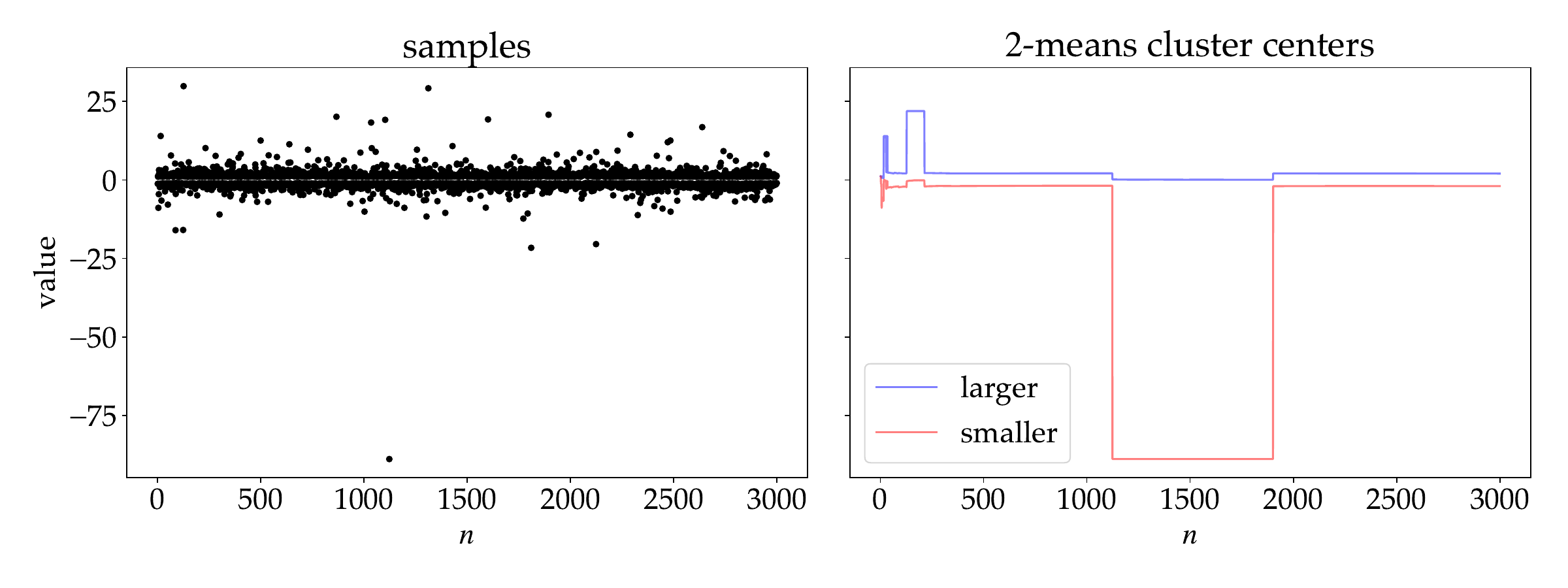}
\caption{Inconsistency for $k$-means clustering applied to a simulated data
set, for $k=2$. The samples come from a distribution $X$ which is symmetric
and satisfies $\mathbb{P}(|X|\ge t) = t^{-2}$ for all $t \ge 1$ (left).
The empirical $k$-means cluster centers do not converge (right).}
\label{fig:two_sided_pareto}
\end{figure}

Before getting into the details, let us consider the simplest possible
version of this problem. That is, we consider $k$-means clustering in
$\mathbb{R}^{d}$ when $k=2$ and $d=1$, and we suppose that $X$ is symmetric
and satisfies $\mathbb{P}(|X|\ge t) = t^{-2}$ for all $t \ge 1$. (In other
words, $|X|$ has a \textit{Pareto distribution with exponent 2}, denoted
$\textnormal{Par}(2)$.) Importantly, note that $X$ satisfies
$\mathbb{E}|X|^{2} = \infty $ but $\mathbb{E}|X|^{p}<\infty $ for all
$1\le p < 2$, so it just barely fails to have finite variance. In Figure~\ref{fig:two_sided_pareto}
we show a simulation of this data and the optimal cluster centers for
$n=1,\ldots , 3000$ samples, and we observe some large jumps of the empirical
cluster centers. As we rigorously prove later, $X$ admits a unique set
of population $k$-means cluster centers given by $\{-2,2\}$ (Example~\ref{ex:2-sided-pareto-pop})
but the empirical $k$-means cluster centers are infinitely often of the
form $\{-\sqrt{n},0\}$ or $\{0,\sqrt{n}\}$ up to logarithmic factors (Example~\ref{ex:inconsistency_two_sided}).
This example shows that $k$-means clustering can be inconsistent even when
all population-level objects are well-defined. The core of this paper is
a collection of examples in this spirit, which yield negative answers to
various forms of the motivating question above.

Now we introduce some notation. Recall that for i.i.d. samples
$X_{1},\ldots , X_{n}$ in $\mathbb{R}^{d}$ from the same distribution as
$X$, the \textit{empirical $k$-means clustering problem} is
%
\begin{equation}
\label{eqn:def-emp-k-means}
\begin{cases}
\textnormal{minimize} & \frac{1}{n}\sum _{i=1}^{n}\min _{c\in
\mathcal{C}}\|c-X_{i}\|^{2}
\\
\textnormal{over} & \mathcal{C}\subseteq \mathbb{R}^{d}
\\
\textnormal{with} & 1 \le \#\mathcal{C}\le k,
\end{cases}
%
\end{equation}
any solution of which is called a
\textit{set of empirical $k$-means cluster centers}. When
$\mathbb{E}\|X\|^{2}<\infty $, it is natural to expect that
\eqref{eqn:def-emp-k-means} should converge (in some suitable sense) to
the \textit{population $k$-means clustering problem} defined as
%
\begin{equation}
\label{eqn:def-pop-k-means}
\begin{cases}
\textnormal{minimize} & \mathbb{E}\left [\min _{c\in \mathcal{C}}\|c-X
\|^{2}\right ]
\\
\textnormal{over} & \mathcal{C}\subseteq \mathbb{R}^{d}
\\
\textnormal{with} & 1 \le \#\mathcal{C}\le k,
\end{cases}
%
\end{equation}
any solution of which is called a
\textit{set of population $k$-means cluster centers}. The value of
\eqref{eqn:def-pop-k-means} at $\mathcal{C}$ is called its
\textit{distortion}, denoted $D_{X}(\mathcal{C})$.

\subsection{Problem statement}
\label{sec1.1}

While it may seem that $\mathbb{E}\|X\|^{2}<\infty $ is necessary in order
to define the problem~\eqref{eqn:def-pop-k-means}, it turns out that there
is a canonical extension to the case where we have
$\mathbb{E}\|X\|^{2}=\infty $ but $\mathbb{E}\|X\|<\infty $. That is, we
may consider the ``renormalized'' problem
\begin{equation}
\label{eqn:def-pop-k-means-minimal}
\begin{cases}
\textnormal{minimize} & \mathbb{E}\left [\min _{c\in \mathcal{C}}\|c-X
\|^{2}-\|X\|^{2}\right ]
\\
\textnormal{over} & \mathcal{C}\subseteq \mathbb{R}^{d}
\\
\textnormal{with} & 1 \le \#\mathcal{C}\le k,
\end{cases}
\tag{$\ast $}
\end{equation}
and we say that any of its solutions is a set of population $k$-means cluster
centers. Also, the value of \eqref{eqn:def-pop-k-means-minimal} at
$\mathcal{C}$ is called its \textit{excess distortion} relative to the set
$\{0\}$, denoted $D_{X}(\mathcal{C}\,|\,\{0\})$. (We will later see that
one can similarly define excess distortion relative to other sets.) Two
remarks help to explain why this definition is indeed canonical. First,
note that problems~\eqref{eqn:def-pop-k-means} and
\eqref{eqn:def-pop-k-means-minimal} are equivalent whenever
$\mathbb{E}\|X\|^{2}<\infty $, and that
\eqref{eqn:def-pop-k-means-minimal} is well-defined whenever
$\mathbb{E}\|X\|<\infty $. Second, note that problem~\eqref{eqn:def-pop-k-means-minimal}
captures, in the $k=1$ case, the notion that
$\mathbb{E}\|X\|<\infty $ suffices to define $\mathbb{E}[X]$, whereas problem~\eqref{eqn:def-pop-k-means}
does not. (Note that such renormalization is not required for the empirical
$k$-means clustering problem \eqref{eqn:def-emp-k-means}.) We note that
reductions from distortion to excess distortion have already appeared in
the finite-sample analyses of the papers
\cite{BiauDevroyeLugosi,Zhivotovskiy}, but neither of these works covers
our regime of interest, namely $\mathbb{E}\|X\|^{2}=\infty $.

The excess distortion $D_{X}(\mathcal{C}\,|\,\{0\})$ being well-defined
does not necessarily guarantee that the problem~\eqref{eqn:def-pop-k-means-minimal}
admits a solution, and we indeed show that solutions may not exist. Even
more delicately, we show that a wide range of population-level behaviors
can occur (Section~\ref{sec:cases}). That is, we show that there exist
distributions of a random variable $X$ such that
\begin{itemize}
\item[(i)] problem \eqref{eqn:def-pop-k-means-minimal} has a solution hence
finite optimal excess distortion (Example~\ref{ex:2-sided-pareto-pop}),
\item[(ii)] problem~\eqref{eqn:def-pop-k-means-minimal} has a finite optimal
excess distortion but no solution (Example~\ref{ex:1-sided-pareto-pop}),
and
\item[(iii)] problem~\eqref{eqn:def-pop-k-means-minimal} has infinite optimal
excess distortion hence no solution (Lemma~\ref{lemma:charact_finite_obj}).
\end{itemize}
We emphasize that the landscape of these cases is only visible when
$k\ge 2$ and $\mathbb{E}\|X\|^{2} = \infty $; if $k= 1$ or
$\mathbb{E}\|X\|^{2} < \infty $, then Kolmogorov's SLLN and Pollard's result
in \cite{Pollard} imply that case (i) always obtains and that strong convergence
of the empirical cluster centers to the population cluster centers always
holds.

As we will illustrate in a variety of negative results below, various notions
of consistency for the empirical $k$-means problem may fail even when the
requisite population-level objects are well-defined. For an example in
case (i), we have already seen Figure~\ref{fig:two_sided_pareto} which
shows that convergence of the empirical $k$-means cluster centers may fail
even when the population $k$-means cluster centers are well-defined. A
careful analysis of this example and related examples reveals the exact
mechanism by which this inconsistency occurs. Namely, for heavy-tailed
distributions $X$, there is a non-negligible probability that the largest
value among $X_{1},\ldots , X_{n}$ is so much larger than the other values
that it must be placed into a cluster alone. This can be observed directly
in Figure~\ref{fig:two_sided_pareto}, where jumps in either cluster center
are due to extreme values of samples.

Because inconsistency in $k$-means clustering occurs due to cluster imbalance,
it is natural to consider an alternative procedure in which we impose some
a priori constraint on the balance of the empirical clusters. To precisely
introduce such a procedure, recall \cite[Chapter~9.1]{Bishop} that problem
\eqref{eqn:def-emp-k-means} is an optimization over
\textit{cluster centers} but it may be equivalently stated as an optimization
over \textit{clusters} themselves; mathematically, we regard the set of
clusters as the \textit{Voronoi partition} determined by the cluster centers,
i.e., the partition of $\mathbb{R}^{d}$ into finitely-many regions determined
by which cluster center is closest. Then, for some integer
$0\le \gamma \le n$, we consider the empirical clustering problem
%
\begin{equation}
\label{eqn:def-emp-k-means-balance}
\begin{cases}
\textnormal{minimize} & \frac{1}{n}\sum _{\ell =1}^{k}\sum _{X_{i}
\in \mathcal{V}_{\ell}}\|\bar m_{n}(\mathcal{V}_{\ell})-X_{i}\|^{2}
\\
\textnormal{over} & \textnormal{Voronoi partitions } \mathcal{V}_{1},
\ldots , \mathcal{V}_{k} \textnormal{ of } \mathbb{R}^{d}
\\
\textnormal{with} & \#\{1\le i \le n: X_{i}\in \mathcal{V}_{\ell}\}
\ge \gamma \textnormal{ for all } 1 \le \ell \le k,
\end{cases}
%
\end{equation}
where
$\bar m_{n}(\mathcal{V}):=\sum _{i=1}^{n}X_{i}\mathbbold{1}\{X_{i}
\in \mathcal{V}\}/\sum _{i=1}^{n}\mathbbold{1}\{X_{i}\in \mathcal{V}
\}$ is the conditional sample mean on a measurable region
$\mathcal{V}\subseteq \mathbb{R}^{d}$. Then, we are interested in the convergence
of the set of points
$\{\bar m_{n}(\mathcal{V}_{1}),\ldots , \bar m_{n}(\mathcal{V}_{k})\}$,
in a suitable sense, when $\gamma $ grows with $n$ at some carefully-chosen
rate. (There are some nuances to this choice of balanced clustering problem,
for example that $\{\mathcal{V}_{1},\ldots , \mathcal{V}_{k}\}$ is typically
not equivalent to the Voronoi partition of the points
$\{\bar m_{n}(\mathcal{V}_{1}),\ldots , \bar m_{n}(\mathcal{V}_{k})\}$
unless $\gamma = 0$. See Subsection~\ref{subsec:balance} for further
detail.)

\subsection{Main results}
\label{sec1.2}

Now we may state our positive and negative results concerning the various
desired notions of convergence in various cases of interest.

We consider convergence of the empirical cluster centers in the Hausdorff
metric (Section~\ref{sec:centers}). First, we show that the empirical cluster
centers may fail to converge even when the population cluster centers exist
and are unique (Example~\ref{ex:inconsistency_two_sided}). Second, we show
that imposing a linear balance constraint, that is
$\gamma _{n} = \alpha n$ for some suitable $0<\alpha <1$ depending on the
distribution of $X$, suffices for the convergence of the balanced empirical
$k$-means cluster centers to the (unconstrained) population $k$-means cluster
centers (Theorem~\ref{thm:cluster-consistency}). The choice of
$\alpha $ may be viewed as a form of well-specification, and it is typically
unknown in practice.

Similarly, we have positive and negative results concerning the convergence
of the excess distortion of the empirical cluster centers (Section~\ref{sec:distortion}).
First, we show that the excess distortion of the empirical cluster centers
may fail to converge to the excess distortion of the population cluster
centers (Example~\ref{ex:inconsistency_two_sided_imbalanced}). Second,
we show that imposing a polylogarithmic balance constraint, that is
$\gamma _{n}\ge (\log n)^{4}$, suffices for the convergence of the excess
distortion of the balanced empirical $k$-means cluster centers to the optimal
(unconstrained) excess distortion (Theorem~\ref{thm:consistency_log_points_per_cluster}).

Let us now summarize our main results thus far. For cases (i) and (ii),
we give Table~\ref{tab:summary} which summarizes the form of consistency
that is desired, the counterexample showing that such consistency can fail,
and the balance constraint that leads to an analogous positive result.
For case (iii) we do not know of a meaningful form of consistency that
may be desired, since essentially no interesting population-level object
is well-defined; we comment further on this in Appendix~\ref{sec:case-iii}.
We also show that there is an intermediate form of consistency in case
(i) when polylogarithmic balance is imposed (Theorem~\ref{thm:consistency_unified}),
whereby we have convergence of the subset of non-diverging cluster centers,
and the limiting set of points is an optimal set of $k'$-means cluster
centers, for some $1\le k'\le k$.

\begin{table}[t]
\centering
\begin{tabular}{|c|c|c|}
		\hline
		& (i)\makecell[c]{\textbf{optimal cluster}\\\textbf{centers exist,}\\ \textbf{so optimal excess}\\ \textbf{distortion is finite}} & (ii)\makecell[c]{\textbf{optimal cluster}\\\textbf{centers do not exist,}\\ \textbf{but optimal excess}\\ \textbf{distortion is finite}} \\
		\hline
		\textbf{desired consistency} & \makecell[c]{convergence of empirical clusters\\centers in Hausdorff metric} & \makecell[c]{convergence of excess distortion\\ of empirical cluster centers} \\
		\hline
		\textbf{negative result} & \makecell[c]{symmetric two-sided $\textnormal{Par}(2)$\\distribution (Example~\ref{ex:inconsistency_two_sided})} & \makecell[c]{asymmetric two-sided $\textnormal{Par}(2)$\\distribution (Example~\ref{ex:inconsistency_two_sided_imbalanced})} \\
		\hline
		\makecell[c]{\textbf{positive result via}\\\textbf{balance constraint}} & \makecell[c]{$\gamma _{n} \ge \alpha n$ for $0<\alpha <1$ which is\\well-specified (Theorem~\ref{thm:cluster-consistency})} & $\gamma _{n} \ge (\log n)^{4}$ (Theorem~\ref{thm:consistency_log_points_per_cluster}) \\
		\hline
	\end{tabular}
\caption{Summary of positive and negative results for cases (i) and (ii),
i.e., the relevant population-level objects are well-defined. We do not
include results for case (iii) since it is not clear what a desirable notion
of consistency should be}
\label{tab:summary}
\end{table}

Lastly, we comment on some methodological implications of the theory developed
in this work. On the one hand, we recommend broader adoption of balanced
$k$-means clustering schemes, since their computation can be reduced to
efficient subroutines from computational optimal transport
\cite{Ng, CuturiDoucet, GenevayDulacArnoldVert} and hence is not too cumbersome.
On the other hand, our results have the following implications for the
standard (i.e., unconstrained) empirical $k$-means clustering problem,
simply by considering whether each balance constraint is binding or not.
Here, when we say an empirical cluster is ``meaningful'' we mean that its
center is close to some population cluster center.
\begin{itemize}
\item If all clusters are roughly balanced, then every cluster is meaningful.
\item If some clusters are moderately unbalanced, then, after discarding
them, the remaining clusters are meaningful.
\item If at least one cluster is extremely unbalanced, then all clusters
can fail to be meaningful.
\end{itemize}
An interesting avenue for future work would be to further understand the
use of these recommendations in practice.

\subsection{Related literature}
\label{sec1.3}

In order to situate our results in the context of related literature, we
briefly summarize some relevant literature on $k$-means clustering and
its variants.

In addition to the seminal work by Pollard on strong consistency
\cite{Pollard}, the statistical theory of $k$-means clustering has seen
many recent developments. For example, strong consistency is known for
$k$-means clustering in a variety of geometric settings
\cite{Parna1, Parna2, Lember, Thorpe, Laloe, JaffeClustering} and under
a variety of adaptive modifications
\cite{JiangAriasCastro, JaffeClustering}, but all results rely on some
condition at least as strong as finite variance. For another example, Pollard
also proved a central limit theorem for the empirical $k$-means cluster
centers \cite{PollardCLT}, and this also relies on the condition of finite
variance. There are also various finite-sample rates of convergence known
for $k$-means clustering (typically for the distortion rather than the
cluster centers); these all rely on stronger conditions than finite variance
\cite{BartlettLinderLugosi, BiauDevroyeLugosi, Levrard, Fefferman, appert2021new}
with the exception of \cite{Zhivotovskiy}.

Somewhat related to the present work is existing literature concerning
robustness properties of $k$-means clustering and its variants. For example,
in \cite[Section~2]{TrimmedKMeansII} it is shown that $k$-means clustering
lacks all of the usual notions of robustness from classical robust statistics.
A primary approach to overcome these difficulties is to employ
\textit{$\alpha $-trimmed $k$-means} \cite{TrimmedKMeansI} for some
$0<\alpha <1$, which amounts to minimizing the $k$-means clustering objective
while omitting up to proportion $\alpha $ of the data points; much is known
about trimmed $k$-means clustering, although it is still an active area
of research
\cite{RobustTrimmed, SoftTrimming, RobustBregman, WassersteinTrimmed, TrimmedCLT}.
Other robust variants of $k$-means clustering have been proposed, for example
model-based methods. \cite{GallegosRitter}; see
\cite{RobustClusteringSurvey} for a survey of statistical and computational
aspects of robust clustering.

One of our fundamental observations is that cluster imbalance leads to
statistical problems, and this has indeed been observed in previous literature.
For example, the works
\cite{Ng, CuturiDoucet, GenevayDulacArnoldVert} argue that cluster balance
is desirable from a practical point of view, and they discuss computational
aspects of implementing such balanced clustering procedures. For another
example, the recent work \cite{Zhivotovskiy} identifies (see
\cite[Theorem~3.1 and Theorem 3.10]{Zhivotovskiy}) that it is precisely
cluster imbalance that controls the minimax rate of convergence of the
distortion, when finite variance is assumed. Our results contribute to
this story by showing that extreme cluster imbalance can even lead to inconsistency.

Our various positive and negative results suggest an interesting open problem
of identifying the sharpest possible balance constraint
$\{\gamma _{n}\}_{n\in \mathbb{N}}$ leading to asymptotic consistency in
balanced $k$-means clustering. In case (i) we have seen that
$\gamma _{n} = \alpha n$ is sufficient for certain $0<\alpha <1$, and in
case (ii) we have seen that $\gamma _{n}=(\log n)^{4}$ is sufficient. However,
it is not clear whether either of these rates is sharp. It would be interesting
to try to identify the sharpest possible balance constraint for each case,
and also to address the analogous questions for convergence in probability
rather than convergence almost surely.

\section{Preliminaries}
\label{sec:prelim}

In this section we establish our notation and give some preliminaries that
will be used in the main results.

\paragraph{Notation}
Throughout the paper, we fix $d\in \mathbb{N}$ and for
$k\in \mathbb{N}$ we write
$\mathscr{C}_{k}:=\{\mathcal{C}\subseteq \mathbb{R}^{d}:1\le \#
\mathcal{C}\le k\}$ for the collection of all non-empty sets of at most
$k$ points in $\mathbb{R}^{d}$. We write
$(\Omega ,\mathscr{F},\mathbb{P})$ for a probability space, with expectation
$\mathbb{E}$, supporting an i.i.d. sequence of random variables
$X,X_{1},X_{2},\ldots $ taking values in $\mathbb{R}^{d}$. For
$p\ge 1$, we write $W_{p}$ for the $p$-Wasserstein distance between probability
measures on $\mathbb{R}^{d}$; we refer the reader to \cite{Villani} for
further detail, but we emphasize (see \cite[Definition~6.8]{Villani}) that
convergence in $W_{p}$ is equivalent to weak convergence plus convergence
of the $p$th moment. We write $D_{\textnormal{KL}}$ for the Kullback--Leibler
divergence (i.e., relative entropy) between two probability measures, and
by a slight abuse of notation we also write
$D_{\textnormal{KL}}(x\,\|\,y):=x \log (x/y) + (1-x)\log ((1-x)/(1-y))$
for the relative entropy between Bernoulli distributions with parameters
$0\le x,y\le 1$. We use standard notation for asymptotic comparisons of
sequences, e.g. $a_{n} = O(b_{n})$ to mean
$\limsup _{n\to \infty}a_{n}/b_{n} < \infty $ and $a_{n}\ll b_{n}$ to mean
$\lim _{n\to \infty}a_{n}/b_{n} = 0$.

\subsection{Clustering definitions}
\label{subsec:basic}

We begin with some precise notions that will allow us to describe the population-
and empirical-level $k$-means clustering problems under the minimal moment
assumption of $\mathbb{E}\|X\|<\infty $. For instance, we have the following
elementary observation which follows from some casework; the proof is therefore
omitted for the sake of brevity:

\begin{lemma}
\label{lem2.1}
If $X$ satisfies $\mathbb{E}\|X\|<\infty $, then we have
\begin{equation*}
\mathbb{E}\left |\min _{c\in \mathcal{C}}\|c-X\|^{2}-\min _{c'\in
\mathcal{C}'}\|c'-X\|^{2}\right |<\infty
\end{equation*}
for any $\mathcal{C},\mathcal{C}'\in \mathscr{C}_{k}$.
\end{lemma}

The preceding integrability trick is known in probability theory as
\textit{renormalization}; although the tails of the random variable
$\min _{c\in \mathcal{C}}\|c-X\|^{2}$ are too heavy to be integrable under
the condition $\mathbb{E}\|X\|<\infty $, the random variable given by the
difference
$\min _{c\in \mathcal{C}}\|c-X\|^{2} - \min _{c'\in \mathcal{C}'}\|c'-X
\|^{2}$ experiences a great deal of cancellation and hence its tails are
light enough to be integrable under the condition
$\mathbb{E}\|X\|<\infty $. See \cite{SchoetzSLLN, FrechetMeanInfDim} for
further explanation of the renormalization trick in the context of Fr\'echet
means, which corresponds to $k$-means clustering with $k=1$ in a general
metric space.

The preceding result allows us to make sense of the following.

\begin{definition}%
\label{def:k-means-pop}
Suppose that $X$ satisfies $\mathbb{E}\|X\|<\infty $. The
\textit{(population) excess distortion} of $\mathcal{C}$ over
$\mathcal{C}'$ is defined to be
\begin{equation*}
D_{X}(\mathcal{C}\,|\,\mathcal{C}'):=\mathbb{E}\left [\min _{c\in
\mathcal{C}}\|c-X\|^{2}-\min _{c'\in \mathcal{C}'}\|c'-X\|^{2}\right ].
\end{equation*}
A set $\mathcal{C}\in \mathscr{C}_{k}$ is called a
\textit{set of population $k$-means cluster centers} if it satisfies
$D_{X}(\mathcal{C}\,|\,\mathcal{C}')\le 0$ for all
$\mathcal{C}'\in \mathscr{C}_{k}$. We write $\mathscr{C}_{k}(X)$ for the
collection of all sets of population $k$-means cluster centers.
\end{definition}

We can also formulate a definition similar to equation
\eqref{eqn:def-pop-k-means} under the assumption of
$\mathbb{E}\|X\|<\infty $. Namely, for any fixed
$\mathcal{C}_{0}\in \mathscr{C}_{k}$, it turns out that
$\mathscr{C}_{k}(X)$ is equal to the set of solutions to the optimization
problem
%
\begin{equation}
\label{eqn:extended-def}
\begin{cases}
\textnormal{minimize} &\mathbb{E}\left [\min _{c\in \mathcal{C}}\|c-X
\|^{2}-\min _{c_{0}\in \mathcal{C}_{0}}\|c_{0}-X\|^{2}\right ]
\\
\textnormal{over} &\mathcal{C}\in \mathscr{C}_{k}.
\end{cases}
%
\end{equation}
Interestingly, note that the resultant solution set does not depend on
the choice of $\mathcal{C}_{0}\in \mathscr{C}_{k}$, so we may simply take
$\mathcal{C}_{0} = \{0\}$. This can be easily proven from the elementary
identity
$D_{X}(\mathcal{C}\,|\,\mathcal{C}'') = D_{X}(\mathcal{C}\,|\,
\mathcal{C}') + D_{X}(\mathcal{C}'\,|\,\mathcal{C}'')$ for all
$\mathcal{C}, \mathcal{C}', \mathcal{C}''\in \mathscr{C}_{k}$.

We now define these notions at the empirical level, as follows.
%
\begin{definition}%
\label{def:k-means-emp}
The \textit{empirical excess distortion} of $\mathcal{C}$ over
$\mathcal{C}'$ is defined to be
\begin{equation*}
\bar D_{n}(\mathcal{C}\,|\,\mathcal{C}'):=\frac{1}{n}\sum _{i=1}^{n}
\left (\min _{c\in \mathcal{C}}\|c-X_{i}\|^{2}-\min _{c'\in
\mathcal{C}'}\|c'-X_{i}\|^{2}\right ).
\end{equation*}
A set $\mathcal{C}\in \mathscr{C}_{k}$ is called a
\textit{set of empirical $k$-means cluster centers} if it satisfies
$\bar D_{n}(\mathcal{C}\,|\,\mathcal{C}')\le 0$ for all
$\mathcal{C}'\in \mathscr{C}_{k}$. We write
$\bar {\mathscr{C}}_{n,k}$ for the collection of all sets of empirical
$k$-means cluster centers.
\end{definition}

The definition above is equivalent to \eqref{eqn:def-emp-k-means}, and
we emphasize that no renormalization is needed. However, we choose to formulate
the problem in terms of empirical excess distortion because the parallel
with (population) excess distortion will be useful in our later proofs.

Note that there is a slight abuse of notation in the preceding definitions.
That is, $D_{X}$ depends not on the random variable $X$ but rather on its
distribution. Similarly, $\bar D_{n}$ depends only on the empirical measure
$\frac{1}{n}\sum _{i=1}^{n}\delta _{X_{i}}$ of
$X_{1},\ldots , X_{n}$. Later we will see that there is some analytical
advantage to expressing the excess distortion in terms of an arbitrary
probability measure. (See Lemma~\ref{lem:D-cts}.)

\subsection{Voronoi regions}
\label{sec2.2}

Next we describe the following geometric aspect of $k$-means clustering
which features in our results.

\begin{definition}
\label{defn2.4}
Fix any $\mathcal{C}\in \mathscr{C}_{k}$. For $c\in \mathcal{C}$, the
\textit{open Voronoi region of $c$ in $\mathcal{C}$} is the set
\begin{equation*}
\mathcal{V}_{\mathcal{C}}^{\circ}(c) :=\{x\in \mathbb{R}^{d}: \|c-x\|<
\|c'-x\| \textnormal{ for all }c'\in \mathcal{C}\setminus \{c\}\},
\end{equation*}
consisting of all points in $\mathbb{R}^{d}$ which are (strictly) closer
to $c\in \mathcal{C}$ than to any other point of $\mathcal{C}$. Now fix
a total order $\preceq $ on $\mathbb{R}^{d}$ (e.g., lexicographic order)
and define \textit{the Voronoi partition of $c$ in $\mathcal{C}$} as
\begin{align*}
\mathcal{V}_{\mathcal{C}}(c) &:=\{x\in \mathbb{R}^{d}: \|c-x\|< \|c'-x
\| \textnormal{ for all }c\succ c'\in \mathcal{C}\}
\\
&\qquad \cap \{x\in \mathbb{R}^{d}: \|c-x\|\le \|c'-x\|
\textnormal{ for all }c\prec c'\in \mathcal{C}\},
\end{align*}
In words, $\{\mathcal{V}_{\mathcal{C}}(c)\}_{c\in \mathcal{C}}$ is a partition
of $\mathbb{R}^{d}$ into measurable sets such that we have
\begin{equation*}
\mathcal{V}_{\mathcal{C}}(c) \supseteq \mathcal{V}^{\circ}_{
\mathcal{C}}(c)
\end{equation*}
for all $c\in \mathcal{C}$, and where ties are broken according to the
ordering of the centers under $\preceq $.
\end{definition}

The interpretation is that sets of points represent
\textit{cluster centers}, while Voronoi regions correspond to
\textit{clusters} themselves. Also, we emphasize that the collection of
open Voronoi regions
$\{\mathcal{V}_{\mathcal{C}}^{\circ}(c)\}_{c\in \mathcal{C}}$ (which is
not a partition, since it does not cover $\mathbb{R}^{d}$) is determined
uniquely by $\mathcal{C}$, but this is not a problem in most cases: It is well-known (see
\cite[Theorem~1.5]{GrafLuschgy}) that the boundary of each Voronoi region
has zero Lebesgue measure, and (see \cite[Theorem~4.2]{GrafLuschgy}) that
the boundary of each Voronoi region has probability zero under $X$ whenever
$\mathcal{C}\in \mathscr{C}_{k}(X)$.

Next we introduce some notation for the conditional mean of $X$ given a
Voronoi partition. In the following, when we write
\textit{a Voronoi partition}
$\mathscr{V} := \{\mathcal{V}_{\ell}\}_{1\le \ell \le k}$ we mean a Voronoi
partition corresponding to \textit{some} set of cluster centers, i.e.,
$\{\mathcal{V}_{\ell}\}_{1\le \ell \le k} = \{\mathcal{V}_{
\mathcal{C}}(c)\}_{c\in \mathcal{C}}$ for some
$\mathcal{C}\in \mathscr{C}_{k}$.

\begin{definition}
\label{defn2.5}
If $X$ satisfies $\mathbb{E}\|X\|<\infty $ then for any measurable set
$\mathcal{V}\subseteq \mathbb{R}^{d}$, we write
\begin{equation*}
m_{X}(\mathcal{V}) :=\mathbb{E}[X\,|\,X\in \mathcal{V}]
\end{equation*}
and for any Voronoi partition
$\mathscr{V}=\{\mathcal{V}_{\ell}\}_{1\le \ell \le k},$ we write
\begin{equation*}
m_{X}(\mathscr{V}) := \{m_{X}(\mathcal{V}_{\ell})\}_{1\le \ell \le k},
\end{equation*}
which is an element of $\mathscr{C}_{k}$. We also write
\begin{equation*}
\bar m_{n}(\mathcal{V}) :=
\frac{\sum _{i=1}^{n}X_{i}\mathbbold{1}\{X_{i}\in \mathcal{V}\}}{\sum _{i=1}^{n}\mathbbold{1}\{X_{i}\in \mathcal{V}\}}
\end{equation*}
and $\bar m_{n}(\mathscr{V})$ for the corresponding conditional means with
respect to the empirical distribution of $X_{1},\ldots , X_{n}$.
\end{definition}

The preceding definition does not make sense if $\mathcal{V}$ has zero
probability under the distribution of $X$ or if $\mathcal{V}$ contains
none of the data points $X_{1},\ldots , X_{n}$. In such cases, the conditional
means $m_{X}(\mathcal{V})$ and $\bar m_{n}(\mathcal{V})$ may be arbitrarily
set to $0$.

Lastly, we give a fundamental result relating a given Voronoi partition
to the cluster centers to which it corresponds, in the case that the cluster
centers are an optimal set of population $k$-means cluster centers. This
is usually referred to as the \textit{centroid condition}, and a proof can
be found in \cite[Lemma~A]{PollardCLT} or
\cite[Lemma~4.10]{GrafLuschgy}; note that both proofs rely only on
$\mathbb{E}\|X\|<\infty $ although they have assumed the stronger
$\mathbb{E}\|X\|^{2}<\infty $.

\begin{lemma}%
\label{lem:centroid-condition}
If $X$ satisfies $\mathbb{E}\|X\|<\infty $ and
$\mathcal{C}\in \mathscr{C}_{k}(X)$, then
$\mathscr{V}:=\{\mathcal{V}_{\mathcal{C}}(c)\}_{c\in \mathcal{C}}$ satisfies
$m_{X}(\mathscr{V}) = \mathcal{C}$. In other words, we have
\begin{equation*}
c = \mathbb{E}[X\,|\,X\in \mathcal{V}_{\mathcal{C}}(c)]
\end{equation*}
for all $c\in \mathcal{C}$.
\end{lemma}

\subsection{Balance and constraints}
\label{subsec:balance}

Our positive results concern a modification of $k$-means clustering in
which the clusters are constrained to have some degree of balance, so we
now give precise definitions for various notions related to balance. Because
of Lemma~\ref{lem:centroid-condition}, it is well-known that we may equivalently
formulate the (unconstrained) $k$-means clustering problem in terms of
either the cluster centers or the clusters. However, adding balance constraints
means that the resulting two problems need not agree, and we must choose
which one to take as our definition; in this work, we choose to consider
the problem of optimizing over Voronoi partitions satisfying some balance
constraint.

We begin by introducing some further notation. First, we write
$\mathscr{P}_{k}$ for the set of all Voronoi partitions of
$\mathbb{R}^{d}$ into at most $k$ regions; in other words,
$\mathscr{V}\in \mathscr{P}_{k}$ if and only if we have
$\mathscr{V} = \{\mathcal{V}_{\mathcal{C}}(c)\}_{c\in \mathcal{C}}$ for
some $\mathcal{C}\in \mathscr{C}_{k}$. Then we give the following:

\begin{definition}%
\label{def:Voronoi-distortion}
Suppose that $X$ satisfies $\mathbb{E}\|X\|<\infty $. The
\textit{(population) excess distortion} of
$\mathscr{V}\in \mathscr{P}_{k}$ is
%
\begin{equation}
\label{eqn:Voronoi-distortion}
\begin{split} D_{X}(\mathscr{V}) &:=\mathbb{E}\left [ \sum _{\ell =1}^{k}
\mathbbold{1}\{X\in \mathcal{V}_{\ell }\} \|m_{X}(\mathcal{V}_{\ell })-X
\|^{2} - \|X\|^{2} \right ]
\\
&\qquad =-\sum _{\ell =1}^{k} \mathbb{P}(X\in \mathcal{V}_{\ell}) \|m_{X}(
\mathcal{V}_{\ell})\|^{2}.
\end{split}
%
\end{equation}
(More precisely, $D_{X}(\mathscr{V})$ is the excess distortion of
$\mathscr{V}$ relative to the partition
$\mathscr{V}_{0} := \{\mathbb{R}^{d}\}$, which is generated by
$\{0\}$. However, in this paper, we will not consider the excess distortion
with respect to more general Voronoi partitions.)
\end{definition}

Note that~\eqref{eqn:Voronoi-distortion} provides two equivalent formulations
of the excess distortion of $\mathscr{V}$. Although the equivalence of
these definitions is elementary, we will see throughout the paper that
it turns out to be very useful to have both formulations. This observation
has been made in previous works under the stronger condition of
$\mathbb{E}\|X\|^{2}<\infty $ (e.g.,
\cite[Remark~4.6(c)]{GrafLuschgy}), but we emphasize that it also holds
under the weaker $\mathbb{E}\|X\|<\infty $.

We also note that, although $\mathscr{V}$ need not be equal to the Voronoi
partition generated by the conditional means $m_{X}(\mathscr{V})$, we nonetheless
have the following well-known result which states that the corresponding
unconstrained problems have the same value. We emphasize that this may
be regarded as an extension of Lemma~\ref{lem:centroid-condition} to the
case where optimal clusters need not exist.
%
\begin{lemma}%
\label{lemma:equivalence_distortions}
If $X$ satisfies $\mathbb{E}\|X\|<\infty $, then
$\inf _{\mathcal{C}\in \mathscr{C}_{k}}D_{X}(\mathcal{C}\,|\,
\mathcal{C}_{0}) = \inf _{\mathscr{V}\in \mathscr{P}_{k}} D_{X}(
\mathscr{V})$ for $\mathcal{C}_{0}=\{0\}$.
\end{lemma}

Similarly, we define analogous notions for the empirical-level problem.
As usual, this can be seen as a special case of the population-level problem
when the distribution of $X$ is assumed to be equal to the empirical measure
of the samples $X_{1},\ldots , X_{n}$.

\begin{definition}%
\label{defn2.9}
The \textit{empirical excess distortion} of
$\mathscr{V}\in \mathscr{P}_{k}$ is defined as
%
\begin{equation}
\label{eqn:emp-excess-dist}
\begin{split} \bar D_{n}(\mathscr{V}) &:=\frac{1}{n}\sum _{i=1}^{n}
\sum _{\ell =1}^{k} \mathbbold{1}\{X_{i}\in \mathcal{V}_{\ell}\}
\left ( \|\bar m_{n}(\mathcal{V}_{\ell})-X_{i}\|^{2}-\|X_{i}\|^{2}
\right )
\\
&\qquad = -\frac{1}{n}\sum _{i=1}^{n}\sum _{\ell =1}^{k}
\mathbbold{1}\{X_{i}\in \mathcal{V}_{\ell}\} \|\bar m_{n}(\mathcal{V}_{
\ell})\|^{2}.
\end{split}
%
\end{equation}
\end{definition}

Now we may define the central notion of balanced $k$-means clustering.
Because we will use empirical balanced $k$-means clustering to consistently
estimate the cluster centers in the unconstrained population $k$-means
clustering problem, we focus for now on balance constraints at the empirical
level.
%
\begin{definition}%
\label{def:balance-k-means-emp}
Fix some $0 \le \gamma \le n$. A Voronoi partition $\mathscr{V}$ is called
a \textit{$\gamma $-balanced Voronoi partition} if it satisfies
$\#\{1\le i \le n:X_{i}\in \mathcal{V}\}\ge \gamma $ for all
$\mathcal{V}\in \mathscr{V}$, and we write
$\bar{\mathscr{P}}_{n,k,\gamma}$ for the set of all $\gamma $-balanced
Voronoi partitions with at most $k$ regions. An element
$\mathscr{V}\in \bar{\mathscr{P}}_{n,k,\gamma}$ is called an
\textit{empirical $\gamma $-balanced $k$-means clustering Voronoi partition}
if it satisfies
$\bar D_{n}(\mathscr{V}) \le \bar D_{n}(\mathscr{V}')$ for all
$\mathscr{V}'\in \bar{\mathscr{P}}_{n,k,\gamma}$. We write
$\bar {\mathscr{C}}_{n,k,\gamma}:=\{\bar m_{n}(\mathscr{V}):
\mathscr{V}\in \bar{\mathscr{P}}_{n,k,\gamma}, \mathscr{V}
\text{ is an empirical $\gamma $-balanced $k$-means clustering Voronoi partition}
\}$, any element of which is called a
\textit{$\gamma $-balanced set of empirical $k$-means cluster centers}.
\end{definition}

Although the preceding definition is a bit involved, the interpretation
is straightforward: A set of empirical $\gamma $-balanced $k$-means cluster
centers is nothing more than the set of conditional means of the empirical
distribution of the data, with respect to some empirical $\gamma $-balanced
$k$-means clustering Voronoi partition. Also, we note that there always
exist some set of $\gamma $-balanced empirical $k$-means cluster centers,
since there are only finitely many partitions of the $n$ data points. We
lastly note that there are some subtleties to this definition, explained
in the following:

\begin{remark}%
\label{rem:def-balance}
The Voronoi partition generated by $\bar m_{n}(\mathscr{V})$ need not be
$\gamma $-balanced, even though $\mathscr{V}$ is required to be
$\gamma $-balanced. This discrepancy highlights the fact that there is
not a canonical way to impose the balance constraint; while we have chosen
one particular formalism of the balanced clustering problem for convenience
in our proofs, we note that other choices may also be interesting.
\end{remark}

\section{Population-level trichotomy}
\label{sec:cases}

In this section, we show that the population-level $k$-means clustering
problem \eqref{eqn:def-pop-k-means-minimal} is naturally divided into three
cases, as outlined in the introduction. That is, the distribution of
$X$ may be such that problem \eqref{eqn:def-pop-k-means-minimal} has (i)
a solution hence finite optimal excess distortion, (ii) a finite optimal
excess distortion but no solution, or (iii) infinite optimal excess distortion
hence no solution. Throughout this section, $X$ always denotes a random
variable satisfying $\mathbb{E}\|X\|<\infty $.

Most of our explicit examples involve distributions with power-law tails.
The basic starting point is a random variable $X$ that is said to have
a \textit{Pareto distribution with parameter $\alpha >0$}, denoted
$\textnormal{Par}(\alpha )$, meaning
\begin{equation*}
\mathbb{P}(X\ge t) = t^{-\alpha} \qquad \textnormal{ for all } t\ge 1.
\end{equation*}
(We always take the support of $X$ to be $[1,\infty )$ but other definitions
of Pareto distribution may involve support of the form
$[x_{0},\infty )$ for some $x_{0}>0$.) Our examples generally take place
in the simplest possible setting of $d=1$ and $k=2$ because of some explicit
calculations that can be done with Pareto distributions, although we comment
at a few points about how these examples may be generalized, particularly
to general dimension $d\ge 1$.

First we consider case (i), that problem~\eqref{eqn:def-pop-k-means-minimal}
admits a solution, and hence the optimal excess distortion is finite. By
the general results of \cite{Pollard}, this case obtains whenever we have
$\mathbb{E}\|X\|^{2}<\infty $; as we show next, it can also obtain in cases
with $\mathbb{E}\|X\|^{2} = \infty $.
%
\begin{example}%
\label{ex:2-sided-pareto-pop}
Suppose $X$ is a real-valued random variable whose distribution is symmetric
and satisfies
\begin{equation*}
\mathbb{P}(|X|\ge t) = t^{-2} \qquad \textnormal{for all}\qquad t
\ge 1.
\end{equation*}
(In other words, $X$ is symmetric and $|X|\sim \textnormal{Par}(2)$.) Then
$\mathscr{C}_{2}(X) = \{\{-2,2\}\}$.
\end{example}

\begin{proof}
By Lemma~\ref{lem:centroid-condition}, we may parameterize every putative
$\mathcal{C}\in \mathscr{C}_{2}(X)$ by the set of possible Voronoi regions.
Since $d=1$ and $k=2$, the only possible Voronoi regions are of the form
$\{(-\infty ,r),[r,\infty )\}$ for $r\in \mathbb{R}$. (Recall by
\cite[Theorem~4.2]{GrafLuschgy} that we may assign the boundary
$\{r\}$ arbitrarily to either region; we choose the right region for no
particular reason.) Thus it follows that, if there exists some
$\mathcal{C}\in \mathscr{C}_{k}(X)$, then it may be written
$\mathcal{C}= \mathcal{C}(r)$, where
\begin{equation*}
\min (\mathcal{C}(r)) := \mathbb{E}[X\,|X< r]\qquad
\textnormal{ and }\qquad \max (\mathcal{C}(r)) = \mathbb{E}[X\,|X\ge r].
\end{equation*}
We can directly calculate these values, by considering cases for
$r\in \mathbb{R}$. If $-1\le r\le 1$, then we have
$\min (\mathcal{C}(r)) = -2$ and $\max (\mathcal{C}(r)) = 2$. If
$r>1$, we have
\begin{equation*}
\min (\mathcal{C}(r)) =-\frac{r^{-1}}{1-\frac{1}{2}r^{-2}} \qquad
\textnormal{ and }\qquad \max (\mathcal{C}(r)) = 2r,
\end{equation*}
and the $r<-1$ case can be handled by symmetry. Therefore, the function
$D:\mathbb{R}\to \mathbb{R}$ defined via
\begin{equation*}
D(r) := \mathbb{E}\left [\min _{c\in \mathcal{C}(r)}|c-X|^{2} - |X|^{2}
\right ]
\end{equation*}
is given by
\begin{equation*}
D(r)=
\begin{cases}
-4 &\textnormal{ if } |r|\le 1
\\
-2r^{2}/(r^{2}-\frac{1}{2}) &\textnormal{ if } |r|> 1.
\end{cases}
\end{equation*}
This function achieves its global minimum on all of $-1\le r\le 1$, hence
we have $\mathscr{C}_{2}(X) = \{\{-2,2\}\}$ as claimed.
\end{proof}

This example may be extended to higher dimensions using some general results;
for example, if $X$ is an $\mathbb{R}^{d}$-valued random variable whose
coordinates are independent symmetric $\textnormal{Par}(2)$ random variables,
then \cite[Lemma~4.15]{GrafLuschgy} implies that
$\mathscr{C}_{2^{d}}(X) = \{\{-2,2\}^{d}\}$. However, it is not clear what
to expect from $\mathscr{C}_{k}(X)$ when $k$ is not a power of $2$; see
the examples of \cite[Section~4.4]{GrafLuschgy} for an illustration of
the difficulty of the resulting calculations.

When case (i) holds, it is natural to hope that each set of empirical
$k$-means cluster centers should be close to some set of population
$k$-means cluster centers. This question is taken up in detail in Section~\ref{sec:centers}.

Second, we consider case (ii), that problem~\eqref{eqn:def-pop-k-means-minimal}
has a finite optimal excess distortion but no solution. To see that this
is possible, we consider the following.
%
\begin{example}%
\label{ex:1-sided-pareto-pop}
Suppose $X$ is a real-valued random variable whose distribution satisfies
\begin{equation*}
\mathbb{P}(X\ge t) = t^{-2} \qquad \textnormal{for all}\qquad t \ge 1.
\end{equation*}
(In other words, $X$ $\sim \textnormal{Par}(2)$.) Then
$\mathscr{C}_{2}(X) = \varnothing $.
\end{example}

\begin{proof}
The proof is similar to the proof of Example~\ref{ex:2-sided-pareto-pop}
above. By Lemma~\ref{lem:centroid-condition}, we parameterize putative
$\mathcal{C}\in \mathscr{C}_{2}(X)$ by
$\{(-\infty ,r),[r,\infty )\}$ for $r\in \mathbb{R}$. If $r<1$ then
$\mathbb{P}(X<r) = 0$ meaning one cluster center is not well-defined, so
we restrict our attention to the case of $r\ge 1$. By direct calculation
for $r\ge 1$, we get
\begin{equation*}
\min (\mathcal{C}(r)) := \mathbb{E}[X\,|X\le r] =
\frac{2(1-r^{-1})}{1-r^{-2}}\qquad \textnormal{ and }\qquad \max (
\mathcal{C}(r)) := \mathbb{E}[X\,|X>r] = 2r.
\end{equation*}
Next, we calculate that the function $D:[1,\infty )\to \mathbb{R}$ defined
via
\begin{equation*}
D(r) := \mathbb{E}\left [\min _{c\in \mathcal{C}(r)}|c-X|^{2} - |X|^{2}
\right ]
\end{equation*}
is just equal to $D(r) = -8r/(r+1)$. Note that this decreases monotonically
to $-8$ as $r\to \infty $ but that there exists no
$r\in [1,\infty )$ which achieves value $-8$. This shows that no Voronoi
partition is optimal, hence proves the claim.
\end{proof}

This example can also be extended to higher dimensions; if $X$ is an
$\mathbb{R}^{d}$-valued random variable whose coordinates are independent
$\textnormal{Par}(2)$ random variables, then from
\cite[Lemma~4.15]{GrafLuschgy} we obtain
$\mathscr{C}_{2^{d}}(X) = \varnothing $. However, it is still not clear
how to analyze $\mathscr{C}_{k}(X)$ when $k\neq 2^{d}$.

When case (ii) holds, it is natural to hope that each set of empirical
$k$-means cluster centers should achieve a (population) excess distortion
which is similar to the infimal excess distortion. This question is taken
up in detail in Section~\ref{sec:distortion}.

Some remarks about the preceding two examples are due. First, notice that
both examples involved a random variable $X$ satisfying
$\mathbb{E}|X|^{2} = \infty $ but $\mathbb{E}|X|^{p}<\infty $ for all
$1\le p < 2$; thus, cases (i) and (ii) can obtain even in settings where
the finite-variance condition $\mathbb{E}|X|^{2}<\infty $
\textit{just barely} fails. Second, we emphasize that we do not know a simple
way to determine whether problem \eqref{eqn:def-pop-k-means-minimal} admits
a solution (equivalently, which of (i) or (ii) holds); in both previous
examples our analysis relied on exactly calculating the excess distortion,
but this is not possible in more complicated examples.

It is also of interest to study case (iii) that problem~\eqref{eqn:def-pop-k-means-minimal}
has infinite optimal excess distortion hence no solution. In order to not
distract from the main ideas of the paper, we defer this discussion to
Appendix~\ref{sec:case-iii} since things are slightly subtle. While we
can give a characterization of when case (iii) occurs (Lemma~\ref{lemma:charact_finite_obj}),
we do not know an interesting form of consistency that may be desired,
since essentially all population-level objects are ill-defined.

\section{Empirical-level convergence of cluster centers}
\label{sec:centers}

In this section we consider the problem of whether the empirical $k$-means
cluster centers converge, in some sense, to some population $k$-means cluster
centers in case (i). We will show that this naive guess can fail (Example~\ref{ex:inconsistency_two_sided}),
but also that imposing a certain linear balance constraint indeed guarantees
such convergence (Theorem~\ref{thm:cluster-consistency}). Throughout this
section, we let $X$ denote a random variable such that case (i) holds;
that is, there exists a set of population $k$-means cluster centers, and
hence the optimal excess distortion is finite.

Notationally, it is slightly complicated to precisely state a notion of
convergence for $k$-means clustering. One reason is that the objects of
interest are sets of points, so we need a notion of convergence for sets
that is strong enough to be practically meaningful. Another reason is that
there typically does not exist a unique set of $k$-means cluster centers
for either the empirical or population problems, so we need a notion of
convergence which accommodates this ill-posedness. In order to overcome
these challenges, we follow the form of consistency given in most existing
positive results
\cite{Pollard, JaffeClustering, Lember, Parna1, Parna2}, which we now define.

This form of consistency requires the following notion of distance between
non-empty compact subsets of $\mathbb{R}^{d}$, which was used in the earliest
consistency results for $k$-means clustering by Pollard
\cite{Pollard}.

\begin{definition}
\label{defn4.1}
For any non-empty compact sets
$\mathcal{C},\mathcal{C}'\subseteq \mathbb{R}^{d}$, we define
\begin{equation*}
d_{\textnormal{H}}(\mathcal{C},\mathcal{C}') := \max \left \{\max _{c
\in \mathcal{C}}\min _{c'\in \mathcal{C}'}\|c-c'\|,\max _{c'\in
\mathcal{C}'}\min _{c\in \mathcal{C}}\|c-c'\|\right \},
\end{equation*}
called their \textit{Hausdorff distance}.
\end{definition}

It is worth mentioning that there are several equivalent formulations of
the Hausdorff distance. For instance, if
$\#\mathcal{C}= \#\mathcal{C}' = k<\infty $, then we have
$d_{\textnormal{H}}(\mathcal{C},\mathcal{C}') \le \varepsilon $ if and
only if there exist labelings
$\mathcal{C}= \{c_{1},\ldots , c_{k}\}$ and
$\mathcal{C}' = \{c_{1}',\ldots , c_{k}'\}$ satisfying
$\|c_{\ell}-c_{\ell}'\|\le \varepsilon $ for all $1\le \ell \le k$. Alternatively,
if we define the \textit{$\varepsilon $-thickening} of $\mathcal{C}$ via
\begin{equation*}
\mathcal{C}^{\varepsilon}:=\bigcup _{c\in \mathcal{C}}\{x\in
\mathbb{R}^{d}: \|c-x\|\le \varepsilon \}
\end{equation*}
then $d_{\textnormal{H}}(\mathcal{C},\mathcal{C}')$ is exactly the smallest
$\varepsilon \ge 0$ such that we have both
$\mathcal{C}^{\varepsilon}\supseteq \mathcal{C}'$ and
$(\mathcal{C}')^{\varepsilon}\supseteq \mathcal{C}$.

Now we may state that a natural form of consistency, following
\cite{Pollard, JaffeClustering, Lember, Parna1, Parna2}, is
%
\begin{equation}
\label{eqn:centers-consistent}
\max _{\bar{\mathcal{C}}_{n}\in \bar{\mathscr{C}}_{n,k}}\min _{
\mathcal{C}\in \mathscr{C}_{k}(X)}d_{\textnormal{H}}(
\bar{\mathcal{C}}_{n},\mathcal{C})\to 0
\end{equation}
almost surely, as $n\to \infty $. We think of~\eqref{eqn:centers-consistent}
as a ``no false positives'' property, since it shows that every set of
empirical $k$-means cluster centers must be close, with respect to the
Hausdorff distance $d_{\textnormal{H}}$, to some set of population
$k$-means cluster centers. If there indeed exist uniquely defined sets
of empirical and population $k$-means cluster centers, denoted
$\{\bar {\mathcal{C}}_{n}\}_{n\in \mathbb{N}}$ and $\mathcal{C}$ respectively,
then~\eqref{eqn:centers-consistent} is equivalent to
$d_{\textnormal{H}}(\bar {\mathcal{C}}_{n},\mathcal{C})\to 0$ almost surely,
as $n\to \infty $.

Unfortunately, we show in an example that~\eqref{eqn:centers-consistent}
can plainly fail in case (i); see Figure~\ref{fig:two_sided_pareto} for
a visualization of this example, whose proof is somewhat long and is deferred
to Appendix~\ref{sec:proofs}.

\begin{example}%
\label{ex:inconsistency_two_sided}
Suppose $X$ is a real-valued random variable whose distribution is symmetric
and satisfies
\begin{equation*}
\mathbb{P}(|X|\ge t) = t^{-2} \qquad \textnormal{for all}\qquad t
\ge 1.
\end{equation*}
(In other words, $X$ is symmetric and $|X|\sim \textnormal{Par}(2)$.) Then,
there exists a universal constant $c>0$ such that almost surely, for infinitely
many $n\in \mathbb{N}$, for any
$\bar{\mathcal{C}}_{n}\in \bar {\mathscr{C}}_{n,2}$ we have
\begin{equation*}
\max (\bar{\mathcal{C}}_{n}) \ge c\,\frac{\sqrt{n}}{\log n} \quad
\textnormal{or} \quad \min (\bar{\mathcal{C}}_{n}) \le -c\,
\frac{\sqrt{n}}{\log n}.
\end{equation*}
Consequently, we have almost surely that, for infinitely many
$n\in \mathbb{N}$,
\begin{equation*}
\min _{\bar{\mathcal{C}}_{n}\in \bar{\mathscr{C}}_{n,2}} \min _{
\mathcal{C}\in \mathscr{C}_{2}(X)} d_{\textnormal{H}}(
\bar{\mathcal{C}}_{n},\mathcal{C}) \geq c\frac{\sqrt n}{\log n}-2.
\end{equation*}
\end{example}

\begin{remark}
\label{rem4.2}
Although the proof of Example~\ref{ex:inconsistency_two_sided} is somewhat
long, the big picture is, roughly speaking, that the following holds with
high probability:
%
\begin{equation}
\label{eqn:counter-ex-sketch}
\left \{\frac{X_{(n)}}{X_{(n-1)}}\ge \sqrt{n}\right \}\subseteq
\Big\{X_{(n)}\textnormal{ is in a cluster alone}\Big\},
\end{equation}
where $X_{(1)}<\cdots <X_{(n)}$ denote the order statistics of
$X_{1},\ldots , X_{n}$. To gain some intuition for
\eqref{eqn:counter-ex-sketch}, note that if
$X_{(n)}/X_{(n-1)}\ge \sqrt{n}$ but $X_{(n)}$ is not in a cluster alone,
then the distance $\ell $ from $X_{(n)}$ to the nearest cluster center
is at least
$\ell \ge \frac{1}{2}(X_{(n)}-X_{(n-1)}) \ge \frac{1}{2}(\sqrt{n}-1)X_{(n-1)}$;
this implies that the decrease in excess distortion resulting from moving
$X_{(n)}$ into a cluster alone is at least
$\frac{1}{n}\ell ^{2} \ge \frac{1}{4n}(\sqrt{n}-1)^{2}X_{(n-1)}^{2}
\sim \frac{1}{4}X_{(n-1)}^{2}\to \infty $, which is a contradiction since
the excess distortion is itself of constant order. At the same time, the
heavy tails of the $\textnormal{Par}(2)$ distribution ensure that
$X_{(n)}/X_{(n-1)}\ge \sqrt{n}$ happens infinitely often almost surely.
\end{remark}

The upshot of Example~\ref{ex:inconsistency_two_sided} comes from comparing
with Example~\ref{ex:2-sided-pareto-pop}; in the symmetric two-sided Par(2)
distribution, there is a unique set of population-level 2-means cluster
centers, but the empirical 2-means cluster centers do not converge to it,
and even oscillate to $\pm \infty $ as $n\to \infty $.

Lastly, we turn to positive results that can be achieved by imposing some
balance constraint on the empirical $k$-means clusters. To set this up
precisely, we need the following:

\begin{definition}
\label{defn4.3}
For any Voronoi partition $\mathscr{V}$, write
\begin{equation*}
p_{\min}(X,\mathscr{V}) := \min _{\mathcal{V}\in \mathscr{V}}
\mathbb{P}(X\in \mathcal{V})
\end{equation*}
for the smallest probability assigned to any of its regions.
\end{definition}

The quantity $p_{\min}(X,\mathscr{V})$ can be interpreted as the degree
of imbalance of the clusters in $\mathscr{V}$ at the population level.
In fact, it is easy to see that we always have
$0 \le p_{\min}(X,\mathscr{V})\le 1/\#\mathscr{V}$, and the upper bound
is achieved if and only if all clusters have the same mass.

Before giving our positive result, it is instructive to state and prove
a fundamental continuity property for the excess distortion function. Towards
this, observe that the quantity $D_{X}(\mathcal{C}\,|\mathcal{C}')$ depends
not on the random variable $X$ but rather on its distribution; as such,
the next result shows that the excess distortion is jointly continuous
as a function of this distribution and its two sets of points. (In other
words, this is an extension of \cite[Lemma~4]{JaffeClustering} from
$W_{2}$ to $W_{1}$ convergence.)

\begin{lemma}%
\label{lem:D-cts}
The excess distortion $D$ is a continuous function
$(\mathcal{W}_{1}(\mathbb{R}^{d})\times \mathscr{C}_{k}\times
\mathscr{C}_{k},W_{1}\times d_{\textnormal{H}}\times d_{
\textnormal{H}})\to \mathbb{R}$.
\end{lemma}

Finally we get to our positive result. Recall that we are in case (i) so
we are assuming that an optimal set of (unconstrained) $k$-means cluster
centers exists.
%
\begin{theorem}%
\label{thm:cluster-consistency}
Suppose $X$ has a density with respect to Lebesgue measure,
$\mathbb{E}\|X\|<\infty $, and $\#\mathrm{supp}(X)\ge k$, and let
\begin{equation*}
\alpha _{0}(X):=\sup \left \{p_{\min}(X,\mathscr{V}):
\begin{matrix}
\mathscr{V} \textnormal{ is the Voronoi partition for some}
\\
\textnormal{optimal set of $k$-means cluster centers}
\end{matrix}
\right \} > 0.
\end{equation*}
If $0 < \alpha < \alpha _{0}(X)$, then we have
\begin{equation*}
\max _{\bar{\mathcal{C}}_{n}\in \bar{\mathscr{C}}_{n,k,\alpha n}}
\min _{\mathcal{C}\in \mathscr{C}_{k}(X)}d_{\textnormal{H}}(
\bar{\mathcal{C}}_{n},\mathcal{C})\to 0
\end{equation*}
almost surely as $n\to \infty $.
\end{theorem}

\begin{proof}
It suffices to show that, if
$\{\bar {\mathcal{C}}_{n}\}_{n\in \mathbb{N}}$ is any sequence in
$\mathscr{C}_{k}$, with
$\bar {\mathcal{C}}_{n}\in \bar {\mathscr{C}}_{n,k,\alpha n}$, then
$\{\bar {\mathcal{C}}_{n}\}_{n\in \mathbb{N}}$ is
$d_{\textnormal{H}}$-precompact and every $d_{\textnormal{H}}$-subsequential
limit lies in $\mathscr{C}_{k}(X)$, almost surely. For notation, let us
write $\bar{\mathscr{V}}_{n}$ for each $n\in \mathbb{N}$ for the Voronoi
partition such that
$\{\bar m_{n}(\mathcal{V})\}_{\mathcal{V}\in \bar{\mathscr{V}}_{n}}=
\bar {\mathcal{C}}_{n}$.

First, we show the required precompactness. To do this, we make the following
simple calculation for any $\mathcal{V}\in \bar{\mathscr{V}}_{n}$, which
is similar to \cite[equation~(14)]{Zhivotovskiy}:
\begin{equation*}
\|\bar m_{n}(\mathcal{V})\| =
\frac{\|\sum _{i=1}^{n}X_{i}\mathbbold{1}\{X_{i}\in \mathcal{V}\}\|}{\sum _{i=1}^{n}\mathbbold{1}\{X_{i}\in \mathcal{V}\}}
\le
\frac{\sum _{i=1}^{n}\|X_{i}\|\mathbbold{1}\{X_{i}\in \mathcal{V}\}}{\sum _{i=1}^{n}\mathbbold{1}\{X_{i}\in \mathcal{V}\}}
\le
\frac{\sum _{i=1}^{n}\|X_{i}\|}{\sum _{i=1}^{n}\mathbbold{1}\{X_{i}\in \mathcal{V}\}}.
\end{equation*}
Since $\bar{\mathscr{V}}_{n}$ is $\alpha n$-balanced by assumption, we
have shown
\begin{equation*}
\max _{\mathcal{V}\in \bar{\mathscr{V}}_{n}} \|\bar m_{n}(\mathcal{V})
\| \le \frac{1}{\alpha n}\sum _{i=1}^{n}\|X_{i}\|=:Y_{n}
\end{equation*}
almost surely. But $\mathbb{E}\|X\|<\infty $ and the strong law of large
numbers imply $Y_{n}\to \mathbb{E}\|X\|$ almost surely, so we have
\begin{equation*}
\sup _{n\in \mathbb{N}}\max _{\mathcal{V}\in \bar{\mathscr{V}}_{n}}
\|\bar m_{n}(\mathcal{V})\| \le \sup _{n\in \mathbb{N}}Y_{n} <
\infty
\end{equation*}
almost surely. This shows that
$\{\bar {\mathcal{C}}_{n}\}_{n\in \mathbb{N}}$ are almost surely contained
in a bounded set, hence $d_{\textnormal{H}}$-compactness follows from
$k$ iterated applications of the Bolzano-Weierstrass theorem. (Alternatively,
one may directly apply Blaschke's selection theorem
\cite[Theorem~7.3.8]{BuragoBuragoIvanov}.)

Next we suppose that we have a subsequence
$\{n_{j}\}_{j\in \mathbb{N}}$ and a compact set
$\mathcal{C}\subseteq \mathbb{R}^{d}$ such that
$d_{\textnormal{H}}(\bar {\mathcal{C}}_{n_{j}},\mathcal{C})\to 0$ almost
surely, and let us show $\mathcal{C}\in \mathscr{C}_{k}(X)$. To do this,
write $\mathscr{V}=\{\mathcal{V}_{\ell}\}_{1\le \ell \le k}$ for the Voronoi
partition generated by the cluster centers $\mathcal{C}$, and also write
$\bar{\mathscr{V}}_{n_{j}} = \{\mathcal{V}_{j,\ell}\}_{1\le \ell \le k}$.
We claim that we have
\begin{align*}
&\mathbb{E}\left [ \sum _{\ell =1}^{k} \mathbbold{1}\{X\in
\mathcal{V}_{\ell }\} \|m_{X}(\mathcal{V}_{\ell })-X\|^{2} - \|X\|^{2}
\right ]
\\
&\le \liminf _{j\to \infty}\frac{1}{n_{j}}\sum _{i=1}^{n_{j}}\sum _{
\ell =1}^{k}\mathbbold{1}\{X_{i}\in \mathcal{V}_{j,\ell}\}\left (\|
\bar m_{n_{j}}(\mathcal{V}_{j,\ell})-X_{i}\|^{2} - \|X_{i}\|^{2}
\right ).
\end{align*}
To see this, we apply the same uniform integrability trick as in the proof
of Lemma~\ref{lem:D-cts}, where we also use the fact that $X$ possessing
a density implies that we have
$\mathbbold{1}\{X\in \mathcal{V}_{j,\ell}\}\to \mathbbold{1}\{X\in
\mathcal{V}_{\ell}\}$ holding almost surely, after passing to a subsequence
and relabeling, if necessary. Then we let
$\mathscr{V}':=\{\mathcal{V}'_{\ell}\}_{1\le \ell \le k}$ be an arbitrary
Voronoi partition with $p_{\min}(X,\mathscr{V}') \ge \alpha _{0}(X)$, and
we note that $\alpha <\alpha _{0}(X)$ implies that $\mathscr{V}'$ is
$\alpha n$-balanced for sufficiently large $n\in \mathbb{N}$, almost surely.
Thus, we apply the displayed equation above, the optimality of
$\mathscr{V}$, and the strong law of large numbers, to get:
\begin{align*}
D_{X}(\mathscr{V}) &=\mathbb{E}\left [ \sum _{\ell =1}^{k}
\mathbbold{1}\{X\in \mathcal{V}_{\ell }\} \|m_{X}(\mathcal{V}_{\ell })-X
\|^{2} - \|X\|^{2} \right ]
\\
&\le \liminf _{j\to \infty}\frac{1}{n_{j}}\sum _{i=1}^{n_{j}}\sum _{
\ell =1}^{k}\mathbbold{1}\{X\in \mathcal{V}_{j,\ell}\}\left (\|\bar m_{n_{j}}(
\mathcal{V}_{j,\ell})-X_{i}\|^{2} - \|X_{i}\|^{2} \right )
\\
&= \liminf _{j\to \infty}-\frac{1}{n_{j}}\sum _{i=1}^{n_{j}}\sum _{
\ell =1}^{k}\mathbbold{1}\{X\in \mathcal{V}_{j,\ell}\}\|\bar m_{n_{j}}(
\mathcal{V}_{j,\ell})\|^{2}
\\
&\le \liminf _{j\to \infty}-\frac{1}{n_{j}}\sum _{i=1}^{n_{j}}\sum _{
\ell =1}^{k}\mathbbold{1}\{X\in \mathcal{V}'_{\ell}\}\|\bar m_{n_{j}}(
\mathcal{V}'_{\ell})\|^{2}
\\
&=-\sum _{\ell =1}^{k}\mathbb{P}(X\in \mathcal{V}'_{\ell})\| m_{X}(
\mathcal{V}'_{\ell})\|^{2}
\\
&= D_{X}(\mathscr{V}').
\end{align*}
We have shown that $\mathscr{V}$ minimizes $D_{X}$ among all Voronoi partitions
whose value of $p_{\min}$ is at least $\alpha _{0}(X)$, and by Lemma~\ref{lemma:equivalence_distortions}
this shows $\mathcal{C}\in \mathscr{C}_{k}(X)$.
\end{proof}

\begin{remark}
\label{rem4.6}
The condition that the constant $0<\alpha <1$ is sufficiently small should
be interpreted as a type of well-specification. However, note that we require
$\alpha $ to be smaller than
\begin{equation*}
\sup \left \{p_{\min}(X,\mathscr{V}):
\begin{matrix}
\mathscr{V} \textnormal{ is the Voronoi partition for some}
\\
\textnormal{optimal set of $k$-means cluster centers}
\end{matrix}
\right \}
\end{equation*}
rather than
\begin{equation*}
\inf \left \{p_{\min}(X,\mathscr{V}):
\begin{matrix}
\mathscr{V} \textnormal{ is the Voronoi partition for some}
\\
\textnormal{optimal set of $k$-means cluster centers}
\end{matrix}
\right \}.
\end{equation*}
These values may disagree in general, since there need not be a unique
set of population $k$-means cluster centers for $X$. In other words, the
well-specification depends on the population-level ``best-case'' imbalance
rather than the ``worst-case'' imbalance. (Relatedly, see the comments
before \cite[Theorem~3.1]{Zhivotovskiy})
\end{remark}

This result combines with the continuity of Lemma~\ref{lem:D-cts} to also
yield the following:

\begin{corollary}
\label{cor4.7}
If $X$ has a density with respect to Lebesgue measure,
$\mathbb{E}\|X\|<\infty $, $\#\mathrm{supp}(X)\ge k$, and
$0 < \alpha < \alpha _{0}(X)$, then for any fixed
$\mathcal{C}_{0}\in \mathscr{C}_{k}$ we have
%
\begin{equation}
\max _{\bar{\mathcal{C}}_{n}\in \bar{\mathscr{C}}_{n,k,\alpha n}} D_{X}(
\bar{\mathcal{C}}_{n}\,|\,\mathcal{C}_{0}) \to \inf _{\mathcal{C}\in
\mathscr{C}_{k}} D_{X}(\mathcal{C}\,|\,\mathcal{C}_{0})
\label{eq9}
\end{equation}
almost surely as $n\to \infty $.
\end{corollary}

\section{Empirical-level convergence of excess distortion}
\label{sec:distortion}

In this section we consider whether the (population-level) excess distortion
of the empirical $k$-means cluster centers converges to the infimum of
the (population-level) excess distortion in case (ii). We will show that
this guess can fail (Proposition~\ref{ex:inconsistency_two_sided_imbalanced}),
but also that imposing a certain logarithmic balance constraint recovers
such convergence. We let $X$ denote a random variable such that (ii) holds
throughout this section; as a reminder, this means that the optimal excess
distortion is finite but there exists no set of optimal $k$-means cluster
centers. We emphasize that, while convergence of the distortion is well-understood
in many settings
\cite{BartlettLinderLugosi, BiauDevroyeLugosi, Levrard, Zhivotovskiy},
no literature so far has considered the case of
$\mathbb{E}\|X\|^{2}= \infty $, let alone
$\mathscr{C}_{k}(X) = \varnothing $.

We begin with a negative result which is similar in spirit to that of Example~\ref{ex:inconsistency_two_sided}.
Precisely, we consider the following asymmetric version of the two-sided
$\textnormal{Par}(2)$ distribution.

\begin{example}%
\label{ex:population_level_imbalanced_2_sided}
For $0 \le p \le 1$, let $X$ be a real-valued random variable which satisfies
\begin{align*}
\mathbb{P}(X\geq t) = pt^{-2} &\qquad \textnormal{for all} \qquad t
\geq 1,
\\
\mathbb{P}(X\leq -t) = (1-p)t^{-2} &\qquad \textnormal{for all}
\qquad t\geq 1.
\end{align*}
(In other words, $X=SZ$ where $Z$ and $S$ are independent,
$Z\sim \textnormal{Par}(2)$, and $S= +1$ with probability $p$ and
$S=-1$ with probability $1-p$.) For $0 \le p < \sfrac{1}{4}$ or
$\sfrac{3}{4}<p\le 1$, we have $\mathscr{C}_{2}(X)=\varnothing $ and with
$\mathcal{C}_{0}=\{0\}$, we have
\begin{equation*}
\inf _{\mathcal{C}\in \mathscr{C}_{2}} D_{X}(\mathcal{C}\,|\,
\mathcal{C}_{0}) = -4\left ( \max (p,1-p)+(2p-1)^{2} \right ).
\end{equation*}
\end{example}

\begin{proof}
As in Example~\ref{ex:2-sided-pareto-pop}, it suffices to consider Voronoi
regions of the form $\mathscr{V}(r)=\{(-\infty ,r),\break [r,\infty )\}$, which
are generated for instance by the cluster
$\mathcal{C}(r):=\{r-1,r+1\}$. For convenience, we denote by
\begin{equation*}
D(r):=-\sum _{\mathcal{V}\in \mathscr{V}(r)} \mathbb{P}(X\in
\mathcal{V}) m_{X}(\mathcal{V})^{2},
\end{equation*}
the corresponding excess distortion with respect to the cluster center
$\mathcal{C}_{0}=\{0\}$. Using some of the computations from Example~\ref{ex:2-sided-pareto-pop},
by direct calculation, we find
%
\begin{equation}
D(r)=
\begin{cases}
-4p - 4(2p-1)^{2}
\frac{\left ( 1-\frac{p}{(2p-1)r} \right )^{2}}{1-pr^{-2}} &r\geq 1,
\\
-4(1-p) - 4(2p-1)^{2}
\frac{\left ( 1-\frac{1-p}{(2p-1)r} \right )^{2}}{1-(1-p)r^{-2}} &r
\leq -1.
\end{cases}
\label{eq10}
\end{equation}
As a result, we can use that the support of $X$ lies in
$(-\infty ,-1]\cup [1,\infty )$ to get:
\begin{equation*}
\inf _{\mathcal{C}\in \mathscr{C}_{2}} D_{X}(\mathcal{C}\,|\,
\mathcal{C}_{0}) = \inf _{r\in (-\infty ,-1]\cup [1,
\infty )}D(r) = -4\max \left ( 1,\max (p,1-p)+(2p-1)^{2} \right ).
\end{equation*}
Further, $\mathcal{C}_{2}(X)=\varnothing $ whenever the above infimum
$\inf _{r\in (-\infty ,-1]\cup [1,\infty )]}D(r)$ is not attained. This
holds for $p\in [0,1/4)\cup (3/4,1]$, for which we obtain the desired formula
for the infimum excess distortion.
\end{proof}

In the following result (which is similar to Example~\ref{ex:inconsistency_two_sided}),
we show that convergence in excess distortion can fail; its proof is somewhat
long and also appears in Appendix~\ref{sec:proofs}.

\begin{example}%
\label{ex:inconsistency_two_sided_imbalanced}
For $0 \le p < \sfrac{1}{4}$, let $X$ be a real-valued random variable
which satisfies
\begin{align*}
\mathbb{P}(X\geq t) = pt^{-2} &\qquad \textnormal{for all} \qquad t
\geq 1,
\\
\mathbb{P}(X\leq -t) = (1-p)t^{-2} &\qquad \textnormal{for all}
\qquad t\geq 1.
\end{align*}
(In other words, $X$ has the asymmetric $\textnormal{Par}(2)$ distribution
from Example~\ref{ex:population_level_imbalanced_2_sided}.) Then, there
exists a universal constant $c>0$ such that almost surely, for infinitely
many $n\in \mathbb{N}$, for any
$\bar{\mathcal{C}}_{n}\in \bar{\mathscr{C}}_{n,2}$ we have
\begin{equation*}
\max (\bar{\mathcal{C}}_{n})\geq c\frac{\sqrt{pn}}{\log n}.
\end{equation*}
Consequently, we have
\begin{equation*}
\limsup _{n\to \infty} \min _{\bar{\mathcal{C}}_{n}\in
\bar{\mathscr{C}}_{n,2}} D_{X}(\bar{\mathcal{C}}_{n}\,|\, \mathcal{C}_{0})
\geq \inf _{\mathcal{C}\in \mathscr{C}_{2}} D_{X}(\mathcal{C}\,|\,
\mathcal{C}_{0}) +4(1-2p).
\end{equation*}
almost surely.
\end{example}

Next we turn to positive results, which will essentially be consequences
of the following uniform convergence result for the distortion on any Voronoi
partition. The proof of this result uses the fact that any Voronoi region
in a Voronoi partition into $k$ pieces can be written as an intersection
of at most $k-1$ half-spaces, which have finite VC dimension
$O(dk\log k)$ \cite{blumer1989learnability}. We recall the definition of
$\bar{\mathscr{P}}_{n,k,\gamma}$ given in Definition~\ref{def:balance-k-means-emp}.
The proof appears in Appendix~\ref{sec:proofs}.

\begin{lemma}%
\label{lemma:uniform_convergence_voronoi}
There exists a universal constant $c>0$ such that the following holds.
If $X$ satisfies $\mathbb{E}\|X\|<\infty $ and
$\inf _{\mathcal{C}\in \mathscr{C}_{k}}D_{X}(\mathcal{C}\,|\,
\mathcal{C}_{0})>-\infty $ for some (equivalently all)
$\mathcal{C}_{0}\in \mathscr{C}_{k}$, then there exists a constant
$C>0$, depending only on the distribution of $X$, such that for any
$n\in \mathbb{N}$, $k\geq 2$, $\delta \in (0,1]$, and
$\gamma \geq \frac{c}{\delta}(dk\log k\log n+\log \frac{1}{\delta})^{2}$,
we have
\begin{equation*}
\begin{split} \Delta (\mathscr{V}) &:=\max \left ( \left |\bar D_{n}(
\mathscr{V}) - D_{X}(\mathscr{V}) \right |, \sum _{\mathcal{V}\in
\mathscr{V}} \mu _{X}(\mathcal{V})\|\bar m_{n}(\mathcal{V})-m_{X}(
\mathcal{V})\|^{2} \right )
\\
&\qquad \leq Ck
\frac{\sqrt{dk\log k\log n + \log \frac{1}{\delta}}}{(\gamma \delta )^{1/4}},
\end{split}
\end{equation*}
for all $\mathscr{V}\in \bar{\mathscr P}_{n,k,\gamma}$, with probability
at least $1-\delta $. Further, suppose that
$\{\beta _{n}\}_{n\in \mathbb{N}}$ satisfies
$\sum _{n\in \mathbb{N}}(n\beta _{n})^{-1}<\infty $, and that
$\{\gamma _{n}\}_{n\in \mathbb{N}}$ satisfies
$\gamma _{n}\geq \beta _{n}(dk\log k \log n)^{2}$ for all
$n\in \mathbb{N}$; then, almost surely, there exists (random) $N$ such
that $n\ge N$ implies
\begin{equation*}
\Delta (\mathscr{V}) \leq C
\frac{k\beta _{n}^{1/4}\sqrt{\log n}}{\gamma _{n}^{1/4}}
\end{equation*}
for all $\mathscr{V}\in \bar{\mathscr P}_{n,k,\gamma _{n}}$.
\end{lemma}

Now we can state our main positive result in this section. Recall that
we are assuming we are in case (i) or (ii), meaning
$\inf _{\mathcal{C}\in \mathscr{C}_{k}} D_{X}(\mathcal{C}\,|\,
\mathcal{C}_{0})>-\infty $.

\begin{theorem}%
\label{thm:consistency_log_points_per_cluster}
Suppose that $\mathbb{E}\|X\|<\infty $ and that
$\#\mathrm{supp}(X)\ge k$. Then, for any sequence
$\{\gamma _{n}\}_{n\in \mathbb{N}}$ such that
$(\log n)^{2}\ll \gamma _{n}\ll n$ as $n\to \infty $, we have
\begin{equation*}
\max _{\bar{\mathcal{C}}_{n}\in \bar{\mathscr{C}}_{n,k,\gamma _{n}}} D_{X}(
\bar{\mathcal{C}}_{n}\,|\,\mathcal{C}_{0}) \overset{\mathbb{P}}{\to}
\inf _{\mathcal{C}\in \mathscr{C}_{k}} D_{X}(\mathcal{C}\,|\,
\mathcal{C}_{0}).
\end{equation*}
Further, for any sequence $\{\gamma _{n}\}_{n\in \mathbb{N}}$ such that
$(\log n)^{4}\le \gamma _{n}\ll n$, we have
\begin{equation*}
\max _{\bar{\mathcal{C}}_{n}\in \bar{\mathscr{C}}_{n,k,\gamma _{n}}} D_{X}(
\bar{\mathcal{C}}_{n}\,|\,\mathcal{C}_{0}) \to \inf _{\mathcal{C}\in
\mathscr{C}_{k}} D_{X}(\mathcal{C}\,|\,\mathcal{C}_{0}).
\end{equation*}
almost surely, for each $\mathcal{C}_{0}\in \mathscr{C}_{k}$.
\end{theorem}

\begin{proof}
Let $\mathcal{C}_{0}=\{0\}$. Fix a sequence
$\{\gamma _{n}\}_{n\in \mathbb{N}}$ such that $\gamma _{n}=o(n)$. We define
the event
\begin{equation*}
\mathcal{E}:=\left \{ \lim _{n\to \infty } \inf _{\mathscr{V} \in
\bar{\mathscr P}_{n,k,\gamma _{n}}} D_{X}(\mathscr V) =\inf _{
\mathcal{C}\in \mathscr{C}_{k}} D_{X}(\mathcal{C}\,|\,\mathcal{C}_{0})
\right \}.
\end{equation*}
We briefly argue that $\mathcal{E}$ has probability one. First, from Lemma~\ref{lemma:equivalence_distortions},
the right-hand side in the definition of $\mathcal{E}$ is equal to the
infimum of the excess distortion $D_{X}(\mathscr{V})$ over Voronoi partitions
$\mathscr{V}$ with $1\le \#\mathscr{V}\le k$. Then, fix any sequence of
partitions $\mathscr{V}_{\ell} $ for $\ell \geq 1$ with
$1\le \#\mathscr{V}_{\ell}\le k$ such that
$D_{X}(\mathscr{V}_{\ell})\to \inf _{\mathcal{C}\in \mathscr{C}_{k}} D_{X}(
\mathcal{C}\,|\,\mathcal{C}_{0})$ as $\ell \to \infty $. Up to deleting
some centroids, we may also assume without loss of generality that, for
each $\ell \geq 1$, each region of $\mathscr{V}_{\ell}$ satisfies
$p_{\min}(X,\mathscr{V}_{\ell}) > 0$. Then, by the strong law of large
numbers, with probability one, for each $\ell \geq 1$ we have we have
$\liminf _{n\to \infty}\min _{\mathcal{V}\in \mathscr{V}_{\ell}}
\bar\mu _{n}(\mathscr{V}_{\ell})>0$. Therefore, since
$\gamma _{n}=o(n)$, for any $\ell \geq 1$, we have
$\mathscr{V}_{\ell}\in \bar{\mathscr P}_{n,k,\gamma _{n}}$ for $n$ sufficiently
large. Under this event of probability one, the event $\mathcal{E}$ holds.

For convenience, we now define
\begin{equation*}
\Delta _{\max}(n):=\sup _{\mathscr{V}\in \bar{\mathscr P}_{n,k,
\gamma _{n}}} \Delta (\mathscr{V}),
\end{equation*}
where $\Delta (\mathscr{V})$ is defined in Lemma~\ref{lemma:uniform_convergence_voronoi}.
We denote $\mathcal{F}=\{\lim _{n\to \infty}\Delta _{\max}(n)=0\}$. Now
consider any cluster
$\mathcal{C}_{n}\in \bar{\mathscr{C}}_{n,k,\gamma _{n}}$ and denote by
$\bar{\mathscr V}_{n}$ a $\gamma _{n}$-balanced $k$-means clustering Voronoi
partition for which
$\bar{\mathcal{C}}_{n}=\bar m_{n}(\bar{\mathscr V}_{n})$. Then, we have
\begin{align*}
D_{X}(\bar{\mathscr{V}}_{n})-\Delta _{\max}(n) \leq \bar D_{n}(
\bar{\mathscr{V}}_{n}) = \inf _{\mathscr{V} \in \bar{\mathscr P}_{n,k,
\gamma _{n}}} \bar D_{n}(\mathscr{V}) \leq \inf _{\mathscr{V} \in
\bar{\mathscr P}_{n,k,\gamma _{n}}} D_{X}(\mathscr{V}) +\Delta _{\max}(n).
\end{align*}
Since $\bar{\mathscr V}_{n}\in \bar{\mathscr P}_{n,k,\gamma _{n}}$, we
obtained
%
\begin{equation}
\label{eq:bound_distortion_emp}
\left |D_{X}(\bar{\mathscr{V}}_{n}) -\inf _{\mathcal{C}\in
\mathscr{C}_{k}} D_{X}(\mathcal{C}\,|\,\mathcal{C}_{0})\right | \leq 2
\Delta _{\max}(n) + \left |\inf _{\mathscr{V} \in \bar{\mathscr P}_{n,k,
\gamma _{n}}} D_{X}(\mathscr{V}) - \inf _{\mathcal{C}\in \mathscr{C}_{k}}
D_{X}(\mathcal{C}\,|\,\mathcal{C}_{0})\right |.
\end{equation}
Next, note that
\begin{align*}
D_{X}(\bar{\mathcal{C}}_{n}\,|\,\mathcal{C}_{0}) &= \mathbb{E}\left [
\min _{c\in \bar m_{n}(\bar{\mathscr{V}}_{n})} \|c-X\|^{2}-\|X\|^{2}
\right ]
\\
&\leq \mathbb{E}\left [ \sum _{\mathcal{V}\in \bar{\mathscr V}_{n}}
\left ( \|\bar m_{n}(\mathcal{V})-X\|^{2}-\|X\|^{2} \right )
\mathbbold{1}\{X\in \mathcal{V}\} \right ]
\\
&=\mathbb{E}\left [ \sum _{\mathcal{V}\in \bar{\mathscr V}_{n}}
\left ( \|m_{X}(\mathcal{V})-X\|^{2}-\|X\|^{2} \right )\mathbbold{1}
\{X\in \mathcal{V}\} \right ]
\\
&\qquad +\sum _{\mathcal{V}\in \bar{\mathscr V}_{n}} \mathbb{P}(X\in
\mathcal{V}) \|\bar m_{n}(\mathcal{V})-m_{X}(\mathcal{V})\|^{2}
\\
&\leq D_{X}(\bar{\mathscr{V}}_{n}) + \Delta _{\max}(n).
\end{align*}
Together with equation~\eqref{eq:bound_distortion_emp}, we showed that
%
\begin{equation}
\label{eq:final_bound_deviation}
\max _{\bar{\mathcal{C}}_{n}\in \bar{\mathscr{C}}_{n,k,\gamma _{n}}} D_{X}(
\bar{\mathcal{C}}_{n}\,|\,\mathcal{C}_{0}) -\inf _{\mathcal{C}\in
\mathscr{C}_{k}} D_{X}(\mathcal{C}\,|\,\mathcal{C}_{0}) \leq 3\Delta _{
\max}(n) + \left |\inf _{\mathscr{V} \in \bar{\mathscr P}_{n,k,
\gamma _{n}}} D_{X}(\mathscr{V}) - \inf _{\mathcal{C}\in \mathscr{C}_{k}}
D_{X}(\mathcal{C}\,|\,\mathcal{C}_{0})\right |
\end{equation}
On $\mathcal{E}\cap \mathcal{F}$, the right-hand side converges to
$0$ as $n\to \infty $. We recall that $\mathcal{E}$ has probability one.
Also, from Lemma~\ref{lemma:uniform_convergence_voronoi}, we have
$\mathbb{P}(\mathcal{F})=1$ whenever $\gamma _{n}\geq (\log n)^{4}$, which
gives the desired almost-sure result.

For the result on convergence in probability, Lemma~\ref{lemma:uniform_convergence_voronoi}
shows that whenever $\gamma _{n}\gg (\log n)^{2}$, we have
$\Delta _{\max}(n)\to 0$ in probability as $n\to \infty $. Hence, since
$\mathcal{E}$ has probability one, the right-hand side of equation~\eqref{eq:final_bound_deviation}
converges in probability to $0$. This ends the proof.
\end{proof}

The proof of Theorem~\ref{thm:consistency_log_points_per_cluster} reveals
that the condition $\gamma \ge (\log n)^{4}$ guaranteeing almost sure convergence
can be reduced at least to
$\gamma _{n}\ge (\log n)^{3}(\log \log n)^{2}$. It is not clear what to
expect the sharpest possible rate for $\gamma _{n}$ to be, but this would
be an interesting question for future work.

Lastly, we note that Theorem~\ref{thm:cluster-consistency} and Theorem~\ref{thm:consistency_log_points_per_cluster}
can, in some sense, be unified; more precisely, the following result shows
that a mild (i.e., logarithmic) balance constraint guarantees that the
empirical $k$-means cluster centers converge to some set of population
$k'$-means cluster centers, for some $1\le k'\le k$, after discarding all
diverging centers.

\begin{theorem}%
\label{thm:consistency_unified}
Suppose that $\mathbb{E}\|X\|<\infty $ and
$\inf _{\mathcal{C}\in \mathscr{C}_{k}}D_{X}(\mathcal{C}\,|\,
\mathcal{C}_{0})>-\infty $. For any sequence
$\{\gamma _{n}\}_{n\in \mathbb{N}}$ such that
$(\log n)^{4}\le \gamma _{n} \ll n$, the following holds almost surely:
For any sequence $\{\bar{\mathcal{C}}_{n}\}_{n\geq 1}$ such that
$\bar{\mathcal{C}}_{n}\in \bar{\mathscr{C}}_{n,k,\gamma _{n}}$ for all
$n\in \mathbb{N}$, we can decompose
$\bar{\mathcal{C}}_{n} = \mathcal{A}_{n}\sqcup \mathcal{B}_{n}$ so that
$\mathcal{A}_{n}\neq \varnothing $ and
\begin{equation*}
\min _{1\le k'\le k}\min _{\mathcal{C}\in \mathscr{C}_{k'}(X)}d_{
\textnormal{H}}( \mathcal{A}_{n},\mathcal{C})\to 0 \quad \text{and}
\quad \min _{c\in \mathcal{B}_{n}}\|c\|\to \infty
\end{equation*}
as $n\to \infty $.
\end{theorem}

Most of the heavy-lifting is done in the following result, which we remark
is completely deterministic. It states that, among any sequence of sets
of cluster centers which approximate the optimal excess distortion, the
collection of all bounded centers must converge to a set of population-level
$k'$-means cluster centers, for some $1\le k'\le k$.

\begin{lemma}%
\label{lemma:can_delete_diverging_clusters}
Suppose that $\mathbb{E}\|X\|<\infty $ and
$\inf _{\mathcal{C}\in \mathscr{C}_{k}}D_{X}(\mathcal{C}\,|\,
\mathcal{C}_{0})>-\infty $. Consider a sequence of clusters
$\mathcal{C}_{n}\in \mathscr{C}_{k}$ for $n\in \mathbb{N}$ such that
\begin{equation*}
\lim _{n\to \infty}D_{X}(\mathcal{C}_{n}\,|\,\mathcal{C}_{0})= \inf _{
\mathcal{C}\in \mathscr{C}_{k}} D_{X}(\mathcal{C}\,|\,\mathcal{C}_{0}).
\end{equation*}
Then, there exists a decomposition
$\mathcal{C}_{n}=\mathcal{A}_{n}\sqcup \mathcal{B}_{n}$ for
$n\in \mathbb{N}$ so that $\mathcal{A}_{n}\neq \varnothing $ and
\begin{equation*}
\min _{1\le k'\le k}\min _{\mathcal{C}\in \mathscr{C}_{k'}(X)}d_{
\textnormal{H}}(\mathcal{A}_{n},\mathcal{C})\to 0 \quad \text{and}
\quad \min _{c\in \mathcal{B}_{n}}\|c\|\to \infty
\end{equation*}
as $n\to \infty $.
\end{lemma}

We are now ready to prove our unified convergence result.

\begin{proof}[Proof of Theorem~\ref{thm:consistency_unified}]
If the sequence $\{\gamma _{n}\}_{n\in \mathbb{N}}$ satisfies the given
assumptions, then Theorem~\ref{thm:consistency_log_points_per_cluster}
implies that the event
\begin{equation*}
\mathcal{E}:=\left \{ \lim _{n\to \infty } \max _{\bar {\mathcal{C}}_{n}
\in \bar {\mathscr{C}}_{n,k,\gamma _{n}}} D_{X}(\bar {\mathcal{C}}_{n}
\,|\,\mathcal{C}_{0}) =\inf _{\mathcal{C}\in \mathscr{C}_{k}} D_{X}(
\mathcal{C}\,|\,\mathcal{C}_{0}) \right \}
\end{equation*}
has probability one. Then, on $\mathcal{E}$, from Lemma~\ref{lemma:can_delete_diverging_clusters}
we can decompose any
$\bar{\mathcal{C}}_{n}\in \bar{\mathscr{C}}_{n,k,\gamma _{n}}$ into
$\mathcal{C}_{n}=\mathcal{A}_{n}\sqcup \mathcal{B}_{n}$ to satisfy the
desired result.
\end{proof}

Although the conclusion of Theorem~\ref{thm:consistency_unified} is, strictly
speaking, weaker than the conclusion of Theorem~\ref{thm:cluster-consistency},
we find that it is advantageous in the sense that the balance constraint
is weaker and hence will be non-binding in many settings.


\appendix

\section{Computation for Figure~\ref{fig:two_sided_pareto}}
\label{sec:computation}

The code used to produce Figure~\ref{fig:two_sided_pareto} (i.e., the simulation
of Example~\ref{ex:inconsistency_two_sided}) is publicly available at
\url{https://github.com/aqjaffe/kmeans_clustering_inconsistency}. In this
appendix we briefly comment on the computational of this simulation, which
requires computing the $k$-means cluster centers exactly, and not by standard
iterative procedures like the Lloyd-Forgy algorithm or the Hartigan algorithm;
towards this end, we show that some special structure in the case of dimension
$d=1$ allows us to efficiently compute the exact cluster centers.

We propose a procedure in Algorithm~\ref{alg:clustering}. To see that it
correctly produces the optimal set of $k$-means cluster centers, simply
use Lemma~\ref{lem:centroid-condition} to see that the for loop (lines
8-14) must at some point contain the optimal cluster centers, and use equation~\eqref{eqn:emp-excess-dist}
to see that the empirical excess distortion $D$ can indeed be written (lines
7 and 11) as a certain weighted sum of the putative cluster centers
$a$ and $b$. The time complexity of this procedure is linear in the number
of samples $n$, plus the time complexity of sorting (line 4), i.e.,
$O(n+n\log n)$. Although this procedure is certainly well-known and it
is easy to derive in an elementary way, we are not aware of it appearing
in previous literature; the closely-related works
\cite{WangSong, Fisher} study dynamic programming approaches to the empirical
$k$-means clustering for general $k\in \mathbb{N}$ in dimension
$d=1$, although they establish the worse rate $O(n^{2}k+n\log n)$.

\begin{algorithm}[h]
		\caption{Computing $k$-means clustering for $k=2$ in dimension $d=1$}%
\label{alg:clustering}

		\begin{algorithmic}[1]

			\State \textbf{input:} samples $X_{1},\ldots , X_{n}$ in $\mathbb{R}$
			\State \textbf{output:} optimal set of $2$-means cluster centers $(a^{\ast},b^{\ast})\in \bar{\mathscr{C}}_{n,2}$

			\State \,

			\State $Y_{1},\ldots ,Y_{n} \leftarrow \textbf{sort\_increasing}(X_{1},\ldots , X_{n})$
			\State $a^{\ast}\leftarrow Y_{1}$
			\State $b^{\ast}\leftarrow \sum _{i=2}^{n}Y_{i}/(n-1)$
			\State $D^{\ast}\leftarrow -(|a^{\ast}|^{2}+(n-1)|b^{\ast}|^{2})/n$
			\For{$k=2,\ldots , n-1$}
				\State $a\leftarrow ((k-1)a+Y_{k})/k$
				\State $b\leftarrow ((n-k+1)b-Y_{k})/(n-k)$
				\State $D\leftarrow -(k|a|^{2}+(n-k)|b|^{2})/n$
				\If{$D \le D^{\ast}$}
					\State $D^{\ast} \leftarrow D$
					\State $(a^{\ast},b^{\ast})\leftarrow (a,b)$
				\EndIf
			\EndFor

		\end{algorithmic}
\end{algorithm}

In earlier versions of the simulation, we computed the $k$-means cluster
centers with standard implementations of the Lloyd-Forgy algorithm, and
we found that outlying samples were occasionally not placed into a cluster
alone; we believe this reflects some implicit regularization of early stopping
in the algorithm. As such, we find that the Lloyd-Forgy algorithm suffers
from some limitations for samples from heavy-tailed distributions; we find
this to be in the same spirit the works
\cite{LloydHighDimI,LloydHighDimII} which identify some other limitations
for samples in high dimensions.

\section{Discussion of case (iii)}
\label{sec:case-iii}

In this appendix we give a brief discussion of case (iii) in which problem~\eqref{eqn:def-pop-k-means-minimal}
has infinite optimal excess distortion hence no solution.

We begin by giving a characterization of when case (iii) occurs; the answer
is that it occurs if and only if $\|X\|$ is not stochastically dominated
by a scalar multiple of $\textnormal{Par}(2)$.

Towards this end we require the following, which we will also use in some
later proofs.

\begin{lemma}%
\label{lemma:computation_bound_norm}
Let $X$ be a random variable such that $\|X\|$ is stochastically dominated
by $hZ$ where $h>0$ and $Z\sim \textnormal{Par}(2)$. Then, for any measurable
set $A\subseteq \mathbb{R}^{d}$, we have
\begin{equation*}
\|\mathbb{E}[X\mid X\in A] \| \leq \mathbb{E}[\|X\|\mid X\in A] \leq
\frac{2h}{\sqrt{\mathbb{P}(X\in A)}}.
\end{equation*}
\end{lemma}

\begin{proof}
The first inequality is immediate by convexity. The second is true if
$\mathbb{P}(X\in A) = 0$ since then there is nothing to prove, so let us
assume $\mathbb{P}(X\in A) > 0$. To do this, let $X$ and $Z$ be coupled
on the same probability space in such a way that we have
$\|X\|\le hZ$ almost surely, so that
\begin{equation*}
\mathbb{P}(X\in A)\cdot \mathbb{E}[\|X\|\mid X\in A] =\mathbb{E}[\|X
\|\cdot \mathbbold{1}\{ X\in A\}] \le h\mathbb{E}[Z\cdot
\mathbbold{1}\{ X\in A\}].
\end{equation*}
Now consider maximizing the random variable
$Z\cdot \mathbbold{1}_{\mathcal{E}}$ over all
$\mathcal{E}\in \sigma (X,Z)$ with $\mathbb{P}(\mathcal{E}) = r$, for fixed
$0<r<1$. It is clearly maximized when $Z$ takes on its largest possible
value, i.e., $\mathcal{E}= \{Z\ge F_{c}^{-1}(r)\}$ where
$F_{c}(x)=\mathbb{P}(Z\geq x)=x^{-2}$ is the complementary cumulative distribution
function of $Z$. Thus we take $r:=\mathbb{P}(X\in A)$ and get
\begin{align*}
\mathbb{P}(X\in A)\cdot \mathbb{E}[\|X\|\mid X\in A] &\leq h
\mathbb{E}[Z \cdot \mathbbold{1}\{X\in A\}]
\\
&\leq h\mathbb{E}[Z \cdot \mathbbold{1}\{Z \geq F_{c}^{-1}(r)\}]
\\
&=h\mathbb{E}\left [ Z \cdot \mathbbold{1}\{Z \geq r^{-1/2}\} \right ]
= 2h\sqrt{r}
\end{align*}
This ends the proof.
\end{proof}

Now we have the following result characterizing when we have case (iii).

\begin{lemma}%
\label{lemma:charact_finite_obj}
For $k\geq 2$ and $\mathbb{E}\|X\|<\infty $, the following are equivalent:
\begin{enumerate}
\item[(a)] For some (equivalently, all)
$\mathcal{C}_{0}\in \mathscr{C}_{k}$, we have
$\inf _{\mathcal{C}\in \mathscr{C}_{k}} D_{X}(\mathcal{C}\,|\,
\mathcal{C}_{0}) >-\infty $.
\item[(b)] There exists a constant $r>0$ such that we have
$\mathbb{P}(\|X\|\geq t) \leq rt^{-2}$ for all $t\ge 0$
\end{enumerate}
Consequently, if $\mathbb{E}\|X\|^{p}=\infty $ for some $1\le p < 2$, then
$\inf _{\mathcal{C}\in \mathscr{C}_{k}} D_{X}(\mathcal{C}\,|\,
\mathcal{C}_{0}) =-\infty $ for all
$\mathcal{C}_{0}\in \mathscr{C}_{k}$.
\end{lemma}

\begin{proof}
We start by proving that (a) implies (b). Fix
$\mathcal{C}_{0}=\{0\}$ and consider a random variable $X$ for which property
(a) holds. For convenience, denote
\begin{equation*}
M := -\inf _{\mathcal{C}\in \mathscr{C}_{k}} D_{X}(\mathcal{C}\,|\,
\mathcal{C}_{0})<\infty ,
\end{equation*}
and note that $D_{X}(\mathcal{C}_{0}\,|\, \mathcal{C}_{0})=0$, hence
$M\geq 0$. Now for any vector $v\in \mathbb{R}^{d}$ with $\|v\|=1$ and
any $t\geq 0$, we consider the set of cluster centers
$\mathcal{C}(v,t):=\{0,t v\}$, and we note that the corresponding Voronoi
partition is given by the pair of halfspaces
$\{x\in \mathbb{R}^{d}: v^{\top }X<t/2\}$ and
$\{x\in \mathbb{R}^{d}:v^{\top }X\geq t/2\}$. Consequently, we can compute
the distortion as follows:
\begin{align*}
D_{X}(\mathcal{C}(v,t) \,|\,\mathcal{C}_{0}) &= \mathbb{E}\left [ (\|X-tv
\|^{2}-\|X\|^{2})\mathbbold{1}\{v^{\top }X>t/2\} \right ]
\\
&= \mathbb{E}\left [ (\|tv\|^{2}-2tv^{\top }X)\mathbbold{1}\{v^{\top }X>t/2
\} \right ]
\\
&= t\mathbb{E}\left [ (t-2v^{\top }X)\mathbbold{1}\{v^{\top }X>t/2\}
\right ]
\\
&\leq -t^{2}\mathbb{P}(v^{\top }X \geq t).
\end{align*}
Rearranging this and using the definition of $M$, we get
\begin{equation*}
\mathbb{P}(v^{\top }X\geq t) \leq \frac{M}{t^{2}}.
\end{equation*}
By applying this bound to all the standard basis vectors in
$\mathbb{R}^{d}$, we get
\begin{equation*}
\mathbb{P}(\|X\|\geq t) \leq \mathbb{P}\left ( \|X\|_{\infty }\geq
\frac{t}{\sqrt d} \right ) \leq \frac{2d^{2}M}{t^{2}}.
\end{equation*}
for all $t\ge 0$. This proves the claim, with $r=2d^{2}M$.

Next we show that (b) implies (a). Let $X$ be any random variable satisfying
condition (b). Fix $\mathcal{C}_{0}=\{0\}$ and any set of points
$\mathcal{C}\in \mathscr{C}_{k}$, and compute:
\begin{equation}
\label{eq:lower_bound_distortion_1}
\begin{split} D_{X}(\mathcal{C}\,|\, \mathcal{C}_{0}) &=\sum _{c\in
\mathcal{C}} \mathbb{E}\left [ (\|c-X\|^{2} - \|X\|^{2})\cdot
\mathbbold{1}\{X\in \mathcal{V}_{\mathcal{C}}(c)\} \right ]
\notag
\\
&=\sum _{c\in \mathcal{C}} \mathbb{P}(X\in \mathcal{V}_{\mathcal{C}}(c))
\left ( \|c\|^{2} - 2c^{\top }m_{X}(\mathcal{V}_{\mathcal{C}}(c))
\right )
\notag
\\
&\geq - \sum _{c\in \mathcal{C}} \mathbb{P}(X\in \mathcal{V}_{
\mathcal{C}}(c)) \|m_{X}(\mathcal{V}_{\mathcal{C}}(c))\|^{2}.
\end{split}
\end{equation}
Next, fix any $c\in \mathcal{C}$ such that
$\mathbb{P}(X\in \mathcal{V}_{\mathcal{C}}(c))>0$. By condition (b), there
exists a constant $r>0$ such that the random variable $\|X\|$ is stochastically
dominated by $\sqrt{r}\cdot Z$ where $Z$ is a random variable defined via
$\mathbb{P}(Z\ge t) = t^{-2}$ for all $t\ge 1$. Therefore, Lemma~\ref{lemma:computation_bound_norm}
implies
\begin{align*}
\|m_{X}(\mathcal{V}_{\mathcal{C}}(c))\| \leq \mathbb{E}\left [\|X\|
\mid X\in \mathcal{V}_{\mathcal{C}}(c)\right ] \leq 2\sqrt{
\frac{r}{\mathbb{P}(X\in \mathcal{V}_{\mathcal{C}}(c))}},
\end{align*}
and this gives
\begin{equation*}
D_{X}(\mathcal{C}\,|\, \mathcal{C}_{0}) \geq -4kr.
\end{equation*}
Because this holds for all $\mathcal{C}\in \mathscr{C}_{k}$, we have shown
that condition (a) holds.

Lastly, we show that $\mathbb{E}\|X\|^{p} = \infty $ for some
$1\le p<2$ implies that conditions (a) and (b) both fail. To do this, we
observe the contrapositive, namely that condition (b) implies
$\mathbb{E}\|X\|^{p} < \infty $ for all $1\le p<2$. This follows immediately
from the tail formula, since
\begin{equation*}
\mathbb{E}\|X\|^{p} = \int _{0}^{\infty }\mathbb{P}(\|X\|^{p} \geq t)
\leq 1+\int _{1}^{\infty }\frac{r}{t^{2/p}} dt =1+ \frac{rp}{2-p}<
\infty .
\end{equation*}
This completes the proof.
\end{proof}

Next, we pose the question of what form of consistency may be desired in
case (iii). A natural guess would be to require that the empirical
$k$-means cluster centers achieve a (population) excess distortion which
diverges. However, we will illustrate in the following example that this
notion is actually quite weak.

Indeed, consider the following ``clustering procedure'' that, given samples
$X_{1},\ldots ,X_{n}$, outputs the cluster centers
$\{0,\bar v_{n}\}$, where $\bar v_{n}\in \mathbb{R}^{d}$ is any vector
that minimizes the sample variance of the subset of the data falling in
the region $H_{v}:=\{x\in \mathbb{R}^{d}: v^{\top }x\ge \|v\|^{2}\}$, provided
that $H_{v}$ does not contain too few points. More precisely, for some
integer $\beta _{n}\ge 0$, we let
\begin{equation*}
\bar v_{n}\in \arg \max \left \{ \|v\|^{2} \sum _{i=1}^{n}
\mathbbold{1}\{v^{\top }X_{i}\geq \|v\|^{2}\} : v\in \mathbb{R}^{d},
\, \sum _{i=1}^{n} \mathbbold{1}\{v^{\top }X_{i}\geq \|v\|^{2}\}
\geq \beta _{n} \right \},
\end{equation*}
and the following result shows $\{0,\bar v_{n}\}$ achieves the putative
form of consistency in case (iii). This is not a meaningful set of cluster
centers, however; it always contains the cluster center $0$, and it yields
two clusters even for $k\ge 3$.

\begin{proposition}%
\label{prop:naive_cluster_estimator}
If $X$ satisfies $\mathbb{E}\|X\|<\infty $ and
$\inf _{\mathcal{C}\in \mathscr{C}_{k}}D_{X}(\mathcal{C}\,|\,
\mathcal{C}_{0})=-\infty $ for some $k\in \mathbb{N}$, and if
$\log n \ll \beta _{n}\ll n$, then we have
\begin{equation*}
\lim _{n\to \infty} D_{X}(\{0, \bar v_{n}\} \,|\,\mathcal{C}_{0})=-
\infty
\end{equation*}
almost surely.
\end{proposition}

\begin{proof}
We first give an upper bound on the excess distortion of clusters of the
form $\{0,v\}$ for $v\in \mathbb{R}^{d}$:
%
\begin{equation}
\label{eq:upper_bound_excess_distortion_naive}
\begin{split} D_{X}(\{0,v\} \,|\,\mathcal{C}_{0}) &= \mathbb{E}\left [ (
\|v\|^{2}-2v^{\top }X)\mathbbold{1}\{\|X-v\|<\|X\|\} \right ]
\\
&\leq -\|v\|^{2} \mathbb{P}(v^{\top }X\geq \|v\|^{2})=:-f(v).
\end{split}
%
\end{equation}
Next, let $\mathcal{H}$ be the set of halfspaces in $\mathbb{R}^{d}$. Since
linear separators in $\mathbb{R}^{d}$ have VC dimension $d+1$, a classical
uniform convergence theorem for VC classes (e.g.,
\cite[Theorem 5.1]{boucheron2005theory}) implies that there is a constant
$c_{1}>0$ for which the following holds. On an event
$\mathcal{E}_{n}$ of probability at least $1-n^{-2}$, we have
\begin{equation*}
\min \left \{\frac{1}{n}\sum _{i=1}^{n}\mathbbold{1}\{X_{i}\in H\},
\mathbb{P}(X\in H)\right \} \geq c_{1} \frac{d\log n}{n}
\end{equation*}
implies
\begin{equation*}
\frac{1}{2}\mathbb{P}(X\in H) \le \frac{1}{n}\sum _{i=1}^{n}
\mathbbold{1}\{X_{i}\in H\}\le 2\mathbb{P}(X\in H)
\end{equation*}
for any $H\in \mathcal{H}$. Now use $\log n = o(\beta _{n})$ to get
$n$ sufficiently large so that $\beta _{n}\geq c_{1} d\log n$, and suppose
that $\mathcal{E}_{n}$ holds. By equation~\eqref{eq:upper_bound_excess_distortion_naive}
and the definition of $\mathcal{E}_{n}$, we have
\begin{equation*}
-D_{X}(\{0,\bar v_{n}\} \,|\,\mathcal{C}_{0})\ge f(\bar v_{n})\geq
\frac{\|\bar v_{n}\|^{2}}{2} \sum _{i=1}^{n} \mathbbold{1}\{\bar v_{n}^{
\top }X_{i}\geq \|\bar v_{n}\|^{2}\}.
\end{equation*}
Also, the VC bound, the definition of $\bar v_{n}$, and the definition
of $\mathcal{E}_{n}$ imply
\begin{align*}
\|\bar v_{n}\|^{2}&\sum _{i=1}^{n} \mathbbold{1}\{\bar v_{n}^{\top }X_{i}
\geq \|\bar v_{n}\|^{2}\}
\\
&\geq \max \left \{ \|v\|^{2} \sum _{i=1}^{n}\mathbbold{1}\{v^{\top }X_{i}
\geq \|v\|^{2}\} : v\in \mathbb{R}^{d},\, \mathbb{P}(v^{\top }X\geq
\|v\|^{2}) \geq \frac{2\beta _{n}}{n} \right \}
\\
&\geq \frac{1}{2} \sup \left \{ \|v\|^{2} \mathbb{P}(v^{\top }X\geq
\|v\|^{2}) : v\in \mathbb{R}^{d},\, \mathbb{P}(v^{\top }X\geq \|v\|^{2})
\geq \frac{2\beta _{n}}{n} \right \}.
\end{align*}
Since
$\inf _{\mathcal{C}\in \mathscr{C}_{k}}D_{X}(\mathcal{C}\,|\,
\mathcal{C}_{0})=-\infty $, the proof of Lemma~\ref{lemma:charact_finite_obj}
precisely shows that
\begin{equation*}
\sup \left \{ \|v\|^{2} \mathbb{P}(v^{\top }X\geq \|v\|^{2}) : v\in
\mathbb{R}^{d} \right \} = \infty .
\end{equation*}
In other words, the function
\begin{equation*}
M(\epsilon ):= \sup \left \{ \|v\|^{2} \mathbb{P}(v^{\top }X\geq \|v
\|^{2}) : v\in \mathbb{R}^{d},\,\mathbb{P}(v^{\top }X\geq \|v\|^{2})
\geq \epsilon \right \},
\end{equation*}
satisfies $M(\epsilon )\to \infty $ as $\varepsilon \to 0$. Thus, combining
the two previous inequalities, we have shown
\begin{equation*}
\mathcal{E}_{n}\subseteq \left \{D_{X}(\{0,\bar v_{n}\} \,|\,
\mathcal{C}_{0}) \leq -\frac{1}{4}M\left ( \frac{2\beta _{n}}{n}
\right )\right \}
\end{equation*}
for sufficiently large $n$. Since
$\mathbb{P}(\mathcal{E}_{n}) \ge 1-n^{-2}$ by construction and since
$\beta _{n} = o(n)$ by assumption, Borel-Cantelli implies
$D_{X}(\{0,\bar v_{n}\} \,|\,\mathcal{C}_{0})\to -\infty $ as
$n\to \infty $ almost surely.
\end{proof}

\section{Additional results and proofs}
\label{sec:proofs}

In this appendix we collect some auxiliary results and proofs that were
deferred from the main body of the paper for the sake of the exposition.

\subsection{From Section~\ref{sec:centers}}
\label{appC.1}

Next, we prove the results from Section~\ref{sec:centers} on convergence
and non-convergence of the empirical cluster centers to the population
counterpart, in the well-defined case (i).

We begin with the proof of our flagship example, which was also visualized
in Figure~\ref{fig:two_sided_pareto}.

\begin{proof}[Proof of Example~\ref{ex:inconsistency_two_sided}]
We introduce some notation for convenience. First, write
$X_{(1)}^{n},\ldots , X_{(n)}^{n}$ for the order statistics of the samples
$X_{1},\ldots , X_{n}$. Then, for integers $1\le i<j\le n$, write
\begin{equation*}
S_{n}{[\,:i]} := \sum _{\ell =1}^{i}X_{(\ell )}^{n}, \qquad S_{n}[i:
\,] := \sum _{\ell =i+1}^{n}X_{(\ell )}^{n}, \qquad \textnormal{ and}
\qquad S_{n}[i:j]:=\sum _{\ell =i+1}^{j}X_{(\ell )}^{n}
\end{equation*}
for the partial sums of contiguous blocks of
$X_{(1)}^{n},\ldots , X_{(n)}^{n}$, and
\begin{equation*}
\bar{X}_{n}{[\,:i]} := \frac{1}{i}S_{n}{[\,:i]}, \qquad \bar{X}_{n}[i:
\,] := \frac{1}{n-i}S_{n}{[i:\,]}\qquad \textnormal{ and} \qquad
\bar{X}_{n}[i:j]:=\frac{1}{j-i}S_{n}[i:j]
\end{equation*}
for the corresponding partial averages. (Note that we use non-standard
conventions regarding inclusivity and exclusivity of the endpoints.) As
in the proofs of Example~\ref{ex:2-sided-pareto-pop} and Example~\ref{ex:1-sided-pareto-pop},
we apply Lemma~\ref{lem:centroid-condition} to see that
$\bar{\mathcal{C}}_{n}\in \bar{\mathscr{C}}_{n,2}$ is uniquely described
by its Voronoi regions; moreover, we may assume that the Voronoi regions
are simply
\begin{equation*}
\left \{\big(-\infty , X_{(i)}^{n}\big),\big[X_{(i)}^{n},\infty \big)
\right \},
\end{equation*}
for some integer $1\le i\le n$. Now define
$\bar{D}_{n,2}:\{1,\ldots , n\}\to \mathbb{R}$ via
\begin{align*}
\bar{D}_{n,2}(i) &:= \frac{1}{n}\sum _{\ell =1}^{i} \left (X_{(\ell )}^{n}-
\bar{X}_{n}[\,:i]\right )^{2} + \frac{1}{n}\sum _{\ell =i+1}^{n}
\left (X_{(\ell )}^{n}-\bar{X}_{n}[i:\,]\right )^{2}
\\
&= \frac{1}{n}\sum _{\ell =1}^{n} X_{\ell}^{2} - \frac{i}{n} \left (
\bar{X}_{n}[\,:i]\right )^{2} - \frac{n-i}{n} \left (\bar{X}_{n}[i:\,]
\right )^{2},
\end{align*}
which is simply the (absolute, not excess) empirical distortion of the
empirical $k$-means cluster centers when the Voronoi regions are as above.

The first step of the proof is to show
%
\begin{equation}
\label{eq:useful_comparison}
\left \{\bar{D}_{n,2}(i)>\bar{D}_{n,2}(j)\right \} = \left \{\bar{X}_{n}[j:
\,] - \bar{X}_{n}[i:j] > \sqrt{\frac{i(n-i)}{j(n-j)}} \left (\bar{X}_{n}[i:j]
- \bar{X}_{n}[\,:i]\right )\right \}
\end{equation}
for any $1\le i<j\le n$. Indeed, this follows by expanding
\begin{align*}
&n(\bar{D}_{n,2}(i) - \bar{D}_{n,2}(j) )
\\
&=\frac{1}{j}(S_{n}[\,:j])^{2} + \frac{1}{n-j}(S_{n}[j:\,])^{2} -
\frac{1}{i}(S_{n}[\,:i])^{2} - \frac{1}{n-i}(S_{n}[i:\,])^{2}
\\
&=\left (\frac{1}{j}-\frac{1}{i}\right )(S_{n}[\,:i])^{2} + \left (
\frac{1}{j}-\frac{1}{n-i}\right )(S_{n}[i:j])^{2} + \left (
\frac{1}{n-j}-\frac{1}{n-i}\right )(S_{n}[j:\,])^{2}
\\
& \qquad \qquad + \frac{2}{j} S_{n}[\,:i]S_{n}[i:j] - \frac{2}{n-i} S_{n}[i:j]S_{n}[j:
\,]
\\
&= (j-i) \frac{n-j}{n-i}\left ( \frac{1}{n-j}S_{n}[j:\,] -
\frac{1}{j-i}S_{n}[i:j]\right )^{2} - (j-i) \frac{i}{j}\left (
\frac{1}{j-i}S_{n}[i:j] - \frac{1}{i}S_{n}[\,:i]\right )^{2}.
\end{align*}
and rearranging.

Next, we introduce some events and bound their probability of failure.
To do this, we introduce the notation
\begin{equation*}
N_{+} := \#\left \{ 1\le i\le n: X_{i}>0 \right \}\quad \text{and}
\quad N_{-} := n-N_{+},
\end{equation*}
to count the number of samples with either sign. Then, by Hoeffding's inequality,
the event
\begin{equation*}
\mathcal{E}_{n}:=\left \{ \left | N_{+} - \frac{n}{2}\right | \leq
\frac{n}{4} \right \}
\end{equation*}
has probability at least $1-2e^{-n/8}$. Next, for $p\ge 1$, we consider
i.i.d. samples $Z_{1},\ldots ,Z_{p}$ from the one-sided Pareto distribution
from Example~\ref{ex:1-sided-pareto-pop}, and, as above, we use the notation
$Z^{p}_{(1)},\ldots ,Z^{p}_{(p)}$ for the order statistics of the samples
$Z_{1},\ldots , Z_{p}$; by convention, let $Z^{p}_{(0)}:=1$. We claim that
for any integers $0<a\leq b\leq p$ and any $C\geq 1$, we have
%
\begin{equation}
\label{eq:bound_pareto_os}
\mathbb{P}\left (\frac{Z^{p}_{(p-a)} }{Z^{p}_{(p-b)}} \leq \sqrt{
\frac{b}{Ca}}\right ) \leq e^{-\frac{(C-1)^{2}}{C}a}.
\end{equation}
To see this, first use the fact that, conditionally on $X_{(n-b:n)}$, the
ratios
\begin{equation*}
X_{(n-b+1:n)}/X_{(n-b:n)},\ldots ,X_{(n:n)}/X_{(n-b:n)}
\end{equation*}
are distributed exactly as $Z^{b}_{(1)},\ldots , Z^{b}_{(b)}$. Then, apply
Chernoff's inequality to the random variable
$Y\sim \text{Binom}(b,Ca/b)$, and use the well-known relative entropy bound
$D(x\parallel y) \ge (y-x)^{2}/(2y)$ for $x\leq y$. This yields
\begin{equation*}
\mathbb{P}\left (\frac{Z^{p}_{(p-a)} }{Z^{p}_{(p-b)}} \leq \sqrt{
\frac{b}{Ca}}\right ) = \mathbb{P}\left (Z^{b}_{(b-a)} \leq \sqrt{
\frac{b}{Ca}}\right ) =\mathbb{P}(Y\leq a) \leq e^{-b D( \frac{a}{b}
\parallel \frac{Ca}{b})} \leq e^{-\frac{(C-1)^{2}}{C}a}.
\end{equation*}
as claimed. Note also that in particular we have
%
\begin{equation}
\label{eq:bound_proba_event_E}
\mathbb{P}\left (Z^{p}_{(n-a)} > \sqrt{\frac{p}{Ca}}\right ) \geq 1-e^{-
\frac{(C-1)^{2}}{C}a}.
\end{equation}
by taking $b=p$. Next note that, conditionally on $\{N_{+}=n_{+}\}$, the
samples $ \{ X_{i}:X_{i}>0,1\le i \le\break n  \}$ and
$\left \{ -X_{i}:X_{i}<0,1\le i \le n \right \}$ are, respectively, independent
sequences of $n_{+}$ and $n_{-}:=n-n_{+}$ i.i.d.\ samples from
$\textnormal{Par}(2)$. Therefore, equation~\eqref{eq:bound_pareto_os} and
the union bound imply that for any $1\le a\le n$, the event
\begin{multline*}
\mathcal{F}_{n}(a;C) := \left \{ 0<X_{(n-b)}^{n} < \sqrt{\frac{Ca}{b}}
\cdot X_{(n-a)}^{n} \textnormal{ for all }a \le b \le N_{+}-1 \right
\}
\\
\cap \left \{ 0 >X_{(b)}^{n} > \sqrt{\frac{Ca}{b}} \cdot X_{(a+1)}^{n}
\textnormal{ for all } a+1 \le b \le \ N_{-} \right \}
\end{multline*}
has probability
%
\begin{equation}
\label{eq:proba_event_F}
\mathbb{P}(\mathcal{F}_{n}(a;C))\geq 1 - n \cdot e^{-
\frac{(C-1)^{2}}{C}a}.
\end{equation}
We also introduce the following event for $C\geq 1$,
\begin{equation*}
\mathcal{G}_{n}(a;C):=\left \{ X_{(a+1)}^{n} \geq -\sqrt{\frac{Cn}{a}}
\right \}\cap \left \{ X_{(n-a)}^{n} \leq \sqrt{\frac{Cn}{a}} \right
\}.
\end{equation*}
To bound the probability of failure, apply Chernoff's inequality to the
random variable $Y\sim \text{Binom}(n,a/(2Cn))$ to get for $C\geq 8$,
%
\begin{equation}
\label{eq:proba_event_G}
\mathbb{P}(\mathcal{G}_{n}^{c}(a;C)) \leq 2\mathbb{P}(Y\geq a) \leq 2e^{-nD(
\frac{a}{n}\parallel \frac{a}{2Cn})} \leq \frac{2}{(2C-1)^{a/4}}
\leq \frac{2}{C^{a/4}},
\end{equation}
where we used the bound
$D(x+y\parallel x) \geq \frac{y}{2}\log \frac{y}{x}$ whenever
$y\geq 8x$ and $x,y\leq \sfrac{1}{4}$ (e.g., see
\cite[Lemma~16]{blanchard2024tight}).

Now we return to the main proof, where the bulk of the work is to develop
estimates on the value of $\bar D_{n,2}(i)$ for moderate values of the
index $i$. We begin with the positive side. Suppose that the events
$\mathcal{E}_{n}$, $\mathcal{F}_{n}(a;C)$, and
$\mathcal{G}_{n}(a;C)$ hold for some $1\le a \le n$ and $C\geq 8$. If
$a+1\le i\le N_{+}$, then we use $\mathcal{F}_{n}(a;C)$, an elementary
inequality, and $\mathcal{G}_{n}(a;C)$ to get:
%
\begin{align}
\bar X_{n}[n-i:n-1] &= \frac{a-1}{i-1}\bar X_{n}[n-a:n-1] +
\frac{1}{i-1} S_{n}[n-i:n-a]
\notag
\\
&\leq \frac{a}{i}X_{(n)}^{n} + \frac{1}{i-1} X_{(n-a)}^{n} \sum _{j=a}^{i-1}
\sqrt{\frac{Ca}{j}}
\notag
\\
&\leq \frac{a}{i}X_{(n)}^{n} + 2\sqrt{\frac{2Ca}{i}} X_{(n-a)}^{n},
\notag
\\
&\leq \frac{a}{i}X_{(n)}^{n} + 2C\sqrt{\frac{2n}{i}}.
\label{eq:upper_bound_avg}
\end{align}
Similarly, we have
\begin{align*}
\bar X_{n}[:n-i] \geq \bar X_{n}[:N_{-}] &\geq \frac{a}{N_{-}} \bar X_{n}[:a]
+ \frac{1}{N_{-}} X_{(a+1)}^{n} \sum _{j=a+1}^{N_{-}} \sqrt{
\frac{Ca}{j}}
\notag
\\
&\geq \frac{a}{n_{-}} \bar X_{n}[:a] + 2\sqrt{\frac{2Ca}{N_{-}}} X_{(a+1)}^{n}
\notag
\\
&\geq \frac{4a}{n} X_{(1)}^{n} - 4C\sqrt 2.
\end{align*}
If we additionally have the event
\begin{equation*}
\mathcal{Z}_{n}^{-}(a;C):=\left \{ -X_{(1)}^{n} \leq \frac{Cn}{2a}
\right \},
\end{equation*}
then we obtained
%
\begin{equation}
\label{eq:lower_bound_mean-}
\bar X_{n}[:n-i] \geq \bar X_{n}[:N_{-}] \geq -8C.
\end{equation}
Next, using equation~\eqref{eq:useful_comparison} we obtain for
$\ell \geq a$,
\begin{align*}
&\mathcal{A}_{n}(\ell )
\\
&\quad :=\left \{ \bar D_{n,2}(i)>\bar D_{n,2}(n-1)
\textnormal{ for all } n-N_{+}+1\le i \le n-\ell \right \}
\\
&\quad = \left \{ X_{(n)}^{n}>\bar X_{n}[n-i:n-1] + \sqrt{i}\cdot (\bar X_{n}[n-i:n-1]
- \bar X_{n}[:n-i]) \textnormal{ for all } \ell \le i \le N_{+}-1
\right \}
\\
&\quad \supseteq \left \{ X_{(n)}^{n}>2\sqrt i\cdot \bar X_{n}[n-i:n-1] -
\sqrt i\cdot \bar X_{n}[:n-i]) \textnormal{ for all } \ell \le i \le N_{+}-1
\right \}
\end{align*}
Combining this with equation~\eqref{eq:upper_bound_avg} and
\eqref{eq:lower_bound_mean-} shows that, if additionally
$\ell \geq 16a^{2}$ (i.e., $\ell \ge 4a$), then
%
\begin{align}
\mathcal{A}_{n}(\ell ) \supseteq \mathcal{E}_{n}\cap \mathcal{F}_{n}(a;C)&
\cap \mathcal{G}_{n}(a;C) \cap \mathcal{Z}_{n}^{-}(a;C)
\label{eq21}
\\
&\cap \left \{ X_{(n)}^{n}>\frac{2a}{\sqrt i} X_{(n)}^{n} + 4C\sqrt{2n}
+8C\sqrt{i} \textnormal{ for all } \ell \le i \le N_{+}-1 \right \}
\notag
\\
\supseteq \mathcal{E}_{n}\cap \mathcal{F}_{n}(a;C)\cap \mathcal{G}_{n}&(a;C)
\cap \mathcal{Z}_{n}^{-}(a;C)\cap \left \{ X_{(n)}^{n} >25C\sqrt{n}
\right \}.
\label{eq:decomposition_positive_side}
\end{align}
We use a similar argument for the negative side. Suppose that
$\mathcal{E}_{n}$, $\mathcal{F}_{n}(a;C)$, and
$\mathcal{G}_{n}(a;C)$ hold. Then, for $a+1\le i\le N_{-}$, we compute
%
\begin{align}
\bar X_{n}[:i] &= \frac{a}{i} \bar X_{n}[:a] + \frac{1}{i} S_{n}[a:i]
\notag
\\
&\geq \frac{a}{i} \bar X_{n}[:a] + \frac{1}{i} X_{(a+1)}^{n} \sum _{j=a+1}^{i}
\sqrt{\frac{Ca}{j}}
\notag
\\
&\geq \frac{a}{i} X_{(1)}^{n} + 2\sqrt{\frac{2Ca}{i}} X_{(a+1)}^{n}
\notag
\\
&\geq \frac{a}{i} X_{(1)}^{n} - 2C\sqrt{\frac{2n}{i}}.
\label{eq:lower_bound_avg-}
\end{align}
Equation~\eqref{eq:useful_comparison} shows that, for any
$\ell \geq a$, we have
\begin{align*}
\mathcal{B}_{n}(\ell ) &:= \left \{ \bar D_{n,2}(i)>\bar D_{n,2}(n-1)
\textnormal{ for all } \ell +1 \le i \le N_{-} \right \}
\\
&\supseteq \left \{ X_{(n)}^{n}> 2\sqrt i\cdot \bar X_{n}[n-N_{+}:n-1]
- \sqrt i\cdot \bar X_{n}[:i]\textnormal{ for all } \ell +1 \le i
\le N_{-} \right \}.
\end{align*}
Last, we introduce the event
$\mathcal{Y}_{n}:=\{X_{(n)}^{n} \geq - X_{(1)}^{n}\}$. Then, for $n$ sufficiently
large, equations \eqref{eq:upper_bound_avg} and~\eqref{eq:lower_bound_avg-}
yield the following, for any $a<\sqrt n/8$ and $\ell \geq 16a^{2}$,
%
\begin{align}
\mathcal{B}_{n}(\ell ) &\supseteq \mathcal{E}_{n}\cap \mathcal{F}_{n}(a;C)
\cap \mathcal{G}_{n}(a;C) \cap \mathcal{Y}_{n}
\notag
\\
&\qquad \cap \left \{ X_{(n)}^{n}> \frac{2a\sqrt i}{N_{+}} X_{(n)}^{n}
+4C \sqrt{\frac{2ni}{N_{+}}} - \frac{a}{\sqrt i} X_{(1)}^{n} + 2C
\sqrt{2n} \textnormal{ for all }\ell \le i \le N_{+}-1 \right \}
\notag
\\
&\supseteq \mathcal{E}_{n}\cap \mathcal{F}_{n}(a;C) \cap \mathcal{G}_{n}(a;C)
\cap \mathcal{Y}_{n}\cap \left \{ \frac{1}{2}X_{(n)}^{n} > -
\frac{a}{\sqrt \ell } X_{(1)}^{n} + 6C\sqrt{2n} \right \}
\notag
\\
&\supseteq \mathcal{E}_{n}\cap \mathcal{F}_{n}(a;C)\cap \mathcal{G}_{n}(a;C)
\cap \mathcal{Y}_{n}\cap \left \{ X_{(n)}^{n} > 24C\sqrt{2n} \right
\}.
\label{eq:decomposition_negative_side}
\end{align}
We now combine equations~\eqref{eq:decomposition_positive_side} and \eqref{eq:decomposition_negative_side}
to obtain, for $n$ sufficiently large and $a<\sqrt n/8$ and
$\ell \geq 16a^{2}$, that
%
\begin{align}
\mathcal{D}_{n}(\ell ) &:= \left \{ \max (\bar{\mathcal{C}}_{n}) \ge X_{(n-
\ell +1)}^{n} \textnormal{ for all }\bar {\mathcal{C}}_{n}\in
\bar {\mathscr{C}}_{n,2} \right \} \cup \left \{ \min (
\bar{\mathcal{C}}_{n}) \le X_{(\ell )}^{n} \textnormal{ for all }
\bar {\mathcal{C}}_{n}\in \bar {\mathscr{C}}_{n,2} \right \}
\notag
\\
&\supseteq \left \{ \min _{1\le i \le n-1} \bar D_{n,2}(i) < \min _{
\ell +1\le i\le n-\ell } \bar D_{n,2}(i) \right \}
\notag
\\
&\supseteq \left \{ \bar D_{n,2}(n-1) < \min _{\ell +1\le i \le n-
\ell } \bar D_{n,2}(i) \right \} = \mathcal{A}_{n}(\ell )\cap
\mathcal{B}_{n}(\ell )
\notag
\\
&\supseteq \mathcal{E}_{n}\cap \mathcal{F}_{n}(a;C)\cap \mathcal{G}_{n}(a;C)
\cap \mathcal{Y}_{n}\cap \mathcal{Z}_{n}^{-}(a;C) \cap \left \{ X_{(n)}^{n}
> 35C\sqrt{n} \right \}
\label{eq:final_bound+}
\end{align}
By applying the symmetry $X_{(i)}^{n}\mapsto -X_{(n-i)}^{n}$, we also have
%
\begin{align}
\mathcal{D}_{n}(\ell ) &\supseteq \left \{ \bar D_{n,2}(1) > \min _{
\ell +1 \le i\le n-\ell } \bar D_{n,2}(i) \right \}
\notag
\\
&\supseteq \mathcal{E}_{n}\cap \mathcal{F}_{n}(a;C)\cap \mathcal{G}_{n}(a;C)
\cap \mathcal{Y}_{n}^{c} \cap \mathcal{Z}_{n}^{+}(a;C) \cap \left \{ -X_{(1)}^{n}
>35C\sqrt n \right \},
\label{eq:final_bound-}
\end{align}
where
\begin{equation*}
\mathcal{Z}_{n}^{+}(a;C) := \left \{ X_{(n)}^{n} \leq \frac{Cn}{2a}
\right \}.
\end{equation*}
For convenience, set
$\mathcal{Z}_{n}(a;C):= \mathcal{Z}_{n}^{+}(a;C) \cup \mathcal{Z}_{n}^{-}(a;C)$.
We will next put together equations~\eqref{eq:final_bound+} and~\eqref{eq:final_bound-},
but first we claim that
\begin{align*}
\mathcal{Z}_{n}(a;C) \cap \mathcal{Y}_{n}^{c} &= \mathcal{Z}_{n}^{+}(a;C)
\cap \mathcal{Y}_{n}^{c}
\\
\mathcal{Z}_{n}(a;C) \cap \mathcal{Y}_{n} &= \mathcal{Z}_{n}^{-}(a;C)
\cap \mathcal{Y}_{n}
\end{align*}
Indeed, for the first equality, if $\mathcal{Z}_{n}^{+}(a;C)$ does not
hold, because of $\mathcal{Y}_{n}^{c}$ we have
$-X_{(1)}^{n} > X_{(n)}^{n} > \frac{Cn}{2a}$, that is
$\mathcal{Z}_{n}^{-}(a;C)$ also does not hold. This contradicts
$\mathcal{Z}_{n}(a;C)$. Similarly, for the second equality, if
$\mathcal{Z}_{n}^{-}(a;C)$ did not hold, we would have
$X_{(n)}^{n} \geq -X_{(1)}^{n} > \frac{Cn}{2a}$, contradicting
$\mathcal{Z}_{n}(a;C)$. Therefore,
%
\begin{align}
&\mathcal{Z}_{n}(a;C) \cap \left \{ X_{(n)}^{n}>35C\sqrt n \right \}
\notag
\\
&\quad = \left ( \mathcal{Z}_{n}^{-}(a;C) \cap \mathcal{Y}_{n} \cap \left
\{ X_{(n)}^{n}>35C\sqrt n \right \} \right ) \cup \left ( \mathcal{Z}_{n}^{+}(a;C)
\cap \mathcal{Y}_{n}^{c} \cap \left \{ X_{(n)}^{n}>35C\sqrt n \right
\} \right )
\notag
\\
&\quad \subseteq \left ( \mathcal{Z}_{n}^{-}(a;C) \cap \mathcal{Y}_{n}
\cap \left \{ X_{(n)}^{n}>35C\sqrt n \right \} \right ) \cup \left (
\mathcal{Z}_{n}^{+}(a;C) \cap \mathcal{Y}_{n}^{c} \cap \left \{ -X_{(1)}^{n}>35C
\sqrt n \right \} \right ).
\label{eq:useful_event_property}
\end{align}
We also introduce the last event
\begin{equation*}
\mathcal{H}_{n}(\ell ;C) = \left \{ X_{(n-\ell +1)}^{n} \geq \sqrt{
\frac{n}{4C\ell }} \right \} \cap \left \{ X_{(\ell )}^{n} \leq -
\sqrt{\frac{n}{4C\ell }} \right \}.
\end{equation*}
Using equation~\eqref{eq:bound_proba_event_E}, we obtain
%
\begin{equation}
\label{eq:proba_event_H}
\mathbb{P}(\mathcal{H}_{n}(\ell ;C)\mid \mathcal{E}_{n}) \geq 1-2e^{-
\frac{(C-1)^{2}}{C}\ell}
\end{equation}
provided that $\ell <n/4$.

Now we put everything together. Combining equations~\eqref{eq:final_bound+},~\eqref{eq:final_bound-},
and~\eqref{eq:proba_event_H}, we proved
%
\begin{align}
\mathcal{C}_{n}(\ell ) &:=\left \{ \max (\bar{\mathcal{C}}_{n}) \ge
\sqrt{\frac{n}{4C\ell }} \textnormal{ for all }\bar {\mathcal{C}}_{n}
\in \bar {\mathscr{C}}_{n,2} \right \} \cup \left \{ \min (
\bar{\mathcal{C}}_{n}) \le -\sqrt{\frac{n}{4C\ell }}
\textnormal{ for all }\bar {\mathcal{C}}_{n}\in \bar {\mathscr{C}}_{n,2}
\right \}
\notag
\\
&\supseteq \mathcal{D}_{n}(\ell ) \cap \mathcal{H}_{n}(\ell ;C)
\notag
\\
&\supseteq \mathcal{E}_{n}\cap \mathcal{F}_{n}(a;C)\cap \mathcal{G}_{n}(a;C)
\cap \mathcal{H}_{n}(\ell ;C)\cap \mathcal{Z}_{n}(a;C) \cap \left \{ X_{(n)}^{n}>35C
\sqrt n \right \},
\label{eq:final_decomposition}
\end{align}
To bound the probability of the right-hand side, we first compute
%
\begin{equation}
\label{eq:proba_event_Z}
\mathbb{P}(\mathcal{Z}_{n}^{c}(a;C)) \leq \mathbb{P}(Y\geq 2) \leq
\frac{16a^{4}}{C^{4}n^{2}}.
\end{equation}
where $Y\sim \text{Binom}(n,4a^{2}/(C^{2}n^{2}))$. Then, using all the failure
probability bounds from \eqref{eq:proba_event_F},
\eqref{eq:proba_event_G}, \eqref{eq:proba_event_H}, and~\eqref{eq:proba_event_Z},
there exists a constant $c_{0}>0$ such that choosing
$a_{n}=\left \lceil c_{0}\log n \right \rceil $ and
$\ell _{n}=16 a_{n}^{2}$, and $C=8$, yields
\begin{equation*}
\sum _{n= 1}^{\infty} \bigg(\mathbb{P}(\mathcal{E}_{n}^{c}) +
\mathbb{P}(\mathcal{F}_{n}^{c}(a_{n};C)\mid \mathcal{E}_{n}) +
\mathbb{P}(\mathcal{G}_{n}^{c}(a_{n};C)) + \mathbb{P}(\mathcal{H}_{n}^{c}(l_{n};C)
\mid \mathcal{E}_{n}) + \mathbb{P}(\mathcal{Z}_{n}^{c}(a_{n};C))
\bigg) <\infty .
\end{equation*}
Hence, by the Borel-Cantelli lemma, the event $\mathcal{E}^{(0)}$ that
for large enough $n$, the events $\mathcal{E}_{n}$,
$\mathcal{F}_{n}(a_{n};C)$, $\mathcal{G}_{n}(a_{n};C)$,
$\mathcal{H}_{n}(\ell _{n};C)$, $\mathcal{Z}_{n}(a_{n};C)$ are all satisfied,
has probability one. Next, let $n_{k}:=2^{k}$ and note that for any
$k\geq 1$,
\begin{equation*}
\mathbb{P}\left (\max _{n_{k-1}<i\le n_{k}} X_{i}> 35C\sqrt{n_{k}}
\right ) =1-\left ( 1-\frac{1}{(35C)^{2} n_{k}} \right )^{n_{k}-n_{k-1}}
\geq 1-e^{-1/(2(35C)^{2})} \geq \frac{1}{2^{11}C^{2}}.
\end{equation*}
Further, all these events are independent. Hence, using the Borel-Cantelli
lemma, the event $\mathcal{E}^{(1)}$ where these occur infinitely often
has probability one as well. Combining this with equation~\eqref{eq:final_decomposition}
implies
\begin{align*}
\mathcal{E}^{(0)}\cap \mathcal{E}^{(1)} &\subseteq \mathcal{E}^{(0)}
\cap \left \{ X_{(n)}^{n} >35C\sqrt n \text{ for infinitely many }n
\in \mathbb{N} \right \}
\\
&\subseteq \left \{ \mathcal{C}_{n}(\ell _{n})
\text{ holds for infinitely many }n\in \mathbb{N} \right \}.
\end{align*}
Because $\mathcal{E}^{(0)}$ and $\mathcal{E}^{(1)}$ have full probability
and $\ell _{n}=O((\log n)^{2})$, this ends the proof.
\end{proof}

In order to prove our positive results on convergence of the empirical
cluster centers, when imposing suitable balance constraints, we needed
an auxiliary result on continuity of the excess distortion; its proof is
as follows.

\begin{proof}[Proof of Lemma~\ref{lem:D-cts}]
For convenience, we write $D_{\mu}(\mathcal{C}\,|\,\mathcal{C}')$ for the
excess distortion when $X$ has distribution $\mu $. Now, suppose that
$\{(\mu _{n},\mathcal{C}_{n},\mathcal{C}_{n}')\}_{n\in \mathbb{N}}$ and
$(\mu ,\mathcal{C},\mathcal{C}')$ in
$\mathcal{W}_{1}(\mathbb{R}^{d})\times \mathscr{C}_{k}\times
\mathscr{C}_{k}$ satisfy $W_{1}(\mu _{n},\mu )\to 0$,
$d_{\textnormal{H}}(\mathcal{C}_{n},\mathcal{C})\to 0$, and
$d_{\textnormal{H}}(\mathcal{C}_{n}',\mathcal{C}')\to 0$ as
$n\to \infty $. Then construct a probability space
$(\tilde{\Omega},\tilde{\mathscr{F}},\tilde{\mathbb{P}})$ on which are
defined random variables $\{Z_{n}\}_{n\in \mathbb{N}}$ and $Z$ such that
$Z_{n}$ has distribution $\mu _{n}$ for all $n\in \mathbb{N}$ and
$Z$ has distribution $\mu $, and such that we have $Z_{n}\to Z$ holding
$\tilde{\mathbb{P}}$-almost surely as well as
$\tilde{\mathbb{E}}\|Z_{n}-Z\|\to 0$. (This is essentially an application
of Skorokhod's representation theorem; see
\cite[Lemma~3]{JaffeClustering} for details.) We now claim that
%
\begin{equation}
\label{eqn:min-diffs}
\left \{\min _{c_{n}\in \mathcal{C}_{n}}\|c_{n}-Z_{n}\|^{2}-\min _{c_{n}'
\in \mathcal{C}_{n}'}\|c_{n}'-Z_{n}\|^{2}\right \}_{n\in \mathbb{N}}
\end{equation}
is uniformly integrable. To show this, note that the assumptions
$d_{\textnormal{H}}(\mathcal{C}_{n},\mathcal{C})\to 0$ and
$d_{\textnormal{H}}(\mathcal{C}_{n}',\mathcal{C}')\to 0$ imply that
\begin{equation*}
\left \{\sum _{c_{n}\in \mathcal{C}_{n}}\|c_{n}\|^{2}\right \}_{n\in
\mathbb{N}},\qquad \left \{\sum _{c_{n}'\in \mathcal{C}_{n}'}\|c_{n}'
\|^{2}\right \}_{n\in \mathbb{N}},\qquad \textnormal{ and } \qquad
\left \{d_{\textnormal{H}}(\mathcal{C}_{n},\mathcal{C}_{n}')\right \}_{n
\in \mathbb{N}}
\end{equation*}
are uniformly bounded, and observe that we have
\begin{align*}
\left |\min _{c_{n}\in \mathcal{C}_{n}}\|c_{n}-Z_{n}\|^{2}-\min _{c_{n}'
\in \mathcal{C}_{n}'}\|c_{n}'-Z_{n}\|^{2}\right | \le \sum _{c_{n}
\in \mathcal{C}_{n}}\|c_{n}\|^{2} + \sum _{c_{n}'\in \mathcal{C}_{n}'}
\|c_{n}'\|^{2} + 2d_{\textnormal{H}}(\mathcal{C}_{n},\mathcal{C}_{n}')
\|Z_{n}\|.
\end{align*}
Since $\{Z_{n}\}_{n\in \mathbb{N}}$ is uniformly integrable by construction,
we conclude that \eqref{eqn:min-diffs} is indeed uniformly integrable.
Since we also have
\begin{equation*}
\min _{c_{n}\in \mathcal{C}_{n}}\|c_{n}-Z_{n}\|^{2}-\min _{c_{n}'\in
\mathcal{C}_{n}'}\|c_{n}'-Z_{n}\|^{2}\to \min _{c\in \mathcal{C}}\|c-Z
\|^{2}-\min _{c'\in \mathcal{C}'}\|c'-Z\|^{2}
\end{equation*}
holding $\tilde{\mathbb{P}}$-almost surely, it follows that we can compute:
\begin{align*}
D_{\mu _{n}}(\mathcal{C}_{n}\,|\,\mathcal{C}_{n}') &= \int _{
\mathbb{R}^{d}}\left (\min _{c_{n}\in \mathcal{C}_{n}}\|c_{n}-x\|^{2}-
\min _{c_{n}'\in \mathcal{C}_{n}'}\|c_{n}'-x\|^{2}\right )
\textnormal{d}\mu _{n}(x)
\\
&= \tilde{\mathbb{E}}\left [\min _{c_{n}\in \mathcal{C}_{n}}\|c_{n}-Z_{n}
\|^{2}-\min _{c_{n}'\in \mathcal{C}_{n}'}\|c_{n}'-Z_{n}\|^{2}\right ]
\\
&= \tilde{\mathbb{E}}\left [\min _{c\in \mathcal{C}}\|c-Z\|^{2}-\min _{c'
\in \mathcal{C}'}\|c'-Z\|^{2}\right ]
\\
&= \int _{\mathbb{R}^{d}}\left (\min _{c\in \mathcal{C}}\|c-x\|^{2}-
\min _{c'\in \mathcal{C}'}\|c'-x\|^{2}\right )\textnormal{d}\mu (x) = D_{
\mu}(\mathcal{C}\,|\,\mathcal{C}').
\end{align*}
This finishes the proof.
\end{proof}

\subsection{From Section~\ref{sec:distortion}}
\label{appC.2}

Lastly, we prove the results from Section~\ref{sec:distortion} on convergence
and non-convergence of the excess distortion of the empirical cluster centers,
when it is well-defined in case (ii).

We begin with the proof of Example~\ref{ex:inconsistency_two_sided_imbalanced}
on the asymmetric two-sided $\textnormal{Par}(2)$ distribution, which is
similar to the proof of Example~\ref{ex:inconsistency_two_sided} on the
symmetric two-sided $\textnormal{Par}(2)$ distribution.

\begin{proof}[Proof of Example~\ref{ex:inconsistency_two_sided_imbalanced}]
The proof uses the same tools as in the proof of Example~\ref{ex:inconsistency_two_sided}.
We also use the same notation for the order statistics
$X_{(1)}^{n},\ldots ,X_{(n)}^{n}$ and the quantities $S_{n}$,
$\bar X_{n}$, $\bar D_{n,2}$, and
$N_{+}=\#\left \{ 1\le i \le n:X_{i}>0 \right \}=n-N_{-}$. Since the distribution
of $N_{+}$ is now Binomial with parameters $n$ and $p$, we introduce the
event
\begin{equation*}
\mathcal{E}_{n}:=\left \{ \left |N_{+}-pn\right |\leq \frac{pn}{2}
\right \},
\end{equation*}
which by Hoeffding's inequality has probability at least
$1-2e^{-(1-p)^{2} n/4}$. Next, for $1\le a\le n$ and $C\geq 1$, we introduce
the same event as in Example~\ref{ex:inconsistency_two_sided}:
\begin{multline*}
\mathcal{F}_{n}(a;C) := \left \{ 0<X_{(n-b)}^{n} < \sqrt{\frac{Ca}{b}}
\cdot X_{(n-a)}^{n} \textnormal{ for all } a\le b\le N_{+} - 1
\right \}
\\
\cap \left \{ 0 >X_{(b)}^{n} > \sqrt{\frac{Ca}{b}} \cdot X_{(a+1)}^{n}
\textnormal{ for all } a+1\le b \le N_{-} \right \},
\end{multline*}
which has probability at least $1-n e^{-\frac{(C-1)^{2}}{C}a}$ using the
same arguments (see equation~\eqref{eq:proba_event_F}). Next, let
\begin{equation*}
\mathcal{G}_{n}(a;C):=\left \{ X_{(a+1)}^{n} \geq -\sqrt{
\frac{2(1-p)Cn}{a}} \right \}\cap \left \{ X_{(n-a)}^{n} \leq \sqrt{
\frac{2pCn}{a}} \right \}.
\end{equation*}
Then, equation~\eqref{eq:proba_event_G} still holds, that is,
$\mathbb{P}(\mathcal{G}_{n}(a;C))\geq 1-2C^{-a/4}$.

We suppose from now that $\mathcal{E}_{n}$, $\mathcal{F}_{n}(a;C)$,
$\mathcal{G}_{n}(a;C)$, and
\begin{equation*}
\mathcal{Z}_{n}^{-}(a;C):=\left \{ -X_{(1)}^{n} \leq
\frac{C(1-p)n}{2a} \right \}
\end{equation*}
hold. Then fix $a+1\le i\le N_{+}$. Slightly adjusting the arguments for
equation~\eqref{eq:upper_bound_avg}, we obtain
%
\begin{equation}
\label{eq:upper_n_v2}
\bar X_{n}[n-i:n-1] \leq \frac{a}{i}X_{(n)}^{n} + 4C\sqrt{
\frac{pn}{i}}.
\end{equation}
Adjusting the arguments for equation~\eqref{eq:lower_bound_mean-} and using
$\mathcal{E}_{n}$, we also obtain
\begin{equation*}
\bar X_{n}[:n-i] \geq \bar X_{n}[:N_{-}] \geq \frac{a}{N_{-}}X_{(1)}^{n}
- 4C\sqrt{\frac{(1-p)n}{N_{-}}}\overset{(i)}{\geq} \frac{4a}{n}X_{(1)}^{n}
- 8C \geq -10C,
\end{equation*}
Therefore, following the arguments from equation~\eqref{eq:decomposition_positive_side}
we have for $\ell \geq 16a^{2}$,
\begin{align*}
\mathcal{A}_{n}(\ell )&:=\left \{ \bar D_{n,2}(i)>\bar D_{n,2}(n-1)
\textnormal{ for all } n-N_{+} + 1\le i \le n-\ell \right \}
\\
&\supseteq \mathcal{E}_{n}\cap \mathcal{F}_{n}(a;C)\cap \mathcal{G}_{n}(a;C)
\cap \mathcal{Z}_{n}^{-}(a;C)
\\
&\qquad \cap \left \{ X_{(n)}^{n} >\frac{2a}{\sqrt i} X_{(n)}^{n} + 8C
\sqrt{pn} + 10C\sqrt i \textnormal{ for all } \ell \le i \le N_{+}-1
\right \}
\\
&\subseteq \mathcal{E}_{n}\cap \mathcal{F}_{n}(a;C)\cap \mathcal{G}_{n}(a;C)
\cap \mathcal{Z}_{n}^{-}(a;C)\cap \left \{ X_{(n)}^{n} >50C\sqrt{pn}
\right \}.
\end{align*}
Next, we turn to indices $i\leq N_{-}$. For $\ell \le i\leq N_{-}$, equation~\eqref{eq:lower_bound_avg-}
becomes
\begin{equation*}
-\sqrt i\cdot \bar X_{n}[:i] \leq -\frac{a}{\sqrt i} X_{(1)}^{n} + 4C
\sqrt{(1-p)n} \leq -\frac{1}{2} X_{(1)}^{n} + 4C\sqrt{n} ,
\end{equation*}
where in the last inequality we used $\ell \geq 16a^{2}$. Next, we introduce
the event
%
\begin{multline}
\label{eq:defnition_tilde_G_n}
\tilde{\mathcal{G}}_{n}(\ell ;C):= \left \{ X_{(b)}^{n} \geq -\sqrt{
\frac{2(1-p)Cn}{b}} \textnormal{ for all }1\le b\le \ell \right \}
\\
\supseteq \left \{ -\sqrt i\cdot \bar X_{n}[:i] \leq 2\sqrt{2Cn}
\textnormal{ for all } 1 \le i \le \ell \right \}
\end{multline}
for $C\geq 8$, which by Chernoff's inequality and the union bound satisfies
%
\begin{equation}
\label{eq:bound_tilde_G_event}
\mathbb{P}(\tilde{\mathcal{G}}_{n}^{c}(\ell ;C)) \leq \sum _{1\le b
\le \ell}e^{-nD(\frac{b}{n}\parallel \frac{b}{2Cn})} \leq
\frac{2}{C^{1/4}-1},
\end{equation}
where in the last inequality we used the same bound as in equation~\eqref{eq:proba_event_G}.
In particular, under $\tilde{\mathcal{G}}_{n}(\ell ;C)$ the two previous
equations show that for $\ell \geq 16a^{2}$, we have
%
\begin{equation}
\label{eq:bound_last_term_i}
-\sqrt i\cdot \bar X_{n}[:i] \leq 5C\sqrt n\textnormal{ for all } 1
\le i\le N_{-}.
\end{equation}
For $n$ sufficiently large, $a<p\sqrt n/8$, and $\ell \geq 16a^{2}$, we
obtain the following, by using \eqref{eq:useful_comparison},~\eqref{eq:upper_n_v2},~\eqref{eq:bound_last_term_i},
and $\mathcal{E}_{n}$:
\begin{align*}
\mathcal{B}_{n}(\ell ) &:= \left \{ \bar D_{n,2}(i)>\bar D_{n,2}(n-1)
\textnormal{ for all }1 \le i \le N_{-} \right \}
\\
&\supseteq \left \{ X_{(n)}^{n}> 2\sqrt i\cdot \bar X_{n}[n-N_{+}:n-1]
- \sqrt i\cdot \bar X_{n}[:i] \textnormal{ for all }1 \le i \le N_{-}
\right \}
\\
&\supseteq \mathcal{E}_{n}\cap \mathcal{F}_{n}(a;C)\cap \mathcal{G}_{n}(a;C)
\cap \tilde {\mathcal{G}}_{n}(\ell ;C)
\\
&\qquad \cap \left \{ X_{(n)}^{n} > \frac{2a\sqrt n}{N_{+}} X_{(n)}^{n}
+8C \sqrt{\frac{pn^{2}}{N_{+}}} +5C\sqrt{n} \textnormal{ for all } 1
\le i \le N_{-} \right \}
\\
&\supseteq \mathcal{E}_{n}\cap \mathcal{F}_{n}(a;C)\cap \mathcal{G}_{n}(a;C)
\cap \tilde {\mathcal{G}}_{n}(\ell ;C)\cap \left \{ X_{(n)}^{n} > 40C
\sqrt{n} \right \}.
\end{align*}
Finally, we define the event
\begin{equation*}
\mathcal{H}_{n}(\ell ;C):=\left \{ X_{(n-\ell +1)}^{n} \geq \sqrt{
\frac{pn}{2C\ell }} \right \},
\end{equation*}
which, using the same arguments as in equation~\eqref{eq:proba_event_H},
satisfies
$\mathbb{P}(\mathcal{H}_{n}(\ell ;C)\mid \mathcal{E}_{n}) \geq 1-e^{-
\frac{(C-1)^{2}}{C}\ell}$ provided that $\ell <pn/2$.

In summary, we have
\begin{align*}
\mathcal{C}_{n}(\ell )&:=\left \{ \max (\bar {\mathcal{C}}_{n})\geq
\sqrt{\frac{pn}{2C\ell }} \textnormal{ for all } \bar {\mathcal{C}}_{n}
\in \bar {\mathscr{C}}_{n,2} \right \}
\\
&\supseteq \left \{ \max (\bar {\mathcal{C}}_{n})\geq X_{(n-\ell +1)}^{n}
\textnormal{ for all } \bar {\mathcal{C}}_{n}\in \bar {\mathscr{C}}_{n,2}
\right \} \cap \mathcal{H}_{n}(\ell ;C)
\\
&\supseteq \left \{ \bar D_{n,2}(n-1) < \min _{1\le i \le n-\ell }
\bar D_{n,2}(i) \right \}\cap \mathcal{H}_{n}(\ell ;C) = \mathcal{A}_{n}(
\ell )\cap \mathcal{B}_{n}(\ell )\cap \mathcal{H}_{n}(\ell ;C)
\\
&\supseteq \mathcal{E}_{n}\cap \mathcal{F}_{n}(a;C)\cap \mathcal{G}_{n}(a;C)
\cap \tilde {\mathcal{G}}_{n}(\ell ;C)\cap \mathcal{H}_{n}(\ell ;C)
\cap \mathcal{Z}_{n}^{-}(a;C)\cap \left \{ X_{(n)}^{n} > 50C\sqrt{n}
\right \}.
\end{align*}
Using the failure probability bounds for each event, there exists
$c_{0}>0$ such that choosing
$a_{n}=\left \lceil c_{0}\log n \right \rceil $,
$\ell _{n}=16a_{n}^{2}$, and $C=8$, we have
\begin{equation*}
\sum _{n\in \mathbb{N}} \bigg(\mathbb{P}(\mathcal{E}_{n}^{c}) +
\mathbb{P}(\mathcal{F}_{n}^{c}(a_{n};C)\,|\,\mathcal{E}_{n}) +
\mathbb{P}(\mathcal{G}_{n}^{c}(a_{n};C)) + \mathbb{P}(\mathcal{H}_{n}^{c}(
\ell _{n};C)\mid \mathcal{E}_{n})\bigg) <\infty .
\end{equation*}
Hence, by the Borel-Cantelli lemma, there is an event
$\mathcal{E}^{(0)}$ of probability one such that for large enough
$n$, all of the events $\mathcal{E}_{n}$, $\mathcal{F}_{n}(a_{n};C)$,
$\mathcal{G}_{n}(a_{n};C)$, $\mathcal{H}_{n}(\ell _{n};C)$ are satisfied.
Next, for $n$ sufficiently large, we have
$\tilde {\mathcal{G}}_{n}(\ell ;C)\subseteq \mathcal{Z}_{n}^{-}(a;C)$.
Therefore, we will focus on the events
$\mathcal{I}(n):=\tilde {\mathcal{G}}_{n}(\ell _{n};C) \cap \left \{ X_{(n)}^{n}
> 50C\sqrt{n} \right \}$ from now onward. We construct by induction a random
sequence $\{N_{k}\}_{k\in \mathbb{N}}$, where $N_{1}=1$, and, having defined
$N_{k-1}$, we let $N_{k}$ be the smallest integer such that
$N_{k}\geq 4N_{k-1}$ and
\begin{equation*}
X_{(1)}^{2N_{k-1}} > -\sqrt{\frac{2(1-p)CN_{k}}{\ell _{N_{k}}}}.
\end{equation*}
Such an integer $N_{k}$ exists since the right-hand side diverges to
$-\infty $ as $N_{k}\to \infty $.

By construction, we use equation~\eqref{eq:bound_tilde_G_event} to bound,
for any $k\in \mathbb{N}$:
\begin{equation*}
\mathbb{P}\left (\tilde {\mathcal{G}}_{N_{k}}(\ell _{n_{k}};C) \mid N_{k},X_{1},
\ldots ,X_{2N_{k-1}}\right ) \geq \mathbb{P}\left (
\tilde {\mathcal{G}}_{N_{k}}(\ell _{N_{k}};C)\mid N_{k}\right ) \geq 1-
\frac{2}{C^{1/4}-1},
\end{equation*}
Next, for $k\in \mathbb{N}$ sufficiently large, use
$N_{k}/2\geq 2N_{k-1}$ to get
\begin{align*}
&\mathbb{P}\left ( X_{(N_{k})}^{N_{k}}>50C\sqrt {N_{k}} \mid N_{k},X_{1},
\ldots ,X_{2N_{k-1}},\tilde {\mathcal{G}}_{N_{k}}(\ell _{N_{k}};C)
\right )
\\
&\geq \mathbb{P}\left ( \max _{N_{k}/2<i\leq N_{k}} X_{i}>50C\sqrt {N_{k}}
\mid N_{k},X_{1},\ldots ,X_{2N_{k-1}},\tilde {\mathcal{G}}_{N_{k}}(
\ell _{N_{k}};C)\right )
\\
&= \mathbb{P}\left ( \max _{N_{k}/2<i\leq N_{k}} X_{i}>50C\sqrt {N_{k}}
\mid N_{k}\right ) =\left ( 1-\frac{p}{(50C)^{2} N_{k}} \right )^{N_{k}/2}
\geq \frac{p}{2^{16}C^{2}}.
\end{align*}
As a result, for $k$ sufficiently large, we obtained
\begin{equation*}
\mathbb{P}(\mathcal{I}(N_{k})\mid N_{k},X_{1},\ldots ,X_{2N_{k-1}})
\geq c_{0}p
\end{equation*}
for some universal constant $c_{0}>0$. Note that all the events
$\mathcal{I}(N_{1}),\ldots ,\mathcal{I}(N_{k-1})$ are
$\sigma (N_{k},X_{1},\ldots ,X_{2N_{k-1}})$-measurable. Hence, a variant
of the Borel--Cantelli lemma (or the Azuma--Hoeffding inequality) implies
that the event $\mathcal{E}^{(1)}$ on which the events
$\mathcal{I}(N_{k})$ occur for infinitely many $k\in \mathbb{N}$ has probability
one. On $\mathcal{E}^{(0)}\cap \mathcal{E}^{(1)}$ which has probability
one, the event $\mathcal{C}_{n}(\ell _{n})$ occurs infinitely often. Consequently,
writing $D(r)$ for the excess distortion for the partition
$\{(-\infty ,r),[r,\infty )\}$, as defined in Example~\ref{ex:population_level_imbalanced_2_sided},
and using the computations from the proof of Example~\ref{ex:population_level_imbalanced_2_sided},
we get that the following holds on
$\mathcal{E}^{(0)}\cap \mathcal{E}^{(1)}$:
\begin{equation*}
\limsup _{n\to \infty} \min _{\bar{\mathcal{C}}_{n}\in
\bar{\mathscr{C}}_{n,2}} D_{X}(\bar{\mathcal{C}}_{n}\,|\,\mathcal{C}_{0})
\geq \lim _{r\to \infty} D(r)=-4(p+(2p-1)^{2}) = \inf _{\mathcal{C}
\in \mathscr{C}_{2}}D_{X}(\mathcal{C}\,|\,\mathcal{C}_{0}) + 4(1-2p),
\end{equation*}
This ends the proof.
\end{proof}

Finally, we turn to the positive results that can be obtained in case (ii).
All of our results are based on the uniform convergence guaranteed by Lemma~\ref{lemma:uniform_convergence_voronoi},
whose proof we now give.

\begin{proof}[Proof of Lemma~\ref{lemma:uniform_convergence_voronoi}]
By assumption, we have
$\inf _{\mathcal{C}\in \mathscr{C}_{k}} D(\mathcal{C}\, | \,
\mathcal{C}_{0})>-\infty $. From Lemma~\ref{lemma:charact_finite_obj} this
implies that there exists a constant $r>0$ (depending only on the distribution
of $X$, and, specifically, not on $k$) such that we have
$\mathbb{P}(\|X\|\geq t)\leq r t^{-2}$ for all $t\geq 0$.

Let $\mathcal{H}_{k}$ be the set of all intersections of at most $k$ half-spaces
in $\mathbb{R}^{d}$. From \cite{blumer1989learnability}, the VC dimension
of $\mathcal{H}_{k}$ is bounded above by $c_{0} dk\log k$ for some universal
constant $c_{0}>0$. Also, let $\bar \mu _{n}$ denote the empirical measure
of the samples $X_{1},\ldots , X_{n}$. Then, the classical uniform convergence
theorem for VC classes (e.g. \cite[Theorem 5.1]{boucheron2005theory}) implies
that, with probability $1-\delta $,
%
\begin{align}
&\left |\mathbb{P}(X\in H) - \bar \mu _{n}(H)\right |
\nonumber \\
&\quad  \leq c_{1}\sqrt{
\bar\mu _{n}(H) \frac{dk\log k \log n + \log \frac{1}{\delta}}{n} } + c_{1}
\frac{dk\log k \log n + \log \frac{1}{\delta}}{n},\quad H\in
\mathcal{H}_{k},
\label{eq:main_H_event}
\end{align}
for some universal constant $c_{1}>0$. We denote this event by
$\mathcal{E}_{n}(\delta )$. For convenience, we denote
$C_{\delta}(n):=dk\log k\log n+\log \frac{1}{\delta}$. In particular, there
exists a sufficiently large universal constant $c_{2}\geq 1$ such that
the following holds. Under $\mathcal{E}_{n}(\delta )$, any
$H\in \mathcal{H}_{k}$ with
$\bar\mu _{n}(H)\geq c_{2} \frac{C_{\delta}(n)}{n}$ satisfies
%
\begin{equation}
\label{eq:uniform_concentration_H_k}
\frac{\mathbb{P}(X\in H)}{1+ 2c_{1}\sqrt{\frac{C_{\delta}(n)}{n\bar\mu _{n}(H)}}}
\leq \bar\mu _{n}(H) \leq \left ( 1+ 2c_{1}\sqrt{
\frac{C_{\delta }(n)}{n\bar \mu _{n}(H)}} \right ) \mathbb{P}(X\in H).
\end{equation}
We use equation~\eqref{eq:uniform_concentration_H_k} to bound the deviation
of the excess distortion for any Voronoi partition $\mathscr{V}$ with
$1\le \#\mathscr{V}\le k$. Indeed, denoting by
$\mathcal{C}\in \mathscr{C}_{k}$ a set of clusters centers associated with
the Voronoi partition $\mathscr{V}$, note that any Voronoi region
$\mathcal{V}=\mathcal{V}_{\mathcal{C}}(c)$ (for $c\in \mathcal{C}$) is
an intersection of at most $k-1$ half-spaces:
\begin{equation*}
\mathcal{V}= \bigcap _{c'\in \mathcal{C}\setminus \{c\}} \left \{ x
\in \mathbb{R}^{d}: 2(c'-c)^{\top }x \leq \|c'\|^{2} - \|c\|^{2}
\right \}\in \mathcal{H}_{k-1}.
\end{equation*}
We now compare $\bar D_{n}(\mathscr{V})$ and $D_{X}(\mathscr{V})$ under
$\mathcal{E}_{n}(\delta )$. Suppose that each Voronoi region from
$\mathscr{V}$ contains at least
$\max (c_{2},4c_{1}^{2})C_{\delta}(n)$ samples, that is,
$\mathscr{V}\in \bar{\mathscr P}_{n,k,\gamma}$ for
$\gamma = \max (c_{2},4c_{1}^{2})C_{\delta}(n)$. Then, under
$\mathcal{E}_{n}(\delta )$, since $\mathcal{V}\in \mathcal{H}_{k}$, equation~\eqref{eq:uniform_concentration_H_k}
gives
%
\begin{equation}
\label{eq:lower_bound_nb_points_per_cluster}
\frac{\mathbb{P}(X\in \mathcal{V})}{1+ 2c_{1}\sqrt{\frac{C_{\delta}(n)}{n\bar\mu _{n}(\mathcal{V})}}}
\leq \bar\mu _{n}(\mathcal{V}) \leq \left ( 1+ 2c_{1}\sqrt{
\frac{C_{\delta }(n)}{n\bar \mu _{n}(\mathcal{V})}} \right )
\mathbb{P}(X\in \mathcal{V}).
\end{equation}
In particular, since $\mathcal{V}$ contains at least
$4c_{1}^{2} C_{\delta}(n)$ samples, we also have
$\mathbb{P}(X\in \mathcal{V})/2 \leq \bar\mu _{n}(\mathcal{V}) \leq 2
\mathbb{P}(X\in \mathcal{V})$.

Most of the work to compare $\bar D_{n}(\mathscr{V})$ and
$D_{X}(\mathscr{V})$ will be in comparing $m_{X}(\mathcal{V})$ and
$\bar m_{n}(\mathcal{V})$ for any fixed $\mathcal{V}\in \mathscr{V}$. This
requires some preliminary estimates. To do so, consider a unit-norm vector
$e\in \mathbb{R}^{d}$ to be fixed later, and divide the samples in
$\mathcal{V}$ into the two following sets:
\begin{align*}
\mathcal{S}^{+}(\mathcal{V};e)&:= \{e^{\top }X_{i}: 1\le i \le n,\, X_{i}
\in \mathcal{V},\, e^{\top }X_{i}> 0\},
\\
\mathcal{S}^{-}(\mathcal{V};e)&:= \{-e^{\top }X_{i}: 1\le i \le n,\, X_{i}
\in \mathcal{V},\, e^{\top }X_{i}<0\}.
\end{align*}
For convenience, we order these sets writing
$\mathcal{S}^{+}(\mathcal{V};e)=:\{z_{1}^{+}<\ldots <z_{\ell ^{+}}^{+}
\}$ and
$\mathcal{S}^{-}(\mathcal{V};e)=:\{z_{1}^{-}<\ldots <z_{\ell ^{-}}^{-}
\}$, where $\ell ^{+}:=|\mathcal{S}^{+}(\mathcal{V};e)|$ and
$\ell ^{-}:=|\mathcal{S}^{-}(\mathcal{V};e)|$. In words, we order the samples
falling in the Voronoi partition of $c$ by their projection onto the unit-vector
$e$. The upshot of these definitions is that
\begin{equation*}
e^{\top }\bar m_{n}(\mathcal{V}) =
\frac{1}{n\bar\mu _{n}(\mathcal{V})} \sum _{i=1}^{n} e^{\top }X_{i}
\cdot \mathbbold{1}\{X_{i}\in \mathcal{V}\},
\end{equation*}
and as a result, letting $z_{0}^{+}=z_{0}^{-}=0$, we have
%
\begin{multline}
\label{eq:bar_m_expression_useful}
\bar\mu _{n}(\mathcal{V})\cdot e^{\top }\bar m_{n}(\mathcal{V}) =
\underbrace{\sum _{j=1}^{\ell ^{+}} (z_{j}^{+}-z_{j-1}^{+}) \cdot \bar\mu _{n}\left ( \mathcal{V}\cap \{x\in \mathbb{R}^{d}: e^{\top }x \geq z_{j}^{+}\} \right )}_{:=
\bar A^{+}_{n}(\mathcal{V};e)}
\\
-
\underbrace{\sum _{j=1}^{l^{-}} (z_{j}^{-}-z_{j-1}^{-}) \cdot \bar\mu _{n}\left ( \mathcal{V}\cap \{x\in \mathbb{R}^{d}: -e^{\top }x \geq z_{j}^{-}\} \right )}_{:=
\bar A^{-}_{n}(\mathcal{V};e)}.
\end{multline}
By construction, each term of both sums is non-negative. Further, all sets
on the right-hand side belong to $\mathcal{H}_{k}$. We treat each sum separately.
Fix a parameter $\alpha \geq \max (1,4c_{1}^{2}/c_{2})$. Next, we treat
separately the sets in the sum that do not contain at least
$\alpha c_{2} C_{\delta}(n)$ samples; for convenience, we denote by
$\tilde \ell ^{+}$ the largest index $j$ for which
$n\bar\mu _{n}\left ( \mathcal{V}\cap \{x\in \mathbb{R}^{d}: e^{\top }x
\geq z_{j}^{+}\} \right )\geq \alpha c_{2}C_{\delta}(n)$. Then, under
$\mathcal{E}_{n}(\delta )$, equation~\eqref{eq:uniform_concentration_H_k}
yields
\begin{align*}
\bar A^{+}_{n}(\mathcal{V};e)\leq \sum _{j=1}^{\tilde \ell ^{+}} (z_{j}^{+}-z_{j-1}^{+})
\left ( 1+\frac{2c_{1}}{\sqrt{c_{2}\alpha }} \right ) \mu \left (
\mathcal{V}\cap \{x\in \mathbb{R}^{d}: e^{\top }x \geq z_{j}^{+}\}
\right ) + z_{l^{+}}^{+} \cdot \alpha \frac{c_{2}C_{\delta}(n)}{n}.
\end{align*}
Now note that
\begin{align*}
\sum _{j=1}^{\tilde \ell ^{+}} (z_{j}^{+}-z_{j-1}^{+}) \mu \left (
\mathcal{V}\cap \{x\in \mathbb{R}^{d}: e^{\top }x \geq z_{j}^{+}\}
\right ) \leq \mathbb{E}\left [ e ^{\top }X \cdot \mathbbold{1}\{X
\in \mathcal{V}\} \mathbbold{1}\{e^{\top }X> 0\} \right ].
\end{align*}
On the other hand, equation~\eqref{eq:uniform_concentration_H_k} implies
that we have the following under $\mathcal{E}_{n}(\delta )$:
\begin{align*}
\bar A^{+}_{n}(\mathcal{V};e) =\sum _{j=1}^{\ell ^{+}} &(z_{j}^{+}-z_{j-1}^{+})
\cdot \bar\mu _{n}\left ( \mathcal{V}\cap \{x\in \mathbb{R}^{d}: e^{
\top }x > z_{j-1}^{+}\} \right )
\notag
\\
&\geq \sum _{j=1}^{\tilde \ell ^{+}} (z_{j}^{+}-z_{j-1}^{+})\left ( 1-
\frac{2c_{1}}{\sqrt{c_{2}\alpha }} \right ) \mu \left ( \mathcal{V}
\cap \{x\in \mathbb{R}^{d}: e^{\top }x > z_{j-1}^{+}\} \right ).
\end{align*}
Also, letting $y^{+}:= z^{+}_{\tilde \ell ^{+}}$ we have
\begin{align*}
\sum _{j=1}^{\tilde \ell ^{+}} (z_{j}^{+}-z_{j-1}^{+}) \mu \left (
\mathcal{V}\cap \{x\in \mathbb{R}^{d}: e^{\top }x > z_{j-1}^{+}\}
\right ) \geq \mathbb{E}\left [ e ^{\top }X \cdot \mathbbold{1}\{X
\in \mathcal{V}\} \mathbbold{1}\{0< e^{\top }X\leq y^{+}\} \right ].
\end{align*}
Therefore, defining
$A^{+}(\mathcal{V};e):= \mathbb{E}\left [ e ^{\top }X \cdot
\mathbbold{1}\{X\in \mathcal{V}\}\mathbbold{1}\{e^{\top }X > 0\}
\right ]$ and putting all four previous equations together shows that,
under $\mathcal{E}_{n}(\delta )$, we have
%
\begin{equation}
\label{eq:bound_A+_term_1}
\begin{split} &\left |\bar A^{+}_{n}(\mathcal{V};e) - A^{+}(
\mathcal{V};e) \right |
\\
&\leq \frac{4c_{1}}{\sqrt{c_{2}\alpha}} A^{+}(\mathcal{V};e) + \max _{1
\le i\le n} \|X_{i}\| \cdot \alpha \frac{c_{2}C_{\delta}(n)}{n} +
\mathbb{E}\left [ \|X\| \cdot \mathbbold{1}\{X\in \mathcal{V}\}
\mathbbold{1}\{e^{\top }X> y^{+}\} \right ],
\end{split}
%
\end{equation}
where we upper-bounded some terms $|e^{\top }X_{i}|$ with
$\|X_{i}\|$ for $1\le i \le n$ We continue bounding the terms on the right-hand
side. First, we recall by Lemma~\ref{lemma:charact_finite_obj} that
$\|X\|$ is stochastically dominated by $\sqrt{r}\cdot Z$, where
$Z\sim \textnormal{Par}(2)$. Then, Lemma~\ref{lemma:computation_bound_norm}
implies
%
\begin{align}
\label{eq:upper_bound_Pareto_argument}
A^{+}(\mathcal{V};e) \leq \mathbb{E}\left [ \|X\| \cdot \mathbbold{1}
\{X\in \mathcal{V}\} \right ] \leq 2\sqrt{r \mathbb{P}(X\in
\mathcal{V})}.
\end{align}
Next, by construction of $y^{+}$, the region
$\mathcal{V}\cap \{x:e^{\top }x>y^{+}\}$ contains at most
$\alpha c_{2}C_{\delta}(n)$ samples. Therefore, under
$\mathcal{E}_{n}(\delta )$, equation~\eqref{eq:main_H_event} implies that
\begin{equation*}
\mu \left ( \mathcal{V}\cap \{x\in \mathbb{R}^{d}:e^{\top }x>y^{+}\}
\right ) \leq c_{3} \frac{C_{\delta}(n)}{n},
\end{equation*}
for some universal constant $c_{3}>0$, and where $\mu $ denotes the law
of $X$. Then, we again use Lemma~\ref{lemma:computation_bound_norm} to
obtain
%
\begin{equation}
\label{eq:bound_edge_case}
\begin{split} &\mathbb{E}\left [ \|X\| \cdot \mathbbold{1}\{X\in
\mathcal{V}\}\mathbbold{1}\{e^{\top }X> y^{+}\} \right ]
\\
&\qquad \leq 2\sqrt{r \cdot \mu \left ( \mathcal{V}\cap \{x\in
\mathbb{R}^{d}:e^{\top }x>y^{+}\} \right )} \leq 2 \sqrt{
\frac{c_{3} r C_{\delta}(n)}{n}}.
\end{split}
%
\end{equation}
In summary, combining equations~\eqref{eq:bound_A+_term_1},~\eqref{eq:upper_bound_Pareto_argument},
and~\eqref{eq:bound_edge_case} we obtained
%
\begin{equation}
\label{eq:bound_A+_term}
\left |\bar A^{+}_{n}(\mathcal{V};e) - A^{+}(\mathcal{V};e) \right |
\leq 4c_{1}\sqrt{\frac{r \mathbb{P}(X\in \mathcal{V})}{c_{2}\alpha}} +
2 \sqrt{\frac{c_{3} r C_{\delta}(n)}{n}} + \max _{1\le i \le n} \|X_{i}
\| \cdot \alpha \frac{c_{2}C_{\delta}(n)}{n}.
\end{equation}
under the event $\mathcal{E}_{n}(\delta )$. If we define
$A^{-}(\mathcal{V};e):= \mathbb{E}\left [ -e ^{\top }X \cdot
\mathbbold{1}\{X\in \mathcal{V}\}\mathbbold{1}\{e^{\top }X < 0\}
\right ]$, then, by symmetry, the same arguments as above show that we
have
%
\begin{equation}
\label{eq:bound_A-_term}
\left |\bar A^{-}_{n}(\mathcal{V};e) - A^{-}(\mathcal{V};e) \right |
\leq 4c_{1}\sqrt{\frac{r\mathbb{P}(X\in \mathcal{V})}{c_{2}\alpha}} + 2
\sqrt{\frac{c_{3} r C_{\delta}(n)}{n}} + \max _{1\le i \le n} \|X_{i}
\| \cdot \alpha \frac{c_{2}C_{\delta}(n)}{n}.
\end{equation}
under $\mathcal{E}_{n}(\delta )$.

Now we turn to comparing $m_{X}(\mathcal{V})$ and
$\bar m_{n}(\mathcal{V})$, for fixed $\mathcal{V}\in \mathscr{V}$. Now
note that by construction and from equation~\eqref{eq:bar_m_expression_useful},
we have
%
\begin{align}
\label{eq:relate_to_original_sums}
\begin{cases}
A^{+}(\mathcal{V};e) - A^{-}(\mathcal{V};e) =\mathbb{E}\left [ e ^{
\top }X \cdot \mathbbold{1}\{X\in \mathcal{V}\} \right ] = \mathbb{P}(X
\in \mathcal{V})\cdot e ^{\top }m_{X}(\mathcal{V});
\\
\bar A^{+}_{n}(\mathcal{V};e) - \bar A^{-}_{n}(\mathcal{V};e) =
\bar \mu _{n}(\mathcal{V})\cdot e ^{\top }\bar m_{n}(\mathcal{V}).
\end{cases}
\end{align}
We next define
$e:=
\frac{\bar m_{n}(\mathcal{V}) - m_{X}(\mathcal{V}) }{\|\bar m_{n}(\mathcal{V}) - m_{X}(\mathcal{V}) \|}$
if $\bar m_{n}(\mathcal{V}) \neq m_{X}(\mathcal{V}) $ and otherwise fix
$e\in \mathbb{R}^{d}$ arbitrarily. Then, under
$\mathcal{E}_{n}(\delta )$ and for any $\alpha \geq 1$, we may use equations~\eqref{eq:lower_bound_nb_points_per_cluster},~\eqref{eq:relate_to_original_sums}, \eqref{eq:bound_A+_term}, and~\eqref{eq:bound_A-_term} to get
\begin{equation*}
|\bar\mu _{n}(\mathcal{V}) - \mu _{X}(\mathcal{V})|\leq
\frac{4c_{1}}{\sqrt n}
\frac{\mu _{X}(\mathcal{V})}{\sqrt{\bar\mu _{n}(\mathcal{V})}}\leq
\frac{4c_{1}}{\sqrt n} \sqrt{2\mu (\mathcal{V})}.
\end{equation*}
Hence
%
\begin{align}
\bar\mu _{n}(\mathcal{V}) &\left \|\bar m_{n}(\mathcal{V}) - m_{X}(
\mathcal{V}) \right \| =\left |\bar \mu _{n}(\mathcal{V}) \cdot e^{
\top }\left ( \bar m_{n}(\mathcal{V}) - m_{X}(\mathcal{V}) \right )
\right |
\notag
\\
&\leq \left |\bar \mu _{n}(\mathcal{V}) \cdot e^{\top }\bar m_{n}(
\mathcal{V}) - \mu (\mathcal{V})\cdot e ^{\top }m_{X}(\mathcal{V})
\right |+ |\bar\mu _{n}(\mathcal{V}) - \mu (\mathcal{V})|\cdot \|m_{X}(
\mathcal{V})\|
\notag
\\
&\leq 8c_{1}\sqrt{\frac{r \mu _{X}(\mathcal{V}) }{c_{2}\alpha}} + 4
\sqrt{\frac{c_{3} r C_{\delta}(n)}{n}} + \max _{1\le i \le n} \|X_{i}
\| \cdot \alpha \frac{2c_{2}C_{\delta}(n)}{n}
\notag
\\
&\qquad + 8c_{1} \sqrt{\frac{\mu _{X}(\mathcal{V})}{n}}\cdot \|m_{X}(
\mathcal{V})\|.
\label{eq:bound_deviation_m}
\end{align}
Now let $\gamma \geq \max (c_{2},4c_{1}^{2})C_{\delta}(n)$. Then, under
$\mathcal{E}_{n}(\delta )$, for any
$\mathscr{V}\in \bar{\mathscr P}_{n,k,\gamma}$ and $\alpha \geq 1$,
\begin{align*}
&\left | \bar D_{n}(\mathscr{V}) - D_{X}(\mathscr{V})\right |
\\
&\leq \sum _{\mathcal{V}\in \mathscr{V}} \left |\mu (\mathcal{V}) \|m_{X}(
\mathcal{V})\|^{2} - \bar \mu _{n}(\mathcal{V}) \|\bar m_{n}(
\mathcal{V})\|^{2}\right |
\\
&\leq \sum _{\mathcal{V}\in \mathscr{V}} \bigg(|\bar\mu _{n}(
\mathcal{V}) - \mu (\mathcal{V})| \cdot \|m_{X}(\mathcal{V})\|^{2}
\\
&\qquad \qquad + \bar\mu _{n}(\mathcal{V}) \|m_{X}(\mathcal{V})\|
\cdot \|\bar m_{n}(\mathcal{V}) - m_{X}(\mathcal{V})\| + \bar\mu _{n}(
\mathcal{V}) \|\bar m_{n}(\mathcal{V}) - m_{X}(\mathcal{V})\|^{2}
\bigg)
\end{align*}
Now use equation~\eqref{eq:bound_deviation_m}, the fact that
$ \mu (\mathcal{V})/2 \leq \bar\mu _{n}(\mathcal{V})\leq 2 \mu (
\mathcal{V})$ from \eqref{eq:lower_bound_nb_points_per_cluster}, and the
previous bound on $|\bar\mu _{n}(\mathcal{V}) - \mu (\mathcal{V})|$, to
get that there exists a universal constant $c_{4}>0$ with
\begin{multline*}
\Delta (\mathscr{V})\leq c_{4}\sum _{\mathcal{V}\in \mathscr{V}}
\left (\sqrt{\frac{\bar \mu _{n}(\mathcal{V})}{n}} \|m_{X}(
\mathcal{V})\|^{2} + \sqrt{r\left (
\frac{\bar \mu _{n}(\mathcal{V})}{\alpha } + \frac{C_{\delta }(n)}{n}
\right )}\|m_{X}(\mathcal{V})\|+ \frac{r}{\alpha} +
\frac{r C_{\delta}(n)}{n\mu (\mathcal{V})} \right .
\\
\left . + \frac{\|m_{X}(\mathcal{V})\|^{2}}{n} +\max _{1\le i \le n}
\|X_{i}\|\cdot \frac{\alpha C_{\delta}(n)}{n} \|m_{X}(\mathcal{V})\|+
\max _{1\le i\le n}\|X_{i}\|^{2}\cdot
\frac{\alpha ^{2} C_{\delta}(n)^{2}}{n^{2}\mu (\mathcal{V})} \right ).
\end{multline*}
Next, since $\|X\|$ is stochastically dominated by $\sqrt{r}\cdot Z$ where
$Z\sim \textnormal{Par}(2)$, Lemma~\ref{lemma:computation_bound_norm} implies
that for any $c\in \mathcal{C}$,
\begin{equation*}
\|m_{X}(\mathcal{V})\| \leq 2\sqrt{\frac{r}{\mu (\mathcal{V})}}.
\end{equation*}
Together with the identity
$ \mu (\mathcal{V})/2 \leq \bar\mu _{n}(\mathcal{V})\leq 2 \mu (
\mathcal{V})$, this implies that under $\mathcal{E}_{n}(\delta )$, letting
$Z_{n}^{\max}:=\max _{1\le i\le n}\|X_{i}\|/\sqrt{r}$, we have
%
\begin{align}
\Delta (\mathscr{V}) &\leq c_{5} r\sum _{\mathcal{V}\in \mathscr{V}}
\left ( \sqrt{\frac{C_{\delta }(n)}{n\bar \mu _{X}(\mathcal{V})}} +
\frac{C_{\delta }(n)}{n\bar \mu _{X}(\mathcal{V})} +\frac{1}{\alpha } +
\frac{\alpha C_{\delta }(n) Z_{n}^{\max }}{n\sqrt{\bar \mu _{n}(\mathcal{V})}}
+
\frac{\alpha ^{2} (C_{\delta }(n))^{2} (Z_{n}^{\max })^{2}}{n^{2}\bar \mu _{n}(\mathcal{V})}
\right )
\notag
\\
&\leq c_{6} kr\left ( \sqrt{\frac{C_{\delta }(n)}{ \gamma }} +
\frac{1}{\alpha } +
\frac{\alpha C_{\delta }(n) Z_{n}^{\max }}{\sqrt{\gamma n}} + \left (
\frac{\alpha C_{\delta }(n) Z_{n}^{\max }}{\sqrt{\gamma n}} \right )^{2}
\right ),
\label{eq:distortion_bound_v2}
\end{align}
for some universal constants $c_{5},c_{6}>0$. It now only remains to bound
$Z_{n}^{\max}$ which is stochastically dominated by the maximum of
$n$ i.i.d. Pareto random variables. To bound these, we introduce the following
event
\begin{equation*}
\mathcal{F}_{n}(\delta ):=\left \{ Z_{n}^{\max } \leq \sqrt{
\frac{n}{\delta }} \right \},
\end{equation*}
Since Lemma~\ref{lemma:charact_finite_obj} implies
$\mathbb{P}(\|X\|\geq \sqrt{r x})\leq 1/x$, we have
$\mathbb{P}(\mathcal{F}_{n}(\delta ))\geq 1-\delta $. In summary, under
$\mathcal{E}_{n}(\delta )\cap \mathcal{F}_{n}(\delta )$ which holds with
probability $1-2\delta $, we can further the bound from equation~\eqref{eq:distortion_bound_v2}
as follows. There exists universal constants $c_{8},c_{9}>0$ such that
for $\gamma \geq c_{8} C_{\delta}(n)^{2}/\delta $, using
$\alpha = (\gamma \delta /C_{\delta}(n)^{2})^{1/4}$, we have
\begin{align*}
\Delta (\mathscr{V}) \leq \frac{c_{9}}{2}kr \left ( \frac{1}{\alpha } +
\frac{\alpha C_{\delta }(n)}{\sqrt{\gamma \delta }} \right ) = c_{9} k
r \left ( \frac{(C_{\delta }(n))^{2}}{\gamma \delta } \right )^{1/4}
\end{align*}
for all $\mathscr{V} \in \bar{\mathscr P}_{n,k,\gamma}$. This ends the
proof of the first claim.

We next turn to almost-sure bounds. To do so, let
$\{\beta _{n}\}_{n\in \mathbb{N}}$ be any increasing and non-negative sequence
such that $\sum _{n\in \mathbb{N}} (n\beta _{n})^{-1}<\infty $. By assumption,
\begin{equation*}
\sum _{n\in \mathbb{N}} \mathbb{P}\left (\|X_{n}\| \geq \sqrt{r n
\beta _{n}}\right ) = \sum _{n\in \mathbb{N}} \frac{1}{n\beta _{n}}<
\infty .
\end{equation*}
So, Borel-Cantelli implies that, on an event $\mathcal{G}$ of probability
one, there exists $N$ such that for all $n\geq N$, we have
$\|X_{n}\|\leq \sqrt{rn\beta _{n}}$ and $\mathcal{E}_{n}(n^{-2})$ holds.
Note that under $\mathcal{G}$, there exists $N_{1}$ such that we have
\begin{equation*}
Z_{n}^{\max} \leq \sqrt{n\beta _{n}} +Z_{\hat n}^{\max} \leq 2\sqrt{n
\beta _{n}},\qquad \textnormal{ for all } n\geq N_{1}.
\end{equation*}
Therefore, there exists universal constants $c_{10},c_{11}>0$ such that
for any sequence $\gamma _{n}\geq \beta _{n}(\log n)^{2}$, using
$\alpha _{n} = \gamma _{n}^{1/4}\beta _{n}^{-1/4}(dk\log k \log n)^{-1/2}$,
the following holds. Under $\mathcal{G}$, for $n$ sufficiently large, we
have
\begin{equation*}
\Delta (\mathscr{V}) \leq c_{10} kr \left ( \frac{1}{\alpha _{n}} +
\frac{\alpha _{n} C_{1/n^{2}}(n)\sqrt{\beta _{n}}}{\sqrt{\gamma _{n}}}
\right ) \leq c_{11} k r \left (
\frac{\beta _{n} (dk\log k \log n)^{2}}{\gamma _{n}} \right )^{1/4}
\end{equation*}
for all $\mathscr{V} \in \bar{\mathscr P}_{n,k,\gamma _{n}}$. Up to dividing
$\beta _{n}$ by a factor $(dk\log k)^{2}$ (which still preserves finiteness
of $\sum _{n\geq 1}(n\beta _{n})^{-1}$, this ends the proof.
\end{proof}

Another technical result we needed along the way is the following deterministic
result which is used in the proof of Theorem~\ref{thm:consistency_unified}.

\begin{proof}[Proof of Lemma~\ref{lemma:can_delete_diverging_clusters}]
In this proof, we denote by
$B_{0}(M):=\{x\in \mathbb{R}^{d}: \|x\|\leq M\}$ the centered ball of radius
$M$. By continuity of the distortion objective (Lemma~\ref{lem:D-cts}),
for any $M\geq 1$ and $1\le k'\le k$, and for any $\epsilon >0$ there exists
$\eta (M;\epsilon ,k')\in (0,1]$ such that for any
$\mathcal{C}\in \mathscr{C}_{k'}$ with
$\mathcal{C}\subseteq B_{0}(M)$,
\begin{equation*}
\min _{\mathcal{C}'\in \mathscr{C}_{k'}(X)}d_{\textnormal{H}}(
\mathcal{C},\mathcal{C}')\geq \epsilon \implies D_{X}(\mathcal{C}\,|
\,\mathcal{C}_{0}) \geq \inf _{\mathcal{C}'\in \mathscr{C}_{k'}} D_{X}(
\mathcal{C}'\,|\,\mathcal{C}_{0}) + \eta (M;\epsilon ,k').
\end{equation*}
We denote
$\eta (M;\epsilon ):=\min _{1\le k'\le k}\eta (M;\epsilon ,k')/2$. Then,
the previous equation implies that there exists
$\phi (M;\epsilon )\geq M$ such that for any $1\le k'\le k$ and
$\mathcal{C}\in \mathscr{C}_{k'}$ with
$\mathcal{C}\subseteq B_{0}(M)$, the condition
\begin{equation*}
\min _{\mathcal{C}'\in \mathscr{C}_{k'}(X)}d_{\textnormal{H}}(
\mathcal{C},\mathcal{C}')\geq \epsilon
\end{equation*}
implies that there exists $\tilde{\mathcal{C}}\in \mathscr{C}_{k'}$ satisfying
$\tilde{\mathcal{C}}\subseteq B_{0}(\phi (M;\epsilon ))$ and
%
\begin{equation}
\label{eq:construction_better_cluster_candidate}
D_{X}(\mathcal{C}\,|\,\mathcal{C}_{0}) \geq D_{X}(
\tilde {\mathcal{C}}\,|\,\mathcal{C}_{0})+ \eta (M;\epsilon ).
\end{equation}
Without loss of generality, we may also assume that
$\eta (M;\epsilon )$ is non-increasing in $M$. Next, we define
$\psi (M;\epsilon )\geq \phi (M;\epsilon )+1$ such that denoting by
$\mu $ the law of $X$:
%
\begin{equation}
\label{eq:def_psi_function}
\mu \left ( \mathbb{R}^{d}\setminus B_{0}\left (
\frac{\psi (M;\epsilon )-\phi (M;\epsilon )-1}{2} \right ) \right )
\leq \frac{(\eta (M;\epsilon )/3)^{2}}{20r\phi (M;\epsilon )^{2}}.
\end{equation}
Again, without loss of generality, assume that $\psi (M;\epsilon )$ is
non-increasing in $M$.

Note that for any set of points $\mathcal{A}\in \mathscr{C}_{k}$ and measurable
region $B\subseteq \mathbb{R}^{d}$, we have
\begin{align*}
&\left |\mathbb{E}\left [ \min _{c\in \mathcal{A}}(\|X-c\|^{2}-\|X\|^{2})
\mathbbold{1}\{X\notin B\} \right ] - D_{X}(\mathcal{A}\,|\,
\mathcal{C}_{0})\right |
\\
&\qquad \leq \max _{a\in \mathcal{A}}\|a\|^{2} \mathbb{P}(X\in B) + 2
\max _{a\in \mathcal{A}}\|a\|\cdot \mathbb{E}\left [ \|X\|
\mathbbold{1}\{X\in B\} \right ]
\notag
\\
&\qquad \leq \max _{a\in \mathcal{A}}\|a\|^{2} \mathbb{P}(X\in B) + 4
\max _{a\in \mathcal{A}}\|a\|\cdot \sqrt{r \mathbb{P}(X\in B)}.
\end{align*}
In the last inequality, we used the fact that $\|X\|$ is stochastically
dominated by $\sqrt{r} \cdot Z$ where $Z\sim \textnormal{Par}(2)$ from
Lemma~\ref{lemma:charact_finite_obj} and Lemma~\ref{lemma:computation_bound_norm}.
Without loss of generality, we suppose $r\geq 1$. In particular, for any
$\mathcal{A}\sqcup \mathcal{B}\in \mathcal{C}_{k}$ such that
$\mathcal{A}\subseteq B_{0}(M)$ and
$\mathbb{P}(X\in B)\leq \eta ^{2}/(20rM^{2})$ for some
$\eta \in (0,1]$, we have
%
\begin{equation}
\label{eq:comparison_deleting_far_clusters}
\left |\mathbb{E}\left [ \min _{c\in \mathcal{A}}(\|X-c\|^{2}-\|X\|^{2})
\mathbbold{1}\{X\notin B\} \right ] - D_{X}(\mathcal{A}\,|\,
\mathcal{C}_{0}) \right | \leq \eta .
\end{equation}
Now let $\mathcal{C}\in \mathscr{C}_{k}$. We suppose that we can decompose
$\mathcal{C}=\mathcal{A}\sqcup \mathcal{B}$ with
$\mathcal{A}\neq \emptyset $ such that letting $k'=\#\mathcal{A}$, we have
%
\begin{equation}
\label{eq:assumptions_intermediary}
\min _{\mathcal{C}'\in \mathscr{C}_{k'}(X)} d_{\textnormal{H}}(
\mathcal{A},\mathcal{C}')\geq \epsilon ,\quad \max _{c\in \mathcal{A}}
\|c\|\leq M,\quad \text{and} \quad \min _{c\in \mathcal{B}} \|c\|\geq
\psi (M;\epsilon ).
\end{equation}
Since $\mathcal{A}\neq \emptyset $, we can fix $a\in \mathcal{A}$. Next,
let $c\in \mathcal{B}$. In particular, by
\eqref{eq:assumptions_intermediary} we have $\|a\|\leq M$ and
$\|c\|\geq \psi (M;\epsilon )$. Therefore,
\begin{align*}
\mathcal{V}_{\mathcal{C}}(c)&\subseteq \{x\in \mathbb{R}^{d}: \|x-c\|
\leq \|x-a\|\}
\\
&\subseteq \{x\in \mathbb{R}^{d}: \psi (M;\epsilon ) - \|x\|\leq \|x
\|+M\} \subseteq \mathbb{R}^{d}\setminus B_{0}\left (
\frac{\psi (M;\epsilon )-M-1}{2} \right ).
\end{align*}
Hence, using \eqref{eq:def_psi_function} and the fact that
$\phi (M;\epsilon )\geq M$, we obtained
\begin{equation*}
\mu \left ( \bigcup _{b\in \mathcal{B}}\mathcal{V}_{\mathcal{C}}(b)
\right ) \leq \mu \left ( \mathbb{R}^{d}\setminus B_{0}\left (
\frac{\psi (M;\epsilon )-\phi (M;\epsilon )-1}{2} \right ) \right )
\leq \frac{(\eta (M;\epsilon )/3)^{2}}{20 r M^{2}}.
\end{equation*}
Then, equation~\eqref{eq:comparison_deleting_far_clusters} applied to
$\mathcal{A}$ and
$B:=\bigcup _{b\in \mathcal{B}}\mathcal{V}_{\mathcal{C}}(b)$ shows that
\begin{align*}
D_{X}(\mathcal{C}\,|\,\mathcal{C}_{0}) &- \mathbb{E}\left [ \sum _{b
\in \mathcal{B}}(\|X-b\|^{2} -\|X\|^{2})\mathbbold{1}\{X\in
\mathcal{V}_{\mathcal{C}}(b)\} \right ] - D_{X}( \mathcal{A}\,|\,
\mathcal{C}_{0})
\\
&=\mathbb{E}\left [ \sum _{a\in \mathcal{A}}(\|X-a\|^{2}-\|X\|^{2})
\mathbbold{1}\{X\in \mathcal{V}_{\mathcal{C}}(a)\} \right ] - D_{X}(
\mathcal{A}\,|\,\mathcal{C}_{0})
\\
&\geq \mathbb{E}\left [ \min _{c\in \mathcal{A}}(\|X-c\|^{2}-\|X\|^{2})
\mathbbold{1}\{X\notin B\} \right ] - D_{X}( \mathcal{A}\,|\,
\mathcal{C}_{0})
\\
&\geq -\frac{\eta (M;\epsilon )}{3}.
\end{align*}
Next, from equation~\eqref{eq:construction_better_cluster_candidate} there
exists $\tilde{\mathcal{A}}\in \mathscr{C}_{k'}$ such that
$\tilde{\mathcal{A}}\subseteq B_{0}(\phi (M;\epsilon ))$ and
$D_{X}(\mathcal{A}\,|\,\mathcal{C}_{0}) \geq D_{X}(
\tilde{\mathcal{A}}\,|\,\mathcal{C}_{0})+ \eta (\epsilon ;M)$. We now define
$\tilde{\mathcal{C}}:=\tilde{\mathcal{A}}\cup \mathcal{B}$. Applying equation~\eqref{eq:comparison_deleting_far_clusters}
to $\tilde{\mathcal{A}}$ and the same set
$B\subseteq \tilde{\mathcal{C}}$, we get
\begin{align*}
D_{X}(\tilde {\mathcal{C}}\,|\,\mathcal{C}_{0}) &- \mathbb{E}\left [
\sum _{b\in \mathcal{B}}(\|X-b\|^{2} -\|X\|^{2})\mathbbold{1}\{X\in
\mathcal{V}_{\mathcal{C}}(b)\} \right ] - D_{X}( \tilde{\mathcal{A}}
\,|\,\mathcal{C}_{0})
\\
&= \mathbb{E}\left [ \min _{c\in \tilde {\mathcal{C}}}(\|X-c\|^{2}-\|X
\|^{2}) \right ] - \mathbb{E}\left [ \sum _{b\in \mathcal{B}}(\|X-b\|^{2}
-\|X\|^{2})\mathbbold{1}\{X\in \mathcal{V}_{\mathcal{C}}(b)\} \right ]
\\
&\qquad - D_{X}( \tilde{\mathcal{A}}\,|\,\mathcal{C}_{0})
\\
&\leq \mathbb{E}\left [ \min _{c\in \tilde {\mathcal{C}}}(\|X-c\|^{2}-
\|X\|^{2}) \mathbbold{1}\{X\notin B\} \right ] - D_{X}(
\tilde{\mathcal{A}}\,|\,\mathcal{C}_{0})
\\
&\leq \mathbb{E}\left [ \min _{c\in \tilde {\mathcal{A}}}(\|X-c\|^{2}-
\|X\|^{2}) \mathbbold{1}\{X\notin B\} \right ] - D_{X}(
\tilde{\mathcal{A}}\,|\,\mathcal{C}_{0}) \leq
\frac{\eta (M;\epsilon )}{3}.
\end{align*}
In summary, we obtained
%
\begin{align}
\label{eq:conclusion_for_recursion}
D_{X}(\mathcal{C}\,|\,\mathcal{C}_{0}) - \inf _{\mathcal{C}'\in
\mathscr{C}_{k}}D_{X}(\mathcal{C}'\,|\,\mathcal{C}_{0}) &\geq D_{X}(
\mathcal{C}\,|\,\mathcal{C}_{0}) - D_{X}(\tilde {\mathcal{C}}\,|\,
\mathcal{C}_{0})
\notag
\\
&\geq D_{X}(\mathcal{A}\,|\,\mathcal{C}_{0}) - D_{X}(
\tilde{\mathcal{A}}\,|\,\mathcal{C}_{0}) -
\frac{2\eta (M;\epsilon )}{3}\geq \frac{\eta (M;\epsilon )}{3}.
\end{align}
by combining the two previous inequalities.

With this identity, we are now ready to show the following for some arbitrary
$\epsilon >0$. Fix a sequence of cluster sets
$\{\mathcal{C}_{n}\}_{n\in \mathbb{N}}$ satisfying the assumptions for
this result. We aim to show that there is a decomposition of
$\mathcal{C}_{n} = \mathcal{A}_{n}(\epsilon )\sqcup \mathcal{B}_{n}(
\epsilon )$ for $n\in \mathbb{N}$ so that for sufficiently large $n$,
%
\begin{equation}
\label{eq:desired_decomposition_epsilon}
\min _{1\le k'\le k}\min _{\mathcal{C}\in \mathscr{C}_{k'}(X)}d_{
\textnormal{H}}(\mathcal{A}_{n}(\epsilon ),\mathcal{C})\leq \epsilon
\quad \text{and} \quad \min _{c\in \mathcal{B}_{n}(\epsilon )}\|c\|
\geq \frac{1}{\epsilon}.
\end{equation}
First, by assumption, there exists $n_{\epsilon}$ such that for all
$n\geq n_{\epsilon}$, we have
\begin{equation*}
D_{X}(\mathcal{C}_{n}\,|\,\mathcal{C}_{0}) \leq \inf _{\mathcal{C}
\in \mathscr{C}_{k}} D_{X}(\mathcal{C}\,|\,\mathcal{C}_{0}) +
\frac{1}{4}\eta \left ( \psi ^{(k)}(1/\epsilon ;\epsilon ) ;\epsilon
\right ),
\end{equation*}
where we used the notation where $\psi ^{(0)}$ is the identity function
and
$\psi ^{({\ell}+1)}(M;\epsilon ):=\break \psi \left ( \psi ^{(\ell )}(M;
\epsilon );\epsilon \right )$ for $\ell \in \mathbb{N}$ and any
$M\geq 0$. For convenience, we denote
\begin{equation*}
M_{\ell}:= \psi ^{(\ell )}(1/\epsilon ;\epsilon )\,\quad 1\le \ell
\le k.
\end{equation*}
We can check that since
$\psi (M;\epsilon )\geq \phi (M;\epsilon )\geq M$, for all $M\geq 0$, we
have
$M_{k}\geq M_{k-1}\geq \cdots \geq M_{1}\geq M_{0}:=1/\epsilon $. Now fix
$n\geq n_{\epsilon}$. First, if we have
$\min _{c\in \mathcal{C}_{n}} \|c\| \geq M_{0}$ then, the desired result
\eqref{eq:desired_decomposition_epsilon} is immediate by taking
$\mathcal{B}_{n}(\epsilon )=\mathcal{C}_{n}$. We now suppose that there
is at least one element $c\in \mathcal{C}_{n}$ with $\|c\|<M_{0}$. Then,
by the pigeonhole principle, there must be some $1\le \ell \le k$ for which
\begin{equation*}
\#\left \{ c\in \mathcal{C}_{n}: \|c\| \in (M_{\ell -1},M_{\ell })
\right \} = 0.
\end{equation*}
We now pose
\begin{equation*}
\mathcal{A}_{n}(\epsilon ):= \mathcal{C}_{n}\cap B_{0}(M_{\ell -1})
\quad \text{and} \quad \mathcal{B}_{n}(\epsilon ) = \mathcal{C}_{n}
\setminus B_{0}(M_{\ell -1}).
\end{equation*}
By construction, we have $\mathcal{A}_{n}(\epsilon )\neq \emptyset $ and
\begin{align*}
&\max _{c\in \mathcal{A}_{n}(\epsilon )}\|c\| \leq M_{\ell -1}, \quad
\min _{c\in \mathcal{B}_{n}(\epsilon )} \|c\| \geq M_{\ell}, \quad \text{and}
\\
& D_{X}(\mathcal{C}_{n}\,|\,\mathcal{C}_{0}) \leq \inf _{
\mathcal{C}\in \mathscr{C}_{k}}D_{X}(\mathcal{C}\,|\,\mathcal{C}_{0}) +
\frac{\eta (M_{\ell -1};\epsilon )}{4}.
\end{align*}
Also, by construction,
$M_{\ell} = \psi (\epsilon ;M_{\ell -1})\geq M_{\ell -1} \geq M_{0} =1/
\epsilon $. Hence, equations~\eqref{eq:assumptions_intermediary} and~\eqref{eq:conclusion_for_recursion}
imply that
\begin{equation*}
\min _{\mathcal{C}\in \mathscr{C}_{k'}(X)} d_{\textnormal{H}}(
\mathcal{A}_{n}(\epsilon ),\mathcal{C}) \leq \epsilon ,
\end{equation*}
where $k'=\#\mathcal{A}_{n}(\epsilon )$. This ends the proof of~\eqref{eq:desired_decomposition_epsilon}.

Because this holds for any $\epsilon >0$, this implies that there exists
a decomposition
$\mathcal{C}_{n}=\mathcal{A}_{n}\sqcup \mathcal{B}_{n}$ for which
\begin{equation*}
\min _{1\le k'\le k}\min _{\mathcal{C}\in \mathscr{C}_{k'}(X)}d_{
\textnormal{H}}(\mathcal{A}_{n},\mathcal{C})\to 0 \quad \text{and}
\quad \min _{c\in \mathcal{B}_{n}}\|c\|\to \infty
\end{equation*}
as $n\to \infty $. It remains only to check that we can take
$\mathcal{A}_{n}$ such that $\mathcal{A}_{n}\neq \varnothing $. To prove
that this holds, it suffices to show that
$\min _{c\in \mathcal{C}_{n}}\|c\|$ remains bounded. Indeed, there exists
$c_{n}\in \mathcal{C}_{n}$ for which
$\mu (\mathcal{V}_{\mathcal{C}_{n}}(c_{n}))\geq 1/k$. From Lemma~\ref{lemma:computation_bound_norm},
we have
$\|m_{X}(\mathcal{V}_{\mathcal{C}_{n}}(c_{n}))\| \leq 2\sqrt{kr}$. Next,
consider the cluster centers
$\mathcal{C}'_{n}=(\mathcal{C}_{n}\setminus \{c_{n}\})\cup \{m_{X}(
\mathcal{V}_{\mathcal{C}_{n}}(c_{n}))\}$, where we replaced the point
$c_{n}$ with $m_{X}(\mathcal{V}_{\mathcal{C}_{n}}(c_{n}))$. Then,
\begin{align*}
D_{X}(\mathcal{C}_{n}\,|\,\mathcal{C}_{0}) -D_{X}(\mathcal{C}_{n}'\,|
\,\mathcal{C}_{0}) &= \mathbb{E}\left [ \min _{c\in \mathcal{C}_{n}}
\|X-c\|^{2} - \min _{c\in \mathcal{C}_{n}'}\|X-c\|^{2} \right ]
\\
&\geq \mathbb{E}\left [ \left ( \|X-c_{n}\|^{2} - \|X-m_{X}(
\mathcal{V}_{\mathcal{C}_{n}}(c_{n}))\|^{2} \right ) \mathbbold{1}\{X
\in \mathcal{V}_{\mathcal{C}_{n}}(c_{n})\} \right ]
\\
&=\mu (\mathcal{V}_{\mathcal{C}_{n}}(c_{n}))\|c_{n}-m_{X}(\mathcal{V}_{
\mathcal{C}_{n}}(c_{n}))\|^{2}
\\
&\geq
\frac{\|c_{n}-m_{X}(\mathcal{V}_{\mathcal{C}_{n}}(c_{n}))\|^{2}}{k}.
\end{align*}
As a result, we obtained
\begin{align*}
\min _{c\in \mathcal{C}_{n}}\|c\| &\leq \|m_{X}(\mathcal{V}_{
\mathcal{C}_{n}}(c_{n}))\| + \|c_{n}-m_{X}(\mathcal{V}_{\mathcal{C}_{n}}(c_{n}))
\|
\\
&\leq 2\sqrt{kr} + \sqrt{k\left ( D_{X}(\mathcal{C}_{n}\,|\,
\mathcal{C}_{0}) -\inf _{\mathcal{C}\in \mathscr{C}_{k}} D_{X}(
\mathcal{C}\,|\,\mathcal{C}_{0}) \right )},
\end{align*}
which remains bounded since the right-hand side converges to
$2\sqrt{kr}$ as $n\to \infty $. This ends the proof.
\end{proof}

\bibliography{refs}

@inproceedings{blanchard2024tight,
	title={Tight bounds for local {G}livenko-{C}antelli},
	author={Blanchard, M. and Voracek, V.},
	booktitle={ALT},
	pages={179--220},
	year={2024},
	organization={PMLR}
}

@article{boucheron2005theory,
	title={Theory of classification: {A} survey of some recent advances},
	author={Boucheron, S. and Bousquet, O. and Lugosi, G.},
	journal={ESAIM: P\&S},
	volume={9},
	pages={323--375},
	year={2005}
}

@article{Pollard,
	title={Strong consistency of k-means clustering},
	author={Pollard, D.},
	journal={Ann. Statist.},
	pages={135--140},
	year={1981},
	publisher={JSTOR}
}

@book{Villani,
	title={Optimal Transport: Old and New},
	author={Villani, C.},
	isbn={9783540710509},
	lccn={2008932183},
	series={Grundlehren der mathematischen Wissenschaften},
	year={2008},
	publisher={Springer Berlin Heidelberg}
}

@article {Thorpe,
	AUTHOR = {Thorpe, M. and Theil, F. and Johansen, A. M. and
		Cade, N.},
	TITLE = {Convergence of the {$k$}-means minimization problem using
		{$\Gamma$}-convergence},
	JOURNAL = {SIAM J. Appl. Math.},
	FJOURNAL = {SIAM Journal on Applied Mathematics},
	VOLUME = {75},
	YEAR = {2015},
	NUMBER = {6},
	PAGES = {2444--2474},
	ISSN = {0036-1399},
	MRCLASS = {62G20 (62H30)},
	MRREVIEWER = {Hiroshi Kurata},
	DOI = {10.1137/140974365},
	URL = {https://doi.org/10.1137/140974365},
}

@article {Lember,
	AUTHOR = {Lember, J.},
	TITLE = {On minimizing sequences for {$k$}-centres},
	JOURNAL = {J. Approx. Theory},
	FJOURNAL = {Journal of Approximation Theory},
	VOLUME = {120},
	YEAR = {2003},
	NUMBER = {1},
	PAGES = {20--35},
	ISSN = {0021-9045},
	MRCLASS = {41A65 (41A50 60B05)},
	MRREVIEWER = {Carlos Matr\'{a}n},
	DOI = {10.1016/S0021-9045(02)00010-2},
	URL = {https://doi.org/10.1016/S0021-9045(02)00010-2},
}

@article {Parna1,
	AUTHOR = {P\"{a}rna, K.},
	TITLE = {Strong consistency of {$k$}-means clustering criterion in
		separable metric spaces},
	JOURNAL = {Tartu Riikl. \"{U}l. Toimetised},
	FJOURNAL = {Tartu Riikliku \"{U}likooli Toimetised. Uchenye Zapiski Tartuskogo
		Gosudarstvennogo Universiteta. Acta et Commentationes
		Universitatis Tartuensis},
	NUMBER = {733},
	YEAR = {1986},
	PAGES = {86--96},
	MRCLASS = {62H30},
}

@article {Parna2,
	AUTHOR = {P\"{a}rna, K.},
	TITLE = {On the stability of {$k$}-means clustering in metric spaces},
	JOURNAL = {Tartu Riikl. \"{U}l. Toimetised},
	FJOURNAL = {Tartu Riikliku \"{U}likooli Toimetised. Uchenye Zapiski Tartuskogo
		Gosudarstvennogo Universiteta. Acta et Commentationes
		Universitatis Tartuensis},
	NUMBER = {798},
	YEAR = {1988},
	PAGES = {19--36},
	MRCLASS = {62H30},
	MRREVIEWER = {Ferenc Juh\'{a}sz},
}

@article {Laloe,
	AUTHOR = {Lalo\"{e}, T.},
	TITLE = {{$L_1$}-quantization and clustering in {B}anach spaces},
	JOURNAL = {Math. Methods Statist.},
	FJOURNAL = {Mathematical Methods of Statistics},
	VOLUME = {19},
	YEAR = {2010},
	NUMBER = {2},
	PAGES = {136--150},
	ISSN = {1066-5307},
	MRCLASS = {62H30 (60B11 62G05)},
	MRREVIEWER = {Zailong Wang},
	DOI = {10.3103/S1066530710020031},
	URL = {https://doi.org/10.3103/S1066530710020031},
}

@article {Zhivotovskiy,
	AUTHOR = {Klochkov, Y. and Kroshnin, A. and Zhivotovskiy, N.},
	TITLE = {Robust {$k$}-means clustering for distributions with two
		moments},
	JOURNAL = {Ann. Statist.},
	FJOURNAL = {The Annals of Statistics},
	VOLUME = {49},
	YEAR = {2021},
	NUMBER = {4},
	PAGES = {2206--2230},
	ISSN = {0090-5364},
	MRCLASS = {62F35 (60B11 62F12 62H30)},
	MRREVIEWER = {Yan Liu},
	DOI = {10.1214/20-aos2033},
	URL = {https://doi.org/10.1214/20-aos2033},
}

@article{JaffeClustering,
	author = {A. Q. Jaffe},
	doi = {10.1214/25-AOS2514},
	journal = {Ann. Statist.},
	number = {4},
	pages = {1559 -- 1586},
	publisher = {Institute of Mathematical Statistics},
	title = {{Asymptotic theory of geometric and adaptive k-means clustering}},
	url = {https://doi.org/10.1214/25-AOS2514},
	volume = {53},
	year = {2025}
}

@article {PollardCLT,
	AUTHOR = {D. Pollard},
	TITLE = {A central limit theorem for {$k$}-means clustering},
	JOURNAL = {Ann. Probab.},
	FJOURNAL = {The Annals of Probability},
	VOLUME = {10},
	YEAR = {1982},
	NUMBER = {4},
	PAGES = {919--926},
	ISSN = {0091-1798,2168-894X},
	MRCLASS = {60F05 (62E20 62H30)},
	MRREVIEWER = {Albert\ Raugi},
	URL =
	{http://links.jstor.org/sici?sici=0091-1798(198211)10:4<919:ACLTFC>2.0.CO;2-Q&origin=MSN},
}

@book{GrafLuschgy,
	title={Foundations of Quantization for Probability Distributions},
	author={Graf, S. and Luschgy, H.},
	isbn={9783540455776},
	series={Lecture Notes in Mathematics},
	url={https://books.google.com/books?id=0KZ7CwAAQBAJ},
	year={2007},
	publisher={Springer Berlin Heidelberg}
}

@article{BiauDevroyeLugosi,
	author = {Biau, G. and Devroye, L. and Lugosi, G.},
	year = {2008},
	month = {03},
	pages = {781--790},
	title = {On the performance of clustering in {H}ilbert spaces},
	volume = {54},
	journal = {IEEE Trans. Inform. Theory},
	doi = {10.1109/TIT.2007.913516}
}

@ARTICLE{BartlettLinderLugosi,
	author={Bartlett, P. L. and Linder, T. and Lugosi, G.},
	journal={IEEE Trans. Inform. Theory}, 
	title={The minimax distortion redundancy in empirical quantizer design}, 
	year={1998},
	volume={44},
	number={5},
	pages={1802--1813},
	doi={10.1109/18.705560}
}

@article{Levrard,
	author = {C. Levrard},
	doi = {10.1214/14-AOS1293},
	journal = {Ann. Statist.},
	number = {2},
	pages = {592--619},
	publisher = {Institute of Mathematical Statistics},
	title = {{Nonasymptotic bounds for vector quantization in Hilbert spaces}},
	url = {https://doi.org/10.1214/14-AOS1293},
	volume = {43},
	year = {2015}
}

@book {BuragoBuragoIvanov,
	AUTHOR = {Burago, D. and Burago, Y. and Ivanov, S.},
	TITLE = {A {C}ourse in {M}etric {G}eometry},
	SERIES = {Graduate Studies in Mathematics},
	VOLUME = {33},
	PUBLISHER = {American Mathematical Society, Providence, RI},
	YEAR = {2001},
	PAGES = {xiv+415},
	ISBN = {0-8218-2129-6},
	MRCLASS = {53C23},
	MRREVIEWER = {Mario Bonk},
	DOI = {10.1090/gsm/033},
	URL = {https://doi.org/10.1090/gsm/033},
}

@InProceedings{JiangAriasCastro,
	title = 	 { On the Consistency of Metric and Non-Metric K-Medoids },
	author =       {Jiang, H. and Arias-Castro, E.},
	booktitle = 	 {AISTATS},
	pages = 	 {2485--2493},
	year = 	 {2021},
	organization={PMLR}
}

@article{Fefferman,
	author = {Fefferman, C. and Mitter, S. and Narayanan, H.},
	year = {2016},
	pages = {983--1049},
	title = {Testing the manifold hypothesis},
	volume = {29},
	journal = {J. Amer. Math. Soc.},
	doi = {10.1090/jams/852}
}

@article{Ng,
	author = {M. K. Ng},
	doi = {https://doi.org/10.1016/S0031-3203(99)00057-6},
	issn = {0031-3203},
	journal = {Pattern Recognit.},
	number = {3},
	pages = {515--519},
	title = {A note on constrained k-means algorithms},
	url = {https://www.sciencedirect.com/science/article/pii/S0031320399000576},
	volume = {33},
	year = {2000}
}

@InProceedings{CuturiDoucet,
	title = 	 {Fast Computation of {W}asserstein Barycenters},
	author = 	 {Cuturi, M. and Doucet, A.},
	booktitle = 	 {ICML},
	pages = 	 {685--693},
	year = 	 {2014},
	organization={PMLR}
}

@article{GenevayDulacArnoldVert,
	title={Differentiable Deep Clustering with Cluster Size Constraints}, 
	author={Genevay, A. and Dulac-Arnold , G. and Vert, J.-P.},	
	year={2019},
	journal={arXiv preprint arXiv:1910.09036},
	
}

@article{appert2021new,
	title={New bounds for $ k $-means and information $ k $-means},
	author={Appert, G. and Catoni, O.},
	journal={arXiv preprint arXiv:2101.05728},
	year={2021}
}

@article{blumer1989learnability,
	title={Learnability and the {V}apnik-{C}hervonenkis dimension},
	author={Blumer, A. and Ehrenfeucht, A. and Haussler, D. and Warmuth, M. K.},
	journal={JACM},
	volume={36},
	number={4},
	pages={929--965},
	year={1989},
	publisher={ACM New York, NY, USA}
}

@book{Bishop,
	author = {Bishop, C. M.},
	cn = {Q327 B52 2006; Q327 .B52 2006},
	id = {9910032530902121},
	publisher = {Springer New York},
	title = {Pattern {R}ecognition and {M}achine {L}earning},
	url = {https://search.library.wisc.edu/catalog/9910032530902121},
	year = {{$[$}2006{$]$}}}

@article{FrechetMeanInfDim,
	title={Fr\'echet means in infinite dimensions}, 
	author={A. Q. Jaffe},	
	year={2025},
	journal={arXiv preprint arXiv:2410.17214},
}

@article {SchoetzSLLN,
	AUTHOR = {Sch\"otz, C.},
	TITLE = {Strong laws of large numbers for generalizations of
		{F}r\'echet mean sets},
	JOURNAL = {Statistics},
	FJOURNAL = {Statistics. A Journal of Theoretical and Applied Statistics},
	VOLUME = {56},
	YEAR = {2022},
	NUMBER = {1},
	PAGES = {34--52},
	ISSN = {0233-1888,1029-4910},
	MRCLASS = {60F15},
	DOI = {10.1080/02331888.2022.2032063},
	URL = {https://doi.org/10.1080/02331888.2022.2032063},
}

@article{TrimmedKMeansII,
	ISSN = {01621459, 1537274X},
	URL = {http://www.jstor.org/stable/2670010},
	author = {L. A. García-Escudero and A. Gordaliza},
	journal = {J. Am. Stat. Assoc.},
	number = {447},
	pages = {956--969},
	publisher = {[American Statistical Association, Taylor & Francis, Ltd.]},
	title = {Robustness Properties of k Means and Trimmed k Means},
	urldate = {2025-06-22},
	volume = {94},
	year = {1999}
}

@article{TrimmedKMeansI,
	author = {J. A. Cuesta-Albertos and A. Gordaliza and C. Matr{\'a}n},
	title = {{Trimmed $k$-means: an attempt to robustify quantizers}},
	volume = {25},
	journal = {Ann. Statist.},
	number = {2},
	publisher = {Institute of Mathematical Statistics},
	pages = {553--576},
	year = {1997},
	doi = {10.1214/aos/1031833664},
	URL = {https://doi.org/10.1214/aos/1031833664}
}

@article{RobustTrimmed,
	author = {O. Dorabiala and J. N. Kutz and A. Y. Aravkin},
	doi = {https://doi.org/10.1016/j.patrec.2022.07.007},
	issn = {0167-8655},
	journal = {Pattern Recognit. Lett.},
	pages = {9--16},
	title = {Robust trimmed k-means},
	url = {https://www.sciencedirect.com/science/article/pii/S0167865522002185},
	volume = {161},
	year = {2022}}

@article{RobustClusteringSurvey,
	author = {Garc{\'\i}a-Escudero, L. A. and Gordaliza, A. and Matr{\'a}n, C. and Mayo-Iscar, A.},
	date = {2010/09/01},
	doi = {10.1007/s11634-010-0064-5},
	id = {Garc{\'\i}a-Escudero2010},
	isbn = {1862-5355},
	journal = {Adv. Data Anal. Classif.},
	number = {2},
	pages = {89--109},
	title = {A review of robust clustering methods},
	url = {https://doi.org/10.1007/s11634-010-0064-5},
	volume = {4},
	year = {2010}}

@article{SoftTrimming,
	author = {S. Taheri and A. M. Bagirov and N. Sultanova and B. Ordin},
	doi = {https://doi.org/10.1016/j.patrec.2024.06.032},
	issn = {0167-8655},
	journal = {Pattern Recognit. Lett.},
	pages = {15--22},
	title = {Robust clustering algorithm: The use of soft trimming approach},
	url = {https://www.sciencedirect.com/science/article/pii/S0167865524002022},
	volume = {185},
	year = {2024}}

@article{GallegosRitter,
	author = {M. T. Gallegos and G. Ritter},
	doi = {10.1214/009053604000000940},
	journal = {Ann. Statist.},
	number = {1},
	pages = {347--380},
	publisher = {Institute of Mathematical Statistics},
	title = {{A robust method for cluster analysis}},
	url = {https://doi.org/10.1214/009053604000000940},
	volume = {33},
	year = {2005}}

@article{RobustBregman,
	ISSN = {00905364, 21688966},
	URL = {https://www.jstor.org/stable/27169435},
	author = {C. Brécheteau and A. Fischer and C. Levrard},
	journal = {Ann. Statist.},
	number = {3},
	pages = {1679--1701},
	publisher = {Institute of Mathematical Statistics},
	title = {ROBUST {B}REGMAN CLUSTERING},
	urldate = {2025-06-22},
	volume = {49},
	year = {2021}
}

@article{WassersteinTrimmed,
	author = {del Barrio, E. and Cuesta-Albertos, J. A. and Matr{\'a}n, C. and Mayo-{\'I}scar, A.},
	date = {2019/01/01},
	doi = {10.1007/s11222-018-9800-z},
	id = {del Barrio2019},
	isbn = {1573-1375},
	journal = {Stat. Comput.},
	number = {1},
	pages = {139--160},
	title = {Robust clustering tools based on optimal transportation},
	url = {https://doi.org/10.1007/s11222-018-9800-z},
	volume = {29},
	year = {2019}}

@article{TrimmedCLT,
	ISSN = {00905364, 21688966},
	URL = {http://www.jstor.org/stable/120150},
	author = {L. A. García-Escudero and A. Gordaliza and C. Matrán},
	journal = {Ann. Statist.},
	number = {3},
	pages = {1061--1079},
	publisher = {Institute of Mathematical Statistics},
	title = {A Central Limit Theorem for Multivariate Generalized Trimmed k-Means},
	urldate = {2025-06-22},
	volume = {27},
	year = {1999}
}

@article{LloydHighDimI,
	title={An Observation on {L}loyd’s $k$-Means Algorithm in High Dimensions}, 
	author={D. Silva-S{\'a}nchez and R. R. Lederman},	
	year={2025},
	journal={arXiv preprint 2506.14952},
}

@article{LloydHighDimII,
	title={The Catastrophic Failure of The k-Means Algorithm in High Dimensions, and How {H}artigan's Algorithm Avoids It}, 
	author={Lederman, R. R. and Silva-S{\'a}nchez, D. and Chen, Z. and Mordant, G. and Balanov, A. and Bendory, T.},	
	year={2026},
	journal={arXiv preprint 2602.09936},
}

@article{WangSong,
  author = {Wang, M. and Song, H.},
  title = {Ckmeans.1d.dp: Optimal k-means Clustering in One Dimension by Dynamic Programming},
  journal = {RJ},
  year = {2011},
  volume = {3},
  issue = {2},
  issn = {2073-4859},
  pages = {29-33}
}

@article{Fisher,
	author = {W. D. Fisher},
	eprint = {https://www.tandfonline.com/doi/pdf/10.1080/01621459.1958.10501479},
	journal = {J. Am. Stat. Assoc.},
	number = {284},
	pages = {789--798},
	publisher = {Taylor \& Francis},
	title = {On Grouping for Maximum Homogeneity},
	url = {https://www.tandfonline.com/doi/abs/10.1080/01621459.1958.10501479},
	volume = {53},
	year = {1958}}
\bibliographystyle{plain}

\end{document}